\documentclass[11pt, a4paper]{article}

\usepackage{comment}
\usepackage{indentfirst}
\usepackage{amsfonts}
\usepackage{amsmath}
\usepackage{amssymb}
\usepackage{mathtools}
\usepackage{float}
\usepackage{graphicx}
\usepackage[graphicx]{realboxes}
\graphicspath{{figure/}}
\usepackage{amsthm}
\usepackage{color}
\usepackage{tikz}
\usepackage{xcolor}
\usetikzlibrary{positioning, shapes.geometric}
\usepackage{makecell}
\usepackage{siunitx}
\usepackage{multirow}
\usepackage[reftex]{theoremref}
\usepackage{bm}
\usepackage{wasysym}
\usepackage{leftindex}
\usepackage{rotating}
\usepackage{authblk}

\newtheorem{theorem}{Theorem}[section]
\newtheorem{lemma}[theorem]{Lemma}

\newtheorem{proposition}[theorem]{Proposition}
\newtheorem{definition}[theorem]{Definition}

\newtheorem{example}[theorem]{Example}
\newtheorem{remark}[theorem]{Remark}

\newcommand{\grobner}{Gr\"{o}bner }
\newcommand{\p}[1]{\bm{#1}}
\newcommand{\bases}[1]{\langle #1 \rangle}

\newcommand{\fnum}{\mathbb{F}}
\newcommand{\pset}{\mathcal}

\newcommand{\lrhombus}{\turnbox{58}{$\lozenge$}}
\newcommand{\rrhombus}{\turnbox{122}{$\lozenge$}}
\DeclareRobustCommand{\abinom}{\genfrac{\{}{\}}{0pt}{}}

\DeclareMathOperator{\zero}{S}
\DeclareMathOperator{\var}{V}
\DeclareMathOperator{\gr}{Gr}
\DeclareMathOperator{\wt}{wt}

\definecolor{lightgreen}{RGB}{199,233,192}
\definecolor{lightblue}{RGB}{173,216,230}
\definecolor{lightred}{RGB}{255,182,193}
\definecolor{dlgreen}{RGB}{70,170,70}

\newcommand{\drawuptriangle}[6]{
    
    \filldraw[fill=#1] (#2, #3) -- ++(0.5,0.866) -- ++(0.5,-0.866) -- cycle;
    \node at (#2 + 0.14,#3 + 0.5) {#4};
    \node at (#2 + 0.86,#3 + 0.5) {#5};
    \node at (#2 + 0.5,#3 + 0.15) {#6};
}

\newcommand{\drawdowntriangle}[6]{
    
    \filldraw[fill=#1] (#2, #3) -- ++(0.5,0.866) -- ++(-1,0) -- cycle;
    \node at (#2,#3 + 0.966+0.05) {#4};
    \node at (#2-0.14,#3 + 0.5) {#5};
    \node at (#2 + 0.14,#3 + 0.5) {#6};
}

\newcommand{\drawuprumbus}[3]{
    
    \filldraw[fill=#1] (#2, #3) -- ++(0.5,0.866) -- ++(-0.5,0.866) -- ++ (-0.5,-0.866) -- cycle;
    \node at (#2-0.14,#3 + 0.5) {0};
    \node at (#2 + 0.14,#3 + 0.5) {1};
    \node at (#2 + 0.36,#3 + 0.5+0.866) {0};
    \node at (#2-0.36,#3 + 0.5+0.866) {1};
}

\newcommand{\drawuprumbusReverse}[3]{
    
    \filldraw[fill=#1] (#2, #3) -- ++(0.5,0.866) -- ++(-0.5,0.866) -- ++ (-0.5,-0.866) -- cycle;
    \node at (#2-0.14,#3 + 0.5) {1};
    \node at (#2 + 0.14,#3 + 0.5) {0};
    \node at (#2 + 0.36,#3 + 0.5+0.866) {1};
    \node at (#2-0.36,#3 + 0.5+0.866) {0};
}

\newcommand{\drawNErumbus}[3]{
    
    \filldraw[fill=#1] (#2, #3) -- ++(1,0) -- ++(-0.5,0.866) -- ++ (-1,0) -- cycle;
    \node at (#2 + 0.5 ,#3 + 0.15) {0};
    \node at (#2 + 0.84,#3 + 0.5) {1};
    \node at (#2 ,#3 +0.966+0.05) {0};
    \node at (#2-0.14,#3 + 0.5) {1};
}
\newcommand{\drawNErumbusReverse}[3]{
    
    \filldraw[fill=#1] (#2, #3) -- ++(1,0) -- ++(-0.5,0.866) -- ++ (-1,0) -- cycle;
    \node at (#2 + 0.5 ,#3 + 0.15) {1};
    \node at (#2 + 0.84,#3 + 0.5) {0};
    \node at (#2 ,#3 +0.966+0.05) {1};
    \node at (#2-0.14,#3 + 0.5) {0};
}

\newcommand{\drawSWrumbus}[3]{
    
    \filldraw[fill=#1] (#2, #3) -- ++(1,0) -- ++(0.5,0.866) -- ++ (-1,0) -- cycle;
    \node at (#2 + 0.5 ,#3 + 0.15) {1};
    \node at (#2 + 1.14,#3 + 0.5) {0};
    \node at (#2 + 1,#3 +0.966+0.05) {1};
    \node at (#2 +0.14,#3 + 0.5) {0};
}



\begin{document}
\sloppy


  \title{Puzzle Ideals for Grassmannians}

  \author{Chenqi Mou}

  \author{Weifeng Shang}
  \affil{LMIB -- School of Mathematical Sciences, \\
  Beihang University, Beijing 100191, China \\
  \{chenqi.mou, weifengshang\}@buaa.edu.cn}

\date{}
\maketitle

  \begin{abstract}
Puzzles are a versatile combinatorial tool to interpret the Littlewood-Richardson coefficients for Grassmannians. In this paper, we propose the concept of puzzle ideals whose varieties one-one correspond to the tilings of puzzles and present an algebraic framework to construct the puzzle ideals which works with the Knutson-Tao-Woodward puzzle and its $T$-equivariant and $K$-theoretic variants for Grassmannians. For puzzles for which one side is free, we propose the side-free puzzle ideals whose varieties one-one correspond to the tilings of side-free puzzles, and the elimination ideals of the side-free puzzle ideals contain all the information of the structure constants for Grassmannians with respect to the free side. 

Besides the underlying algebraic importance of the introduction of these puzzle ideals is the computational feasibility to find all the tilings of the puzzles for Grassmannians by solving the defining polynomial systems, demonstrated with illustrative puzzles via computation of \grobner bases. 

  \end{abstract}


\section{Introduction}
The Grassmannian $\gr_k(\mathbb{C}^n)$ is the set of all the $k$-dimensional linear subspaces of $\mathbb{C}^n$, and the cohomology ring $H^*(\gr_k(\mathbb{C}^n))$ of the Grassmannian is one underlying algebraic structure for counting intersections of projective linear spaces in the intersection theory of algebraic geometry or in Schubert calculus \cite{Ful13i,Man01s}. The Schubert classes $\{[X_{\lambda}]: \lambda \mbox{ a partition}\}$ form a $\mathbb{Z}$-linear basis of $H^*(\gr_k(\mathbb{C}^n))$, and thus the product of two Schubert classes can be expanded as $[X_{\lambda}] [X_{\mu}] = \sum_{\nu} c_{\lambda \mu }^{\nu} [X_{\nu}]$, where $\lambda$, $\mu$, and $\nu$ are partitions. The coefficients $c_{\lambda \mu }^{\nu}$ here are called the Littlewood-Richardson coefficients, and the famous Littlewood-Richardson rule based on Young tableaux, among other combinatorial rules, can fully describe these structure constants \cite{FUL97Y}. The standard representations of the Schubert classes for $\gr_k(\mathbb{C}^n)$ are the Schur polynomials defined by semi-standard Young tableaux and thus Littlewood-Richardson rule can also be viewed as with respect to (w.r.t. hereafter) Schur polynomials. We would like to mention that for the cohomology ring of the complete flag variety there is a similar formula $[X_{\mu}] [X_{\nu}] = \sum_{\omega \in S_n} C_{\mu \nu}^{\omega} [X_{\omega}]$ for all the Schubert classes, indexed by permutations in $S_n$, which are also a $\mathbb{Z}$-linear basis of the cohomology ring. The coefficients $C_{\mu \nu}^{\omega}$ here, which are called the Schubert structure constants and have been proved to be non-negative integers \cite{Kle74t}, are one central object in the study of Schubert calculus.

In this paper we are interested in puzzles, another versatile combinatorial tool for interpreting the Littlewood-Richardson coefficients and more. Puzzles are tilings of an equilateral triangle labeled by $\lambda$, $\mu$, and $\nu$ in the Littlewood-Richardson coefficients with a set of specific puzzle pieces. Puzzles are first introduced in the seminal paper of Knutson, Tao, and Woodward \cite{K04h} and it is proved there that the number of puzzles equals the Littlewood-Richardson coefficient by relating puzzles to the honeycombs for solving the Horn conjecture \cite{KT99t}. Then a combinatorial proof for the equality is given in \cite{K03p} by considering a new kind of puzzles for the $T$-equivariant cohomology. In his geometric interpretation of the Littlewood-Richardson rule \cite{V06g}, Vakil also studies puzzles for the $K$-theory of Grassmannians, exploiting the structure constants in the multiplication of $K$-theoretic classes, or equivalently of stable Grothendieck polynomials \cite{FK94g}. A similar puzzle is also proposed by Wheeler and Zinn-Justin for the $K$-theoretic Littlewood-Richardson rule in \cite{W19l}. Then the structure constants in multiplication of the Schubert bases in the $K$-homology ring of Grasssmannians are studied with another kind of puzzles by Pylyavskyy and Yang in \cite{P20p}. The interpretation of structure constants for Grassmannians in different settings above by means of puzzles gives manifestly positive formulae for the corresponding structure constants. 

The Schubert structure constants in the cohomology ring of complete flag varieties are central in Schubert calculus, and the power of puzzles, already fully demonstrated for Grassmannians, also extends to complete flag varieties: the attempts to exploit the Schubert structure constants with puzzles are also successful for special cases like when the indexing permutations have (almost) separate descents in \cite{KZJ23s}. It is worth mentioning here that yet another tool to study the Schubert calculus of the complete flag variety is the (bumpless) pipe dream \cite{LLS21b,FGX23b,HP22b}. 

In this paper we focus on the puzzles for Grassmannians in different settings and propose a general algebraic framework to construct polynomial ideals from the puzzles whose varieties one-one correspond to all the tilings of the puzzles. To construct these \emph{puzzle ideals}, we first refine the $\fnum_2$-valued puzzle pieces into unit triangle puzzle pieces with $\fnum_3$-values by introducing sides of the additional $2$-value. In this way a refined puzzle with the new set of $\fnum_3$-valued puzzle pieces is constructed from the original puzzle, and we define the \emph{atomic puzzle ideal} whose variety one-one corresponds to the tilings of the refined puzzle. The next step is to restrict the tilings of the refined puzzle so that they can be recovered to those of the original one. For this we introduce an intermediate tiling obtained by stitching $\fnum_3$-valued puzzle pieces with shared $2$-sides in the tiling of the refined puzzle, and study the differences between this intermediate tilling against the original tiling, which lie in the appearances of \emph{forbidden} puzzle pieces which are not original pieces and \emph{implicit} puzzle pieces which are hidden in the original pieces. By adding corresponding polynomials for these puzzle pieces, we extend the atomic puzzle ideal for the refined puzzle to the \emph{puzzle ideal} for the original puzzle and prove the one-one correspondence between the variety of the puzzle ideal and the tilings of the original puzzle. 

Besides the underlying algebraic importance of the introduction of puzzle ideals is the computational feasibility now to find all the tilings of the puzzles for Grassmannians by solving the defining polynomial systems of the puzzle ideals, a topic extensively studied in symbolic computation \cite{CLO1997I} and for which one highly efficient method is based on the computation of \grobner bases \cite{B85G}. For the side-free puzzles for which the values of one side of the equilateral triangle are not determined a prior, their tilings cover all the structure constants. Very naturally from puzzle ideals, we introduce the \emph{side-free puzzle ideals} whose varieties one-one correspond to the tilings of side-free puzzles. Then one can reveal all the information of the structure constants w.r.t. the undetermined side by considering the elimination ideals of the side-free puzzle ideals, which again can be efficiently computed by using the elimination property of lexicographic \grobner bases. We show that the prime decomposition of this radical elimination ideal corresponds to all the non-zero structure constants and that the puzzle ideals w.r.t. each value assignment of the undetermined side can be constructed by taking the sum of the side-free puzzle ideal and the corresponding prime ideal in the decomposition. 

The underlying ideas of refining the original puzzles for atomic puzzle ideals and how to deal with the forbidden and implicit puzzle pieces in the stitching are demonstrated in Sections~\ref{sec:kt}, \ref{sec:t}, and \ref{sec:k} respectively in our treatments of Knutson-Tao-Woodward, $T$-equivariant, and $K$-theoretic puzzles. Then we formulate the framework for defining the puzzle ideals for general puzzles with \emph{separable} puzzle pieces and prove the one-one correspondence in Section~\ref{sec:puzzle}. The side-free puzzle ideals and their properties are discussed in Section~\ref{sec:side-free}.

\section{Preliminaries}
\label{sec:pre}

\subsection{Schur polynomial and Littlewood-Richardson coefficient}

Schur polynomials are an important kind of symmetric polynomials defined with semi-standard Young tableaux. A \emph{Young diagram} is a finite collection of boxes arranged in left-justified rows such that the number of boxes in the rows weakly decrease. Young diagrams can be identified with sequences of weakly decreasing non-negative integers, which are called \emph{partitions}. Let $ \lambda$ be a partition. Then a \emph{semi-standard Young tableau of shape $\lambda$} is a filling of the Young diagram identified with $\lambda$ with positive integers such that in each row the integers weakly increase and in each column the integers strictly increase. For each integer in a semi-standard Young tableau $T$, the number of its occurrences is called its \emph{weight} in $T$. 

Let $ \lambda = (\lambda_{1},\lambda_{2},\ldots,\lambda_{n})$ be a partition. Then we denote $|\lambda| := \sum_{i=1}^n \lambda_i$ or equivalently the number of boxes in the Young diagram corresponding to $\lambda$. To an arbitrary semi-standard Young tableau $T$ of shape $\lambda$ we can associate a term $\p{x}^T$ in $n$ variables $x_1, \ldots, x_n$: the variable $x_i$ corresponds to the integer $i$ in $T$, and the term is defined as $ \p{x}^T = x_{1}^{t_{1}}\cdots x_{n}^{t_{n}}$, where $t_i$ is the weight of $i$ in $T$. Then the \emph{Schur polynomial} w.r.t. $ \lambda $ is $ s_{\lambda}(x_1, \ldots, x_n) = \sum_{T} \p{x}^{T} $, where the summation is taken over all the semi-standard Young tableaux $T$ of shape $\lambda$.

\begin{example}\rm
  For the simple partition $(3, 1)$, all the three semi-standard Young tableaux of shape $(3, 1)$ are listed in Figure~\ref{fig:schur}. Then the Schur polynomial is $S_{(3,1)} = x_{1}^{3} x_{2} + x_{1}^{2} x_{2}^{2} +x_{1}x_{2}^{3} $.

\begin{figure}[H]
\begin{center}

\begin{tikzpicture}[scale=0.75]
        \draw (0,0) -- (0,2);
        \draw (1,0) -- (1,2);
        \draw (0,0) -- (1,0);
        \draw (0,1) -- (3,1);
        \draw (0,2) -- (3,2);
        \draw (2,2) -- (2,1);
        \draw (3,2) -- (3,1);
        
        \node at (0.5,0.5) {2};
        \node at (0.5,1.5) {1};
        \node at (1.5,1.5) {1};
        \node at (2.5,1.5) {1};
        
        \draw (0 + 4,0) -- (0 + 4,2);
        \draw (1 + 4,0) -- (1 + 4,2);
        \draw (0+ 4,0) -- (1+ 4,0);
        \draw (0+ 4,1) -- (3+ 4,1);
        \draw (0+ 4,2) -- (3+ 4,2);
        \draw (2+ 4,2) -- (2+ 4,1);
        \draw (3+ 4,2) -- (3+ 4,1);

        \node at (4.5,0.5) {2};
        \node at (4.5,1.5) {1};
        \node at (5.5,1.5) {1};
        \node at (6.5,1.5) {2};

        \draw (0+8,0) -- (0+8,2);
        \draw (1+8,0) -- (1+8,2);
        \draw (0+8,0) -- (1+8,0);
        \draw (0+8,1) -- (3+8,1);
        \draw (0+8,2) -- (3+8,2);
        \draw (2+8,2) -- (2+8,1);
        \draw (3+8,2) -- (3+8,1);

        \node at (8.5,0.5) {2};
        \node at (8.5,1.5) {1};
        \node at (9.5,1.5) {2};
        \node at (10.5,1.5) {2};

\end{tikzpicture}

\caption{All the semi-standard Young tableaux of shape $(3, 1)$}\label{fig:schur}
    
\end{center}
\end{figure}
  
\end{example}

Schur polynomials are symmetric polynomial and they form a $\mathbb{Z}$-linear basis of the space of all symmetric polynomials. Then the product of any two Schur polynomials $s_{\lambda} s_{\mu}$ can be written uniquely as a linear combination of Schur polynomials in the form $ s_{\lambda} s_{\mu} = \sum_{\nu} c_{\lambda \mu}^{\nu} s_{\nu}$. The coefficient $c_{\lambda \mu}^{\nu} $ in the expression above is called the \emph{Littlewood-Richardson coefficient}. It can be shown that $c_{\lambda \mu}^{\nu} \neq 0$ if and only if $|\nu| = |\lambda| + |\mu|$. Then the sum in $\sum_{\nu} c_{\lambda \mu}^{\nu} s_{\nu}$ can be taken for $|\nu| = |\lambda| + |\mu|$, and thus it is a finite sum. These coefficients are also the structure constants for the cohomology ring $H^*(\gr_k(\mathbb{C}^n))$ of the Grassmannian $\gr_k(\mathbb{C}^n)$. With Schubert classes indexed by partitions as the $\mathbb{Z}$-linear basis of this cohomology ring, the product of two Schubert classes can be written as $[X_{\lambda}] [X_{\mu}] = \sum_{\nu} c_{\lambda \mu }^{\nu} [X_{\nu}]$.

The Littlewood-Richardson rule is a rule to determine these coefficients based on skew Young tableaux \cite{FUL97Y}. In this paper, we consider the puzzles for interpreting the Littlewood-Richardson coefficients, and here Schur polynomials are indexed by binary sequences defined in the following way.

Let $ n \ge k > 0 $. Consider a rectangle of $ k \times (n-k) $ boxes and a partition $\lambda$ whose corresponding Young diagram fits inside this rectangle. We locate the Young diagram identified with $\lambda$ at the top left corner of the rectangle. Then we construct a binary sequence for $\lambda$ by moving from the top right corner to the bottom left one of the rectangle along the boundaries of the rectangle and the Young diagram: for each horizontal edge we write a $0$, and for each vertical edge a $1$. At the endpoint, we will have a binary sequence of length $ n $ with $ k $ $1$s and $ n-k $ $0$s, and this is the binary sequence corresponding to $\lambda$. We use $\abinom{n}{k}$ to denote the set of such binary sequences. Clearly when $n$ and $k$ are fixed, there is a one-one correspondence between the Young diagrams inside the $ k \times (n-k) $ rectangle and $\abinom{n}{k}$. 

\begin{example}\rm
Let us continue with the partition $(3, 1)$. Its  corresponding binary sequence in $\abinom{6}{3}$ and how it is constructed are shown in Figure~\ref{fig:binary} below. 

  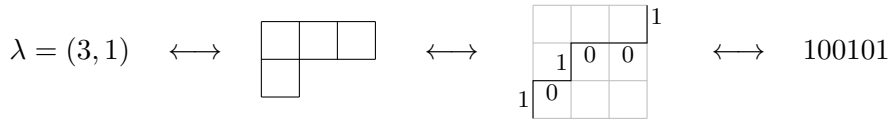
\begin{figure}[H]
  \begin{center}

$\lambda = (3,1) \quad \longleftrightarrow \quad $ \raisebox{-0.5\height}{\begin{tikzpicture}[scale=0.5]
        \draw (0,0) -- (0,2);
        \draw (1,0) -- (1,2);
        \draw (2,1) -- (2,2);
        \draw (3,1) -- (3,2);

        \draw (0,0) -- (1,0);
        \draw (0,1) -- (3,1);
        \draw (0,2) -- (3,2);
        
    \end{tikzpicture} } $\quad \longleftrightarrow  \quad $\raisebox{-0.5\height}{\begin{tikzpicture}[scale=0.5]
        \foreach \x in {0,...,3} 
            \draw[color=gray!50] (\x,0) -- (\x,3);
        \foreach \y in {0,...,3} 
            \draw[color=gray!50] (0,\y) -- (3,\y);
        \draw (0,0) -- (0,1) -- (1,1) -- (1,2) -- (2,2) -- (3,2) -- (3,3);
        \node at (-0.25,0.5) {\footnotesize 1};
        \node at (0.5,0.7) {\footnotesize 0};
        \node at (0.75,1.5) {\footnotesize 1};
        \node at (1.5,1.7) {\footnotesize 0};
        \node at (2.5,1.7) {\footnotesize 0};
        \node at (3.25,2.7) {\footnotesize 1};

    \end{tikzpicture}} $\quad \longleftrightarrow \quad 100101 $
    \end{center}
    \caption{The binary sequence for $(3, 1)$}\label{fig:binary}
  \end{figure}
  
\end{example}

\subsection{Knutson-Tao-Woodward puzzle for the Littlewood-Richardson rule}
In their seminal paper \cite{K04h}, Knutson, Tao, and Woodward introduce puzzles and prove that the numbers of tilings of puzzles equal the Littlewood-Richardson coefficients. Next we first define puzzles in a more general setting. 

\begin{figure}
\begin{center}
\begin{tikzpicture}[scale = 0.8]

    \draw (0,0) -- (3,3*2*0.866) -- (6,0) -- (0,0);
    \node at (0.14,0.5) {0};
    \node at (0.5+0.14,0.866+0.5) {1};
    \node at (0.5*2+0.14,0.866*2+0.5) {0}; 
    \node at (0.5*3+0.14,0.866*3+0.5) {1}; 
    \node at (0.5*4+0.14,0.866*4+0.5) {0}; 
    \node at (0.5*5+0.14,0.866*5+0.5) {1}; 

    \node at (6-0.14,0.5) {1};
    \node at (6-0.14-0.5,0.5+0.866) {0};
    \node at (6-0.14-0.5*2,0.5+0.866*2) {1};
    \node at (6-0.14-0.5*3,0.5+0.866*3) {0};
    \node at (6-0.14-0.5*4,0.5+0.866*4) {1};
    \node at (6-0.14-0.5*5,0.5+0.866*5) {0};

    \node at (0.5,0.1) {1};
    \node at (0.5+1,0.1) {0};
    \node at (0.5+2,0.1) {1};
    \node at (0.5+3,0.1) {1};
    \node at (0.5+4,0.1) {0};
    \node at (0.5+5,0.1) {1};

    \draw [->] (-0.5,0.5) -- (2,0.5+0.866*5); 
    \draw [->] (4,0.5+0.866*5) -- (6.5,0.5); 
    \draw [->] (0.5,-0.5) -- (5.5,-0.5);
    
\end{tikzpicture}    
\end{center}
\caption{An illustrative $\bigtriangleup^{\nu}_{\lambda \mu}$}\label{fig:0-bigtriangle} 
\end{figure}
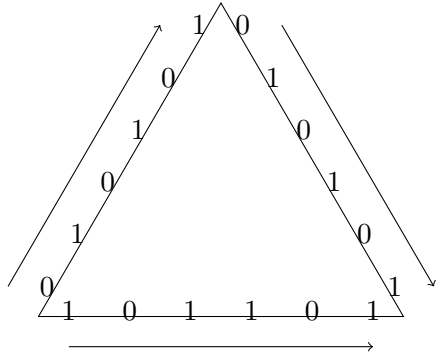

Let $\lambda, \mu, \nu \in \abinom{n}{k}$. Then $\bigtriangleup^{\nu}_{\lambda \mu}$ is the upward equilateral triangle with side-length of $n$ units such that the $n$ units of its left, right, and bottom sides are assigned to the $\fnum_2$-value of $\lambda$, $\mu$, and $\nu$ from left to right respectively (see Figure~\ref{fig:0-bigtriangle} below for an illustration). A ($\fnum_q$-valued) \emph{puzzle piece} is a convex polygon such that each of its sides is parallel to one side of $\bigtriangleup^{\nu}_{\lambda \mu}$ and of length of multiple units and each unit of the side is assigned to an $\fnum_q$-value. We will work with $\fnum_3$-valued puzzle pieces in this paper, and this justifies our definition with $\fnum_q$-valued pieces. 

Let $p$ be an $\fnum_q$-valued puzzle piece (possibly $\bigtriangleup^{\nu}_{\lambda \mu}$ itself after embedding $\fnum_2$ to $\fnum_q$) and $\Omega$ be a set of $\fnum_q$-valued puzzle pieces. Then by a \emph{tiling} of $p$ with $\Omega$ we mean a covering of $p$ using the puzzle pieces in $\Omega$ (possibly multiple times) with no overlaps and no gaps such that the values of the common edges of any two adjacent pieces match. Note that rotating or mirroring a puzzle piece in general results in a different piece, and thus no rotation and mirroring is allowed in the tiling. 

A \emph{puzzle} consists of the $\fnum_2$-valued $\bigtriangleup^{\nu}_{\lambda \mu}$ and a set $\Omega$ of $\fnum_2$-valued puzzle pieces, and we denote this puzzle by $P^{\nu, \Omega}_{\lambda \mu}$. For a puzzle $P^{\nu, \Omega}_{\lambda \mu}$, its \emph{solution} is a tiling of $\bigtriangleup^{\nu}_{\lambda \mu}$ with $\Omega$, and we denote the set of all its solutions by $\zero(P^{\nu, \Omega}_{\lambda \mu})$. Note that in the literature, a puzzle usually refers to a tiling, which is different from our notions. We prefer to view the puzzle as the game we want to solve, with the tilings as its solutions. This is consistent with our treatment with puzzle ideals for puzzles and the varieties of puzzle ideals (or solutions of the defining polynomial systems) for the solutions of puzzles in this paper.

The puzzle pieces in the Knutson-Tao-Woodward puzzle are those shown in Figure~\ref{fig:0-piece} below, and we denote the set of such puzzle pieces by $\Omega_0$ throughout this paper. 

\begin{figure}[H]
    \begin{center}
        \begin{tikzpicture}
            \drawuptriangle{lightred}{0}{0}{0}{0}{0};
            \drawdowntriangle{lightred}{2}{0}{0}{0}{0};
            \drawuptriangle{lightblue}{3}{0}{1}{1}{1};
            \drawdowntriangle{lightblue}{5}{0}{1}{1}{1};
            \drawNErumbus{lightgreen}{9}{0};
            \drawuprumbus{lightgreen}{7}{0};
            \drawSWrumbus{lightgreen}{11}{0};
        \end{tikzpicture}
    \end{center}
      \caption{The puzzle pieces for Knutson-Tao-Woodward puzzles}\label{fig:0-piece}
\end{figure}
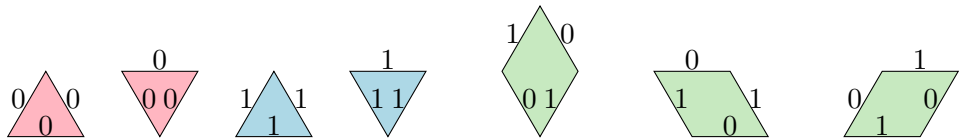

\begin{example}\rm

Consider $\lambda = 0101$, $\mu = 0101$, and $\nu = 0110$ from $\abinom{4}{2}$. One tiling of the Knutson-Tao-Woodward puzzle $P^{\nu, \Omega_0}_{\lambda \mu}$ is illustratived in Figure~\ref{fig:0-solution} below. 

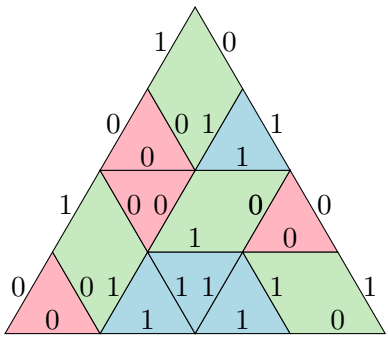
\begin{figure}[H]

\begin{center}

\begin{tikzpicture}[scale=1.25]



\drawuptriangle{lightred}{0}{0}{0}{0}{0};
\drawuptriangle{lightblue}{1}{0}{1}{1}{1};
\drawuptriangle{lightblue}{2}{0}{1}{1}{1};
\drawdowntriangle{lightblue}{2}{0}{1}{1}{1};
\drawuprumbus{lightgreen}{1}{0};
\drawNErumbus{lightgreen}{3}{0};
\drawSWrumbus{lightgreen}{1.5}{0.866};
\drawdowntriangle{lightred}{1.5}{0.866}{0}{0}{0};
\drawuptriangle{lightred}{2.5}{0.866}{0}{0}{0};
\drawuptriangle{lightred}{1}{0.866*2}{0}{0}{0};
\drawuptriangle{lightblue}{2}{0.866*2}{1}{1}{1};
\drawuprumbus{lightgreen}{2}{0.866*2};

\end{tikzpicture}
\end{center}

\caption{One tiling of Knutson-Tao-Woodward puzzle $P^{\nu, \Omega_0}_{\lambda \mu}$}\label{fig:0-solution} 
\end{figure}

\end{example}

The following theorem relates the numbers of solutions of the Knutson-Tao-Woodward puzzle to the Littlewood-Richardson coefficients. 

\begin{theorem}[{\cite[Theorem~1]{K04h}}]
\label{puzzle}
  Let $\lambda$, $\mu$, and $\nu$ be partitions in $\abinom{n}{k}$. Then $\#\zero(P^{\nu, \Omega_0}_{\lambda \mu}) = c_{\lambda \mu}^{\nu}$.
\end{theorem}

\subsection{$T$-equivariant puzzle for equivariant cohomology of Grassmannians}

In the combinatorial proof of Theorem~\ref{puzzle} in \cite{K03p}, the equivariant cohomology ring of Grassmannians is considered and a corresponding puzzle is proposed in the equivariant case. 

The equivariant cohomology ring $H^{*}_{T}(\gr_{k}(\mathbb{C}^{n})) $ is naturally a graded module over the polynomial ring $ \mathbb{Z}[y_{1},\ldots,y_{n}]$ and it has a basis of equivariant Schubert classes $\{[\tilde{X}_{\lambda}]\}$. The equivariant cohomology ring $H^{*}_{T}(\gr_{k}(\mathbb{C}^{n})) $ and the ordinary one $H^{*}(\gr_{k}(\mathbb{C}^{n}))$ are connected by a natural forgetful map $H^{*}_{T}(\gr_{k}(\mathbb{C}^{n})) \rightarrow H^{*}(\gr_{k}(\mathbb{C}^{n}))$ which sets all the $y_{i}$ to $0$. In particular, this forgetful map takes each equivariant Schubert class $[\tilde{X}_{\lambda}]$ to the corresponding counterpart $[X_{\lambda}]$. Then the product of two equivariant Schubert classes can be uniquely expanded as $[\tilde{X}_{\lambda}] [\tilde{X}_{\mu}] = \sum_{\nu} \tilde{c}_{\lambda \mu}^{\nu} [\tilde{X}_{\nu}]$. The coefficients $\tilde{c}_{\lambda \mu}^{\nu} \in \mathbb{Z}[y_{1},\ldots,y_{n}]$ are the equivariant structure constants, and they agree with the Littlewood-Richardson coefficients $c_{\lambda \mu}^{\nu}$ when $ |\nu| = |\lambda| + |\mu| $.

\begin{figure}[H]
\begin{center}
    \begin{tikzpicture}
        \drawuptriangle{lightblue}{0}{0}{1}{1}{1};
        \drawuptriangle{lightblue}{1}{0}{1}{1}{1};
        \drawNErumbus{lightgreen}{2}{0};
        \drawuptriangle{lightblue}{3}{0}{1}{1}{1};
        \drawNErumbus{lightgreen}{4}{0};
        \drawuprumbusReverse{dlgreen}{5}{0};
        \drawuptriangle{lightred}{5}{0}{0}{0}{0};
        \drawdowntriangle{lightblue}{1}{0}{1}{1}{1};
        \drawdowntriangle{lightblue}{3}{0}{1}{1}{1};
        \drawSWrumbus{lightgreen}{0.5}{0.866};
        \drawuptriangle{lightred}{1.5}{0.866}{0}{0}{0};
        \drawSWrumbus{lightgreen}{2.5}{0.866};
        \drawuptriangle{lightred}{3.5}{0.866}{0}{0}{0};
        \drawdowntriangle{lightred}{2.5}{0.866}{0}{0}{0};
        \drawdowntriangle{lightred}{4.5}{0.866}{0}{0}{0};
        \drawSWrumbus{lightgreen}{1}{0.866*2};
        \drawuptriangle{lightred}{2}{0.866*2}{0}{0}{0};
        \drawSWrumbus{lightgreen}{3}{0.866*2};
        \drawuptriangle{lightred}{4}{0.866*2}{0}{0}{0};
        \drawdowntriangle{lightred}{3}{0.866*2}{0}{0}{0};
        \drawuptriangle{lightblue}{1.5}{0.866*3}{1}{1}{1};
        \drawuprumbusReverse{dlgreen}{2.5}{0.866*3};
        \drawuptriangle{lightred}{2.5}{0.866*3}{0}{0}{0};
        \drawuprumbus{lightgreen}{3.5}{0.866*3};
        \drawuptriangle{lightblue}{3.5}{0.866*3}{1}{1}{1};
        \drawdowntriangle{lightblue}{3}{0.866*4}{1}{1}{1};
        \drawuptriangle{lightblue}{2.5}{0.866*5}{1}{1}{1};

        \drawuprumbusReverse{dlgreen}{-2.5}{0.866*2};

        \draw [->,dashed,line width=0.4mm,red] (2.25,2*1.75*0.866) -- (0.5,0);
        \draw [->,dashed,line width=0.4mm,red] (2.75,2*1.75*0.866) -- (4.5,0);
        \draw [->,dashed,line width=0.4mm,red] (4.75,0.5*0.866) -- (4.5,0) ;
        \draw [->,dashed,line width=0.4mm,red] (5.25,0.433) -- (5.5,0); 
    \end{tikzpicture}
\end{center}
\caption{The equivariant puzzle piece (left) and a tiling of a $T$-equivariant puzzle (right)}\label{fig:T-tiling} 
\end{figure}
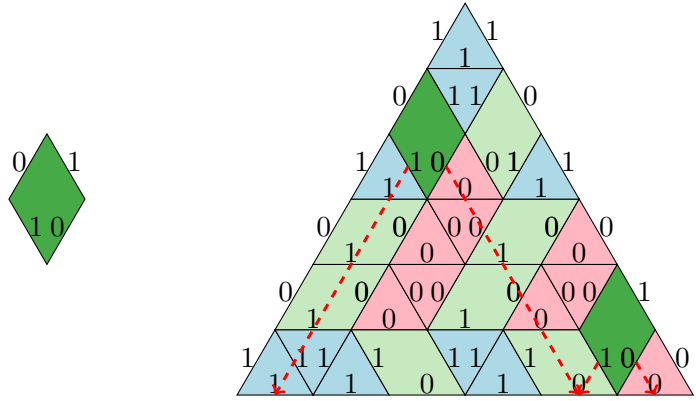

To interpret these equivariant structure constants, the $T$-equivariant puzzle is proposed in \cite{K03p} with an additional rhombus puzzle piece shown in Figure~\ref{fig:T-tiling}~(left). We call this additional piece the \emph{equivariant puzzle piece} and denote by $\Omega_T$ the set obtained by adjoining $\Omega_0$ with it. An illustrative tiling in $\zero(P^{\nu, \Omega_T}_{\lambda \mu})$ is shown in Figure~\ref{fig:T-tiling}~(right). 

For any tiling $t \in \zero(P^{\nu, \Omega_T}_{\lambda \mu})$, to each equivariant puzzle piece $p$ appearing in $t$ we associate a weight $\wt(p)$ in the following way: We drag the piece $p$ in the South-East direction parallel to the right side of $\bigtriangleup^{\nu}_{\lambda \mu}$ until it pokes out the $i$-th interval of the bottom side of $\bigtriangleup^{\nu}_{\lambda \mu}$, and drag it in the South-West direction until it pokes out the $j$-th interval. Then we set the weight $\wt(p) = y_i - y_j$. Taking the product of all the weights of equivariant puzzle pieces appearing in $t$ gives the weight of $t$: $\wt(t) = \prod_{p \in t} \wt(p)$. For the tiling $t$ in Figure~\ref{fig:T-tiling}, one can see that the weights of two equivariant puzzle pieces are $y_5-y_1$ and $y_6-y_5$ respectively, and thus $\wt(t) = (y_5-y_1)(y_6-y_5)$.

The construct structures $\tilde{c}_{\lambda \mu}^{\nu}$ in the equivariant cohomology ring of Grassmannians can be described by the weights of tilings in $\zero(P^{\nu, \Omega_T}_{\lambda \mu})$. 

\begin{theorem}[{\cite[Theorem~2]{K03p}}]
Let $\lambda$, $\mu$, and $\nu$ be partitions in $\abinom{n}{k}$. Then $\sum_{t\in \zero(P^{\nu, \Omega_T}_{\lambda \mu})}\wt(t) = \tilde{c}_{\lambda \mu}^{\nu}$. 
\end{theorem}

\subsection{Puzzles for $K$-theory of Grassmannians}

The $K$-theory of a Grassmannian $\gr_k(\mathbb{C}^n)$ studies its Grothendieck ring $K^{\circ}\gr_k(\mathbb{C}^n)$, and the structure constants in this case are for the multiplication of two classes of the structure sheaves in $K^{\circ}\gr_k(\mathbb{C}^n)$ \cite{B02l}. The representatives of the classes of the structure sheaves in $K^{\circ}\gr_k(\mathbb{C}^n)$ are the stable Grothendieck polynomials \cite{FK94g}. The stable Grothendieck polynomial $ G_{\lambda} $, indexed by a partition $\lambda$, is defined with semi-standard set-valued tableaux of shape $\lambda $ and thus generalization of Schur polynomials. 

It is shown in \cite{B02l}
that the $\mathbb{Z}$-linear span of all the stable Grothendieck polynomials forms a bialgebra with the product given by 
$$ G_\lambda G_\mu = \sum_{\nu} (-1)^{|\nu| - |\lambda| - |\mu|}  \hat{c}^\nu_{\lambda \mu} G_\nu $$
and the coproduct $ \bigtriangleup $ by 
$$ \bigtriangleup(G_\nu) = \sum_{\lambda,\mu} (-1)^{|\nu| - |\lambda| - |\mu|} \hat{d}^\nu_{\lambda\mu} G_\lambda \otimes G_\mu. $$ 
Note that here we explicitly extract the sign $(-1)^{|\nu| - |\lambda| - |\mu|}$ out so that the structure constants $\hat{c}^\nu_{\lambda \mu}$ for the product and $\hat{d}^\nu_{\lambda\mu}$ for the coproduct are non-negative integers. For the product, when $ |\nu| < |\lambda| + |\mu| $, we have $ \hat{c}^\nu_{\lambda \mu} = 0 $ and when $ |\nu| = |\lambda| + |\mu| $, the coefficient $ \hat{c}^\nu_{\lambda \mu} $ equals the  Littlewood-Richardson coefficient $c^\nu_{\lambda \mu}$. For the coproduct, when $ |\nu| > |\lambda| + |\mu|$, we have $ \hat{d}^\nu_{\lambda \mu} = 0 $. 

The following variants of $ \hat{c}^\nu_{\lambda \mu} $ and $ \hat{d}^\nu_{\lambda \mu} $ are studied. Let $ \widetilde{G}_\lambda = G_\lambda \cdot (1 - G_1) $. Then similarly we have 
$$ \widetilde{G}_\lambda \widetilde{G}_\mu = \sum_{\nu} (-1)^{ |\nu| - |\lambda| - |\mu| } \check{c}^\nu_{\lambda \mu} \widetilde{G}_\nu. $$
A corresponding coefficient $\check{d}^\nu_{\lambda \mu}$ can also be defined from $\hat{d}^\nu_{\lambda \mu}$ (see \cite{P20p} for the details).

$K$-theoretic puzzles are proposed in \cite{V06g,W19l,P20p} to describe the structure constants $\hat{c}^\nu_{\lambda \mu} $, $\hat{d}^\nu_{\lambda \mu} $, $\check{c}^\nu_{\lambda \mu} $, and $\check{d}^\nu_{\lambda \mu} $ respectively. These results can be viewed as puzzle interpretation of the $K$-theoretic Littlewood-Richardson coefficients. We refer to the puzzle pieces shown in the figure below as the $A$-, $B$-, $C$-, and $D$-pieces form left to right. The set $\Omega_A$ of puzzle pieces is defined to the union of $\Omega_0$ and the $A$-piece, and $\Omega_B$, $\Omega_C$, and $\Omega_C$ are similarly defined. 

\begin{figure}[H]
  \begin{center}
      \begin{tikzpicture}[scale=0.75]
          \draw (0,0) -- (-1,0.866*2) -- (1,0.866*2) -- cycle;
          \drawdowntriangle{blue}{0}{0.866/2+0.1}{}{}{};

          \draw (2,0) -- (3,0.866*2) -- (4,0) -- cycle;
          \drawuptriangle{orange}{2.5}{0.433-0.1}{}{}{};

          \draw (6,0) -- (6-0.5,0.866) -- (6,0.866*2) -- (7,0.866*2) -- (7.5,0.866) -- (7,0) -- cycle;

          \filldraw[fill=green] (6.25,0.433) -- (6,0.866) -- (6.25,0.866+0.433) -- (7-0.25,0.866+0.433) -- (7,0.866) -- (7-0.25,0.433) -- cycle;
          
          \draw (9,0) -- (9-0.5,0.866) -- (9,0.866*2) -- (10,0.866*2) -- (10.5,0.866) -- (10,0) -- cycle;

          \filldraw[fill=pink] (9.25,0.433) -- (9,0.866) -- (9.25,0.866+0.433) -- (10-0.25,0.866+0.433) -- (10,0.866) -- (10-0.25,0.433) -- cycle;

          \node at (-0.25,0.433) {0};
          \node at (-0.75,0.866+0.433) {1};
          \node at (-0.5, 0.866*2) {0};
          \node at ( 0.5, 0.866*2) {1};
          \node at (0.75, 0.866+0.433) {0};
          \node at (0.25, 0.433) {1};

          \node at (3-0.75,0.433) {0};
          \node at (3-0.25,0.866+0.433) {1};
          \node at (3-0.5, 0) {1};
          \node at (3+0.5, 0) {0};
          \node at (3+0.25, 0.866+0.433) {0};
          \node at (3+0.75, 0.433) {1};

          \node at (6.5,0) {1};
          \node at (5.75,0.433) {0};
          \node at (5.75,0.866+0.433) {1};
          \node at (6.5, 0.866*2) {0};
          \node at (7.25,0.866+0.433) {1};
          \node at (7.25,0.433) {0};

          \node at (9.5,0) {0};
          \node at (8.75,0.433) {1};
          \node at (8.75,0.866+0.433) {0};
          \node at (9.5, 0.866*2) {1};
          \node at (10.25,0.866+0.433) {0};
          \node at (10.25,0.433) {1};
          
      \end{tikzpicture}
  \end{center}
  \caption{Four additional puzzle pieces for $K$-theory of Grassmannians} 
  \label{fig.7}
\end{figure}
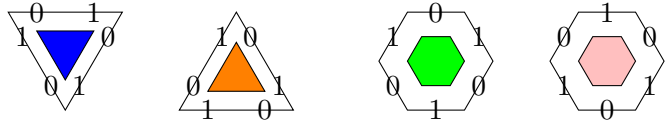

Then the puzzle interpretation of the $K$-theoretic Littlewood-Richardson coefficients are the following theorems for the structure constants in the $K$-theory of Grassmannians. 

\begin{theorem}[\cite{V06g,W19l}]
Let $\lambda,\mu,\nu$ be partitions in $\abinom{n}{k}$ with $|\nu| \ge |\lambda| + |\mu| $. Then $\#\zero(P^{\nu, \Omega_A}_{\lambda \mu}) = \hat{c}^{\nu}_{\lambda \mu}$ and $\#\zero(P^{\nu, \Omega_B}_{\lambda \mu}) = \check{c}^{\nu}_{\lambda \mu}$. 
\end{theorem}

\begin{theorem}[\cite{P20p}]
Let $\lambda,\mu,\nu$ be partitions in $\abinom{n}{k}$ with $|\nu| \le |\lambda| + |\mu| $. Then $\#\zero(P^{\nu, \Omega_C}_{\lambda \mu}) = \hat{d}^{\nu}_{\lambda \mu}$ and $\#\zero(P^{\nu, \Omega_D}_{\lambda \mu}) = \check{d}^{\nu}_{\lambda \mu}$.
\end{theorem}

\section{Atomic puzzle ideal for Knutson-Tao-Woodward puzzle}
\label{sec:kt}

In this section, we first take the original Knutson-Tao-Woodward puzzles for example to show how to find all their solutions by reducing the problem to polynomial systems solving. 

\begin{figure}[H]
\begin{center}
    \begin{tikzpicture}[scale=0.85]
        \drawuptriangle{lightred}{0}{0}{0}{0}{0};
        \drawuptriangle{lightblue}{1.5}{0}{1}{1}{1};
        \drawuptriangle{lightgreen}{3}{0}{1}{0}{2};
        \drawuptriangle{lightgreen}{4.5}{0}{2}{1}{0};
        \drawuptriangle{lightgreen}{6}{0}{0}{2}{1};

        \drawdowntriangle{lightred}{8}{0}{0}{0}{0};
        \drawdowntriangle{lightblue}{9.5}{0}{1}{1}{1};
        \drawdowntriangle{lightgreen}{11}{0}{2}{0}{1};
        \drawdowntriangle{lightgreen}{12.5}{0}{0}{1}{2};
        \drawdowntriangle{lightgreen}{14}{0}{1}{2}{0};
    \end{tikzpicture}
\end{center}
\caption{Atomic puzzle pieces constructed from the Knutson-Tao-Woodward pieces}\label{fig:kt-refined}
\end{figure}
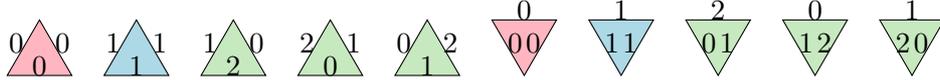

Let $P^{\nu, \Omega_0}_{\lambda \mu}$ be a Knutson-Tao-Woodward puzzle with $\lambda, \mu, \nu \in \abinom{n}{k}$, where $\Omega_0$ consists of the $\fnum_2$-valued puzzle pieces as shown in Figure~\ref{fig:0-piece}. There are two kinds of pieces in $\Omega_0$: four unit triangles and three rhombuses. Our key idea to construct a corresponding polynomial system from $P^{\nu, \Omega_0}_{\lambda \mu}$ is to cut each of the rhombus piece into two unit triangles by adding an extra side with the assigned value $2$ in the middle of the piece. In this way, we turn the $\fnum_2$-valued puzzle pieces in $\Omega_0$ to the following $\fnum_3$-valued ones shown in Figure~\ref{fig:kt-refined}, which are all unit triangles.

An $\fnum_3$-valued puzzle piece is said to be \emph{atomic} if it is a unit triangle whose assigned values sum up to $0$ in $\fnum_3$. Then the four triangular puzzle pieces in $\Omega_0$ are all atomic if we embed its $\fnum_2$-values in $\fnum_3$. Denote by $\overline{\Omega}_0$ this new set of $\fnum_3$-valued atomic puzzle pieces, and we call $\overline{\Omega}_0$ the \emph{atomic refinement} of $\Omega_0$. We will formally define it in a more general setting later in Section~\ref{sec:puzzle}. We also embed the $\fnum_2$-valued $\bigtriangleup^{\nu}_{\lambda \mu}$ in $\fnum_3$ and keep the same notation as long as no ambiguity occurs. In this way, the original puzzle $P^{\nu, \Omega_0}_{\lambda \mu}$ is transformed into a new puzzle $P^{\nu, \overline{\Omega}_0}_{\lambda \mu}$ consisting of the $\fnum_3$-valued $\bigtriangleup^{\nu}_{\lambda \mu}$ and the set $\overline{\Omega}_0$ of atomic puzzle pieces. Similarly for $P^{\nu, \overline{\Omega}_0}_{\lambda \mu}$, we can consider its solutions as the tilings of $\bigtriangleup^{\nu}_{\lambda \mu}$ with $\overline{\Omega}_0$, and $\zero(P^{\nu, \overline{\Omega}_0}_{\lambda \mu})$ denotes the set of all its solutions, in the same way as for $P^{\nu, \Omega_0}_{\lambda \mu}$. It turns out that finding $\zero(P^{\nu, \overline{\Omega}_0}_{\lambda \mu})$ is equivalent to solving a polynomial system associated to the puzzle $P^{\nu, \overline{\Omega}_0}_{\lambda \mu}$. 

\begin{remark}\rm
    The new puzzle pieces in $\overline{\Omega}_0$ shown in Figure~\ref{fig:kt-refined} are exactly the same as those in Section~1.1 of \cite{KZJ17s} for the particular case of Grassmannians, with the difference that the additional value $2 \in \fnum_3$ is used to label the unit triangle pieces here while a non-numeric $10$ is labeled instead in \cite{KZJ17s}. Furthermore, we call $\fnum_3$-valued puzzle pieces atomic if the values of their three sides sum up to 0 in $\fnum_3$, while in \cite{KZJ17s} three directions are introduced to 3 groups of combinations of sides with their labels so that the new puzzle pieces there have their directions sum up to $0$ in the sense of vectors. We feel that our treatment with the introduction of $2\in \fnum_3$ and embedding $\fnum_2$-valued puzzle pieces in $\fnum_3$ is more algebraic and thus natural for constructing polynomial systems to define puzzle ideals, and this technique turns out to be adequate for handling existing puzzles for Grassmannians.  
\end{remark}

Let $P^{\nu, \Omega_0}_{\lambda \mu}$ be a Knutson-Tao-Woodward puzzle with $\lambda, \mu, \nu \in \abinom{n}{k}$ and $P^{\nu, \overline{\Omega}_0}_{\lambda \mu}$ be the puzzle induced by the atomic refinement of $\Omega_0$ as stated above. Then inside $\bigtriangleup^{\nu}_{\lambda \mu}$ there are $n^2$ unit triangles with $N = \frac{3n(n+1)}{2}$ unit intervals as their sides. Out of these $N$ unit intervals, $3n$ ones in the left, right, and bottom sides of $\bigtriangleup^{\nu}_{\lambda \mu}$ have assigned values determined by $\lambda$, $\mu$, and $\nu$. At this step it is not hard to see that a tiling of $\bigtriangleup^{\nu}_{\lambda \mu}$ with $\overline{\Omega}_0$ is equivalent to assigning values to the $N$ unit intervals such that the following conditions hold:
\begin{enumerate}
    \item ($\fnum_3$-valued) The assigned values are all in $\fnum_3$.
    \item (matching $\lambda, \mu, \nu$) The values assigned to the unit intervals in the left, right, and bottom sides of $\bigtriangleup^{\nu}_{\lambda \mu}$ match $\lambda$, $\mu$, and $\nu$.
    \item (atomic) For each unit triangle, the sum of the assigned values of its three sides is $0$ in $\fnum_3$. 
    \item (allowed pieces) Each unit triangle with the assigned values is a puzzle piece in $\overline{\Omega}_0$ (for example, a $(2, 2, 2)$-valued piece is not allowed). 
\end{enumerate}

We introduce one variable for the value of each of these unit intervals so in total we have $N$ variables $x_1, \ldots, x_N$ (at this point we are not particularly interested in how the variables correspond to the unit intervals), and consider the polynomial ring $\fnum_3[\p{x}]$, where $\p{x} := \{x_1, \ldots, x_N\}$. Then the constraints listed above on the assigned values can all be formulated as equations of polynomials in $\fnum_3[\p{x}]$ in the following way. 

\begin{enumerate}
    \item For each $i=1, \ldots, N$, $x_i^3 - x_i = 0$. These are the field equations, and there are $N$ of them. 
    \item Let $x_{l1}, \ldots, x_{ln}$ be the $n$ variables corresponding to the $n$ intervals of the left side of $\bigtriangleup^{\nu}_{\lambda \mu}$. Then $x_{lj} - \lambda_j = 0$ for $j=1, \ldots, n$. The same equations should hold between the other $2n$ variables for the right and bottom sides and $\mu$ and $\nu$. In total we have $3n$ equations in this form. 
    \item Let $\bigtriangleup_k~(k=1, \ldots, n^2)$ be a unit triangle inside $\bigtriangleup^{\nu}_{\lambda \mu}$ and $x_{k1}$, $x_{k2}$, and $x_{k3}$ be the variables corresponding to its three sides in the clockwise direction. Then $x_{k1} + x_{k2} + x_{k3} = 0$. There are $n^2$ such equations. 
    \item For each $\bigtriangleup_k$, let $\overline{x}_{k1}$, $\overline{x}_{k2}$, and $\overline{x}_{k3}$ be the assigned values of $x_{k1}$, $x_{k2}$, and $x_{k3}$. Then we define a polynomial function $f_k(x_{k1}, x_{k2}, x_{k3}) \in \fnum_3[x, y, z]$ such that $f_k(\overline{x}_{k1}, \overline{x}_{k2}, \overline{x}_{k3}) = 0$ if the assigned values correspond to a puzzle piece in $\overline{\Omega}_0$ and $f_k(\overline{x}_{k1}, \overline{x}_{k2}, \overline{x}_{k3}) \neq 0$ (for example $=1$) otherwise.  

    Note that the general form of a polynomial $f(x, y, z) \in \fnum_3[x, y, z]$ has 27 undetermined coefficients (modulo $\bases{x^3-x$, $y^3-y,z^3-z}$), and the constraints on $f(x, y, z)$ above can be translated as 
    \begin{equation*}
    \begin{split}
        &f(0,0,0) = f(1,1,1) = f(0,2,1) = f(2,1,0) = f(1,0,2) = 0,\\
        &f(0,1,2)  = f(1,2,0) = f(2,0,1) = f(2,2,2) =1.
    \end{split}        
    \end{equation*}
    These constraints form a system of $9$ linear equations with $27$ undetermined coefficients as the variables, and thus this is an underdetermined linear system. Picking a specific solution of this linear system, we have $ f(x, y, z) = x+ x^{2} + 2y+ y^{2} +z^{2} + xyz + 2x^{2}yz$ as a feasible choice. We call such a polynomial $f(x, y, z)$ a \emph{distinguishing polynomial} w.r.t. $\overline{\Omega}_0$. 
    
    Then for the specific $\bigtriangleup_k~(k=1, \ldots, n^2)$, the polynomial equation to confine the unit triangle to lie in $\overline{\Omega}_0$ is 
    $$x_{k1}+ x_{k1}^{2} + 2x_{k2}+ x_{k2}^{2} +x_{k3}^{2} + x_{k1}x_{k2}x_{k3} + 2x_{k1}^{2}x_{k2}x_{k3} = 0.$$
    There are $n^2$ such equations in total. 

\end{enumerate}

\begin{remark}\rm
    The underlying reason why we can assign one uniform distinguishing polynomial for all the atomic puzzle pieces in $\overline{\Omega}_0$ is that the pieces in $\overline{\Omega}_0$ are closed under rotation (by 120 degrees) and mirroring (upside down).
\end{remark}

Let $\pset{F}_a$ be the set of all the polynomials in left hands of the equations of four group listed above for the puzzle $P^{\nu, \overline{\Omega}_0}_{\lambda \mu}$, and now consider the ideal $\bases{\pset{F}_a} \subseteq \fnum_3[\p{x}]$ and denote it by $\leftindex_a I^{\nu, \overline{\Omega}_0}_{\lambda \mu}$. Then clearly we have one-one correspondence between the solutions in $\zero(P^{\nu, \overline{\Omega}_0}_{\lambda \mu})$ and the points in $\var(\leftindex_a I^{\nu, \overline{\Omega}_0}_{\lambda \mu}) \subseteq \fnum_3^N$, where $\var(I^{\nu, \overline{\Omega}_0}_{\lambda \mu})$ is the $\fnum_3$-variety of $\leftindex_a I^{\nu, \overline{\Omega}_0}_{\lambda \mu}$. We call the ideal $\leftindex_a I^{\nu, \overline{\Omega}_0}_{\lambda \mu}$ the \emph{atomic puzzle ideal} of the puzzle $P^{\nu, \Omega_0}_{\lambda \mu}$. Finding $\var(\leftindex_a I^{\nu, \overline{\Omega}_0}_{\lambda \mu})$, or equivalently solving the polynomial system $\pset{F}_a = 0$, can be efficiently done by computing a lexicographic \grobner basis of the ideal $\leftindex_a I^{\nu, \overline{\Omega}_0}_{\lambda \mu}$. 

In particular, for any tiling in $\zero(P^{\nu, \Omega_0}_{\lambda \mu})$, if we further tile any puzzle piece in $\Omega_0$ occurring in this tiling with its atomic refinement in $\overline{\Omega}_0$, we will obtain a tiling in $\zero(P^{\nu, \overline{\Omega}_0}_{\lambda \mu})$. This implies that 
$$\#\zero(P^{\nu, \Omega_0}_{\lambda \mu}) \leq \#\zero(P^{\nu, \overline{\Omega}_0}_{\lambda \mu}) = \#\var(\leftindex_a I^{\nu, \overline{\Omega}_0}_{\lambda \mu}).$$
This inequality holds in general between any puzzle and the new puzzle induced by an atomic refinement of its puzzle pieces if such a refinement exists (see Proposition~\ref{prop:refined-more} in Section~\ref{sec:puzzle}). But for the Knutson-Tao-Woodward puzzles, the inequality can be proved to be an equality. To show this, we need to study the reverse process of atomic refinement to turn a tiling in $\zero(P^{\nu, \overline{\Omega}_0}_{\lambda \mu})$ to one in $\zero(P^{\nu, \Omega_0}_{\lambda \mu})$.  

If we compare a tiling in $\zero(P^{\nu, \overline{\Omega}_0}_{\lambda \mu})$ with one in $\zero(P^{\nu, \Omega_0}_{\lambda \mu})$, it is easy to find that the atomic puzzle pieces in $\overline{\Omega}_0 \setminus \Omega_0$ containing $2$-sides make the underlying differences. To turn a tiling $\overline{t}$ in $\zero(P^{\nu, \overline{\Omega}_0}_{\lambda \mu})$ to one in $\zero(P^{\nu, \Omega_0}_{\lambda \mu})$, we need to glue the puzzles pieces in $\overline{\Omega}_0 \setminus \Omega_0$ which share $2$-sides in the tiling together to $\fnum_2$-valued rhombus, which is the only shape of non-atomic puzzle pieces in the original $\Omega_0$. 
We call the new tiling after this process the \emph{stitching} of $\overline{t}$. Then one can easily check that stitching of $\overline{t}$ only results in the following three cases according to the direction of the shared $2$-side of two atomic puzzle pieces. All of the three stitched puzzle pieces belong to $\Omega_0$, and this means that the stitching of any tiling in $\zero(P^{\nu, \overline{\Omega}_0}_{\lambda \mu})$ becomes a tiling in $\zero(P^{\nu, \Omega_0}_{\lambda \mu})$, and thus the following theorem follows.

\begin{figure}[H]
    \begin{center}
        \begin{tikzpicture}[scale=0.8]
            \drawuptriangle{lightgreen}{0}{0}{0}{1}{2};
            \drawdowntriangle{lightgreen}{0.5}{-0.866}{2}{1}{0};
            \draw[->] (1.2,0) -- (1.8,0);
            \drawuprumbus{lightgreen}{2.5}{-0.866};

            \drawuptriangle{lightgreen}{4}{-0.433}{0}{2}{1};
            \drawdowntriangle{lightgreen}{5}{-0.433}{1}{2}{0};
            \draw[->] (5.7,0) -- (6.3,0);
            \drawSWrumbus{lightgreen}{6.5}{-0.433};

            \drawuptriangle{lightgreen}{9.5}{-0.433}{2}{1}{0};
            \drawdowntriangle{lightgreen}{9.5}{-0.433}{0}{1}{2};
            \draw[->] (10.7,0) -- (11.3,0);
            \drawNErumbus{lightgreen}{12}{-0.433};
        \end{tikzpicture}
    \end{center}
    \caption{Gluing together two atomic puzzle pieces with shared $2$-sides}
    \label{fig:enter-label}
\end{figure}
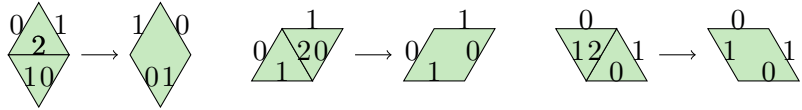

\begin{theorem}\label{thm:KT-11}
    Let $P^{\nu, \Omega_0}_{\lambda \mu}$ be an arbitrary Knutson-Tao-Woodward puzzle with $\lambda$, $\mu$, $\nu$ being partitions in $\abinom{n}{k}$ and $P^{\nu, \overline{\Omega}_0}_{\lambda \mu}$ be the puzzle induced by the atomic refinement $\overline{\Omega}_0$ of $\Omega_0$. Then there exists a one-one correspondence between $\zero(P^{\nu, \Omega_0}_{\lambda \mu})$ and $\zero(P^{\nu, \overline{\Omega}_0}_{\lambda \mu})$. 
\end{theorem}

The one-one correspondence in Theorem~\ref{thm:KT-11} is quite straightforward: For any tiling in $\zero(P^{\nu, \Omega_0}_{\lambda \mu})$, cut each rhombus puzzle piece there into two atomic puzzle pieces by adding a $2$-side in the middle; for any tiling in $\zero(P^{\nu, \overline{\Omega}_0}_{\lambda \mu})$, remove any $2$-side to stitch two atomic puzzle pieces sharing this side to a rhombus puzzle piece in $\Omega_0$. 

\begin{example}
\rm 

Let $\lambda = (0,1,0,1,0,1)$, $\mu = (0,1,0,1,0,1)$, $\nu = (1,0,1,0,1,0) \in \abinom{6}{3}$ which correspond to the partitions $(2,1), (2,1)$, and $(3,2,1)$ respectively. The equilateral triangle $\bigtriangleup^{\nu}_{\lambda \mu}$ has $N=63$ unit intervals and we assign one variable to each interval as specified in Figure~\ref{fig:0-assignment}. 

\begin{figure}[H]
    \begin{center}
        \begin{tikzpicture}[scale = 1.2]
            \foreach \x in {0,...,5} {
            \draw (0+\x/2,0.866*\x) --(6-\x/2,0.866*\x);
            }
            \foreach \x in {0,...,5} {
            \draw (\x,0) -- (3 +0.5*\x,5.196 - 0.866*\x);
            }
            \foreach \x in {0,...,5} {
            \draw (\x*0.5+0.5,0.866*\x+0.866) -- (1+\x,0);
            }

            \foreach \x in {1,...,6}{
            \foreach \y in {1,...,\x}{ 
            \pgfmathtruncatemacro{\lowvalue}{3*(\x*(\x+1)/2) -\x +\y};
            \pgfmathtruncatemacro{\leftvalue}{3*(\x*(\x-1)/2) - 1 + 2*\y};
            \pgfmathtruncatemacro{\rightvalue}{3*(\x*(\x-1)/2)  + 2*\y};
            \node at (2.25-0.5*\x + \y ,0.866*6-0.866*\x+0.433) {$x_{\leftvalue}$};
            \node at (2.75-0.5*\x + \y ,0.866*6-0.866*\x+0.433) {$x_{\rightvalue}$};
            \node at (2.5-0.5*\x + \y ,0.866*6-0.866*\x) {$x_{\lowvalue}$};
            }}
        \end{tikzpicture}
    \end{center}
    \caption{Assigned variables to intervals inside an equilateral triangle of size $6$}\label{fig:0-assignment}
\end{figure}
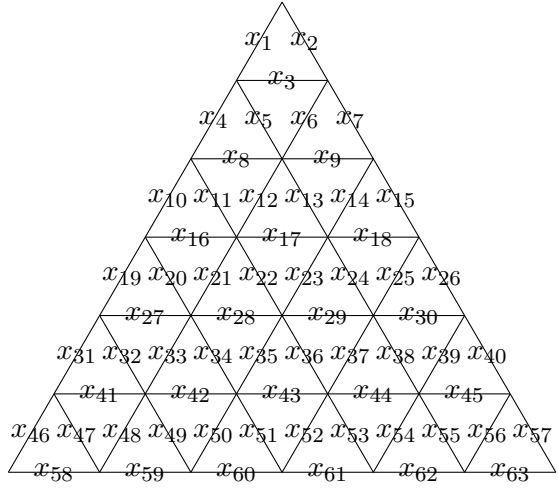

With these variables, the defining polynomials of the atomic puzzle ideal $\leftindex_a I^{\nu, \overline{\Omega}_0}_{\lambda \mu}$ in the four groups are as follows. 

\begin{itemize}
    \item ($\fnum_3$-valued) $x_i^3 - x_i$ for each $ i = 1,\ldots,63$.
    \item (matching $\lambda,\mu,\nu$) $x_{46}$, $x_{31} - 1$, $x_{19}$, $x_{10}-1$, $x_4$, and $x_1 -1$ for the left side;     $x_{2}$, $x_7 - 1$, $ x_{15}$, $x_{26} -1$, $x_{40}$, and $x_{57} -1$ for the right side; $x_{58} - 1 $, $x_{59}$, $x_{60} - 1$, $x_{61}$, $x_{62} -1$, and $x_{63}$ for the bottom side.
    \item (atomic) $x_1 + x_2 +x_3$, $x_4 + x_5 +x_8$, $x_5 + x_3 + x_6, \ldots, x_{56} + x_{57} + x_{63}$ for all the unit triangles inside $\bigtriangleup^{\nu}_{\lambda \mu}$.

    \item (distinguishing polynomial) 
    For example, $x_1 + x_1^2 + 2x_2 + x_2^2 + x_3^2 + x_1x_2x_3 + 2x_1^2x_2x_3 $ for the triangle $\Delta_{(x_1,x_2,x_3)}$ and $ x_5 + x_5^2 + 2x_3 + x_3^2 + x_6^2 + x_3x_5x_6 + 2x_3x_5^2x_6$ for $\Delta_{(x_5,x_3,x_6)}$. For each unit triangle inside $\bigtriangleup^{\nu}_{\lambda \mu}$, there is such a polynomial. 
\end{itemize}

The variety $V(I^{\nu, \overline{\Omega}_0}_{\lambda \mu})$ can be obtained by computing the lexicographic Gr\"{o}bner basis $\pset{G}$ of $\leftindex_a I^{\nu, \overline{\Omega}_0}_{\lambda \mu}$ w.r.t. the variable order $x_1 > \cdots > x_{63} $ as follows.

\begin{equation*}
    \begin{split}
   \pset{G} = \{& x_{63}, x_{62} \!-\! 1 , x_{61}, x_{60} \!-\! 1 ,x_{59}, x_{58} \!-\! 1, x_{57} \!-\! 1, x_{56} \!-\! 2, x_{55} \!-\! 1, \\
   &x_{54} \!-\! 1, x_{53}^2 \!+\! 2x_{53}, x_{53} \!+\! x_{52}, x_{53} \!+\! x_{51} \!+\! 1, 2x_{53} \!+\! x_{50}, x_{49}, x_{48}, \\
   &x_{47} \!-\! 2, x_{46}, x_{45}, x_{53} \!+\! x_{44} \!+\! 1, x_{53} \!+\! x_{43} \!+\! 2, x_{53} \!+\! x_{42},  x_{41} \!-\! 1, \\
   &x_{40}, x_{39}, x_{38} \!+\! x_{53}, x_{53} \!+\! x_{37} \!+\! 2, x_{36} \!+\! x_{53} \!+\! 2, x_{53} \!+\! x_{35} \!+\! 2, \\
   &x_{53} \!+\! x_{34} \!+\! 2, x_{53} \!+\! x_{33} \!+\! 1, x_{32} \!-\! 1, x_{31} \!+\! 2, 2x_{53} \!+\! x_{30}, x_{53} \!+\! x_{29} \!+\! 2, \\
   &x_{53} \!+\! x_{28} \!+\! 2, 2x_{53} \!+\! x_{27}, x_{26} \!-\! 1, x_{53} \!+\! x_{25} \!+\! 1, x_{53} \!+\! x_{24} \!+\! 2 ,\\
   &x_{53} \!+\! x_{23} \!+\! 2, x_{53} \!+\! x_{22} \!+\! 2, x_{53} \!+\! x_{21} \!+\! 2, x_{53} \!+\! x_{20}, x_{19}, x_{53} \!+\! x_{18}, \\
   &x_{53} \!+\! x_{17} \!+\! 2, x_{53} \!+\! x_{16} \!+\! 1, x_{15}, 2x_{53} \!+\! x_{14}, x_{53} \!+\! x_{13} \!+\! 1, x_{53} \!+\! x_{12}, \\
   &2x_{53} \!+\! x_{11}, x_{10} \!-\! 1,  x_{9} \!-\! 1, x_{8}, x_{7} \!-\! 1, x_{6} \!-\! 1, x_{5}, x_{4}, x_{3} \!-\! 2, x_{2}, x_{1} \!-\! 1 \}.     
    \end{split}
\end{equation*}
From $\pset{G}$ one can compute easily the variety $V(I^{\nu, \overline{\Omega}_0}_{\lambda \mu})$ as 
\begin{equation*}
    \begin{split}
\{(&1, 0, 2, 0, 0, 1, 1, 0, 1, 1, 0, 0, 2, 0, 0, 2, 1, 0, 0, 0, 1, 1, 1, 1, 2, 1, 0, 1, 1, 0, 1, 1, \\
&2, 1, 1, 1, 1, 0, 0, 0, 1, 0, 1, 2, 0, 0, 2, 0, 0, 0, 2, 0, 0, 1, 1, 2, 1, 1, 0, 1, 0, 1, 0), \\
&(1, 0, 2, 0, 0, 1, 1, 0, 1, 1, 1, 2, 1, 1, 0, 1, 0, 2, 0, 2, 0, 0, 0, 0, 1, 1, 1, 0, 0, 1, 1, 1, \\
&1, 0, 0, 0, 0, 2, 0, 0, 1, 2, 0, 1, 0, 0, 2, 0, 0, 1, 1, 2, 1, 1, 1, 2, 1, 1, 0, 1, 0, 1, 0 ) \}. 
    \end{split}
\end{equation*}
These two points in $V(I^{\nu, \overline{\Omega}_0}_{\lambda \mu})$ determine the values assigned to the variables and thus correspond to the two tilings in $\zero(P^{\nu, \overline{\Omega}_0}_{\lambda \mu})$ as shown in Figure~\ref{fig:0-exp-atomic} below. 

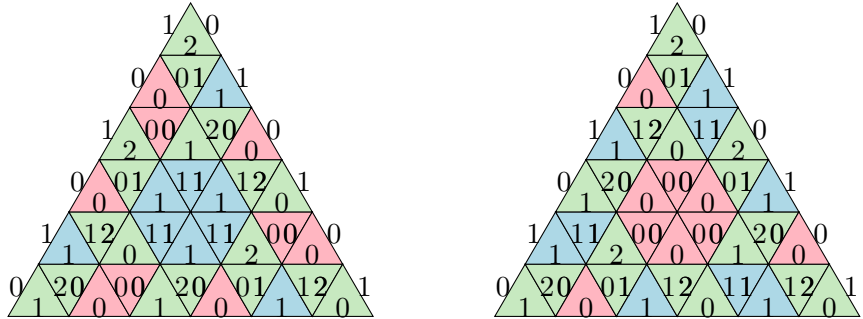
\begin{figure}[H]
\begin{center}

\begin{tikzpicture}[scale=0.8]
    \drawuptriangle{lightgreen}{0}{0}{0}{2}{1};
    \drawdowntriangle{lightgreen}{1}{0}{1}{2}{0};
    \drawuptriangle{lightred}{1}{0}{0}{0}{0};
    \drawdowntriangle{lightred}{2}{0}{0}{0}{0};
    \drawuptriangle{lightgreen}{2}{0}{0}{2}{1};
    \drawdowntriangle{lightgreen}{3}{0}{1}{2}{0};
    \drawuptriangle{lightred}{3}{0}{0}{0}{0};
    \drawdowntriangle{lightgreen}{4}{0}{2}{0}{1};
    \drawuptriangle{lightblue}{4}{0}{1}{1}{1};
    \drawdowntriangle{lightgreen}{5}{0}{0}{1}{2};
    \drawuptriangle{lightgreen}{5}{0}{2}{1}{0};

    \drawuptriangle{lightblue}{0.5}{0.866}{1}{1}{1};
    \drawdowntriangle{lightgreen}{1.5}{0.866}{0}{1}{2};
    \drawuptriangle{lightgreen}{1.5}{0.866}{2}{1}{0};
    \drawdowntriangle{lightblue}{2.5}{0.866}{1}{1}{1};
    \drawuptriangle{lightblue}{2.5}{0.866}{1}{1}{1};
    \drawdowntriangle{lightblue}{3.5}{0.866}{1}{1}{1};
    \drawuptriangle{lightgreen}{3.5}{0.866}{1}{0}{2};
    \drawdowntriangle{lightred}{4.5}{0.866}{0}{0}{0};
    \drawuptriangle{lightred}{4.5}{0.866}{0}{0}{0};

    \drawuptriangle{lightred}{1}{0.866*2}{0}{0}{0};
    \drawdowntriangle{lightgreen}{2}{0.866*2}{2}{0}{1};
    \drawuptriangle{lightblue}{2}{0.866*2}{1}{1}{1};
    \drawdowntriangle{lightblue}{3}{0.866*2}{1}{1}{1};
    \drawuptriangle{lightblue}{3}{0.866*2}{1}{1}{1};
    \drawdowntriangle{lightgreen}{4}{0.866*2}{0}{1}{2};
    \drawuptriangle{lightgreen}{4}{0.866*2}{2}{1}{0};

    \drawuptriangle{lightgreen}{1.5}{0.866*3}{1}{0}{2};
    \drawdowntriangle{lightred}{2.5}{0.866*3}{0}{0}{0};
    \drawuptriangle{lightgreen}{2.5}{0.866*3}{0}{2}{1};
    \drawdowntriangle{lightgreen}{3.5}{0.866*3}{1}{2}{0};
    \drawuptriangle{lightred}{3.5}{0.866*3}{0}{0}{0};


    \drawuptriangle{lightred}{2}{0.866*4}{0}{0}{0};
    \drawdowntriangle{lightgreen}{3}{0.866*4}{2}{0}{1};
    \drawuptriangle{lightblue}{3}{0.866*4}{1}{1}{1};


    \drawuptriangle{lightgreen}{2.5}{0.866*5}{1}{0}{2};


    \drawuptriangle{lightgreen}{0+8}{0}{0}{2}{1};
    \drawdowntriangle{lightgreen}{1+8}{0}{1}{2}{0};
    \drawuptriangle{lightred}{1+8}{0}{0}{0}{0};
    \drawdowntriangle{lightgreen}{2+8}{0}{2}{0}{1};
    \drawuptriangle{lightblue}{2+8}{0}{1}{1}{1};
    \drawdowntriangle{lightgreen}{3+8}{0}{0}{1}{2};
    \drawuptriangle{lightgreen}{3+8}{0}{2}{1}{0};
    \drawdowntriangle{lightblue}{4+8}{0}{1}{1}{1};
    \drawuptriangle{lightblue}{4+8}{0}{1}{1}{1};
    \drawdowntriangle{lightgreen}{5+8}{0}{0}{1}{2};
    \drawuptriangle{lightgreen}{5+8}{0}{2}{1}{0};

    \drawuptriangle{lightblue}{0.5+8}{0.866}{1}{1}{1};
    \drawdowntriangle{lightblue}{1.5+8}{0.866}{1}{1}{1};
    \drawuptriangle{lightgreen}{1.5+8}{0.866}{1}{0}{2};
    \drawdowntriangle{lightred}{2.5+8}{0.866}{0}{0}{0};
    \drawuptriangle{lightred}{2.5+8}{0.866}{0}{0}{0};
    \drawdowntriangle{lightred}{3.5+8}{0.866}{0}{0}{0};
    \drawuptriangle{lightgreen}{3.5+8}{0.866}{0}{2}{1};
    \drawdowntriangle{lightgreen}{4.5+8}{0.866}{1}{2}{0};
    \drawuptriangle{lightred}{4.5+8}{0.866}{0}{0}{0};

    \drawuptriangle{lightgreen}{1+8}{0.866*2}{0}{2}{1};
    \drawdowntriangle{lightgreen}{2+8}{0.866*2}{1}{2}{0};
    \drawuptriangle{lightred}{2+8}{0.866*2}{0}{0}{0};
    \drawdowntriangle{lightred}{3+8}{0.866*2}{0}{0}{0};
    \drawuptriangle{lightred}{3+8}{0.866*2}{0}{0}{0};
    \drawdowntriangle{lightgreen}{4+8}{0.866*2}{2}{0}{1};
    \drawuptriangle{lightblue}{4+8}{0.866*2}{1}{1}{1};

    \drawuptriangle{lightblue}{1.5+8}{0.866*3}{1}{1}{1};
    \drawdowntriangle{lightgreen}{2.5+8}{0.866*3}{0}{1}{2};
    \drawuptriangle{lightgreen}{2.5+8}{0.866*3}{2}{1}{0};
    \drawdowntriangle{lightblue}{3.5+8}{0.866*3}{1}{1}{1};
    \drawuptriangle{lightgreen}{3.5+8}{0.866*3}{1}{0}{2};


    \drawuptriangle{lightred}{2+8}{0.866*4}{0}{0}{0};
    \drawdowntriangle{lightgreen}{3+8}{0.866*4}{2}{0}{1};
    \drawuptriangle{lightblue}{3+8}{0.866*4}{1}{1}{1};


    \drawuptriangle{lightgreen}{2.5+8}{0.866*5}{1}{0}{2};
    
\end{tikzpicture}
\end{center}
\caption{The two tilings of the refined puzzle $P^{\nu, \overline{\Omega}_0}_{\lambda \mu}$}\label{fig:0-exp-atomic}
\end{figure}

Then we glued together all the puzzle pieces in the two tilings in Figure~\ref{fig:0-exp-atomic} which share $2$-sides to get their stitchings as shown in Figure~\ref{fig:0-exp-stitch}, which are also all the tilings in $\zero(P^{\nu, \Omega_0}_{\lambda \mu}) $.

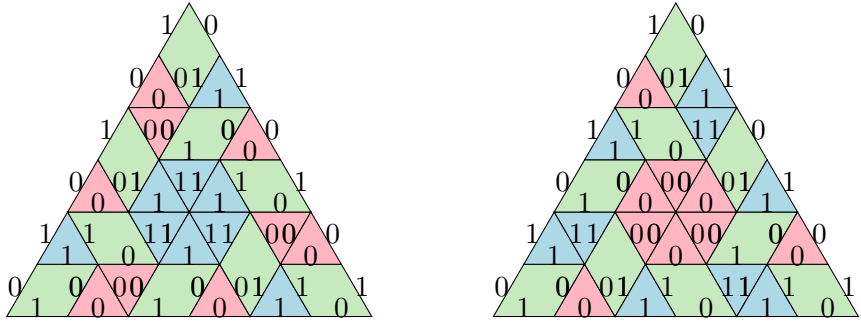
\begin{figure}
\begin{center}

\begin{tikzpicture}[scale=0.8]
    \drawSWrumbus{lightgreen}{0}{0};
    \drawuptriangle{lightred}{1}{0}{0}{0}{0};
    \drawdowntriangle{lightred}{2}{0}{0}{0}{0};
    \drawSWrumbus{lightgreen}{2}{0};
    \drawuptriangle{lightred}{3}{0}{0}{0}{0};
    \drawuprumbus{lightgreen}{4}{0}
    \drawuptriangle{lightblue}{4}{0}{1}{1}{1};
    \drawNErumbus{lightgreen}{5}{0};

    \drawuptriangle{lightblue}{0.5}{0.866}{1}{1}{1};
    \drawNErumbus{lightgreen}{1.5}{0.866};
    \drawdowntriangle{lightblue}{2.5}{0.866}{1}{1}{1};
    \drawuptriangle{lightblue}{2.5}{0.866}{1}{1}{1};
    \drawdowntriangle{lightblue}{3.5}{0.866}{1}{1}{1};
    \drawdowntriangle{lightred}{4.5}{0.866}{0}{0}{0};
    \drawuptriangle{lightred}{4.5}{0.866}{0}{0}{0};

    \drawuptriangle{lightred}{1}{0.866*2}{0}{0}{0};
    \drawuprumbus{lightgreen}{2}{0.866*2};
    \drawuptriangle{lightblue}{2}{0.866*2}{1}{1}{1};
    \drawdowntriangle{lightblue}{3}{0.866*2}{1}{1}{1};
    \drawuptriangle{lightblue}{3}{0.866*2}{1}{1}{1};
    \drawNErumbus{lightgreen}{4}{0.866*2};

    \drawdowntriangle{lightred}{2.5}{0.866*3}{0}{0}{0};
    \drawSWrumbus{lightgreen}{2.5}{0.866*3};
    \drawuptriangle{lightred}{3.5}{0.866*3}{0}{0}{0};


    \drawuptriangle{lightred}{2}{0.866*4}{0}{0}{0};
    \drawuprumbus{lightgreen}{3}{0.866*4};
    \drawuptriangle{lightblue}{3}{0.866*4}{1}{1}{1};




    \drawSWrumbus{lightgreen}{0+8}{0};
    \drawuptriangle{lightred}{1+8}{0}{0}{0}{0};
    \drawuprumbus{lightgreen}{2+8}{0};
    \drawuptriangle{lightblue}{2+8}{0}{1}{1}{1};
    \drawNErumbus{lightgreen}{3+8}{0};
    \drawdowntriangle{lightblue}{4+8}{0}{1}{1}{1};
    \drawuptriangle{lightblue}{4+8}{0}{1}{1}{1};
    \drawNErumbus{lightgreen}{5+8}{0};

    \drawuptriangle{lightblue}{0.5+8}{0.866}{1}{1}{1};
    \drawdowntriangle{lightblue}{1.5+8}{0.866}{1}{1}{1};
    \drawdowntriangle{lightred}{2.5+8}{0.866}{0}{0}{0};
    \drawuptriangle{lightred}{2.5+8}{0.866}{0}{0}{0};
    \drawdowntriangle{lightred}{3.5+8}{0.866}{0}{0}{0};
    \drawSWrumbus{lightgreen}{3.5+8}{0.866};
    \drawuptriangle{lightred}{4.5+8}{0.866}{0}{0}{0};

    \drawSWrumbus{lightgreen}{1+8}{0.866*2};
    \drawuptriangle{lightred}{2+8}{0.866*2}{0}{0}{0};
    \drawdowntriangle{lightred}{3+8}{0.866*2}{0}{0}{0};
    \drawuptriangle{lightred}{3+8}{0.866*2}{0}{0}{0};
    \drawuprumbus{lightgreen}{4+8}{0.866*2};
    \drawuptriangle{lightblue}{4+8}{0.866*2}{1}{1}{1};

    \drawuptriangle{lightblue}{1.5+8}{0.866*3}{1}{1}{1};
    \drawNErumbus{lightgreen}{2.5+8}{0.866*3};
    \drawdowntriangle{lightblue}{3.5+8}{0.866*3}{1}{1}{1};


    \drawuptriangle{lightred}{2+8}{0.866*4}{0}{0}{0};
    \drawuprumbus{lightgreen}{3+8}{0.866*4};
    \drawuptriangle{lightblue}{3+8}{0.866*4}{1}{1}{1};


    
\end{tikzpicture}
    
\end{center}
    \caption{The two tilings of the original puzzle $P^{\nu, \Omega_0}_{\lambda \mu}$}\label{fig:0-exp-stitch}
\end{figure}

\end{example}

\section{Forbidding puzzle ideal for $T$-equivariant puzzle}
\label{sec:t}

The set $\Omega_0$ of puzzle pieces for the Knutson-Tao-Woodward puzzles is the simplest, enabling us to derive the one-one correspondence in Theorem~\ref{thm:KT-11}. Now let us take a look at $\Omega_T$ for the $T$-equivariant puzzles to see what happens if we follow the same process as for the Knutson-Tao-Woodward puzzles. 

First $\Omega_T$ is obtained by adjoining $\Omega_0$ with the equivariant piece, and thus the atomic refinement of $\Omega_T$ results in two additional pieces as shown below. 

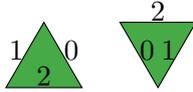
\begin{figure}[H]
    \begin{center}
        \begin{tikzpicture}
            \drawuptriangle{dlgreen}{0}{0}{1}{0}{2};
            \drawdowntriangle{dlgreen}{2}{0}{2}{0}{1};
        \end{tikzpicture}
    \end{center}
    \caption{Additional atomic puzzle pieces in $\overline{\Omega}_T$}
    \label{fig:enter-label}
\end{figure}

Then for any $T$-equivariant puzzle $P^{\nu, \Omega_T}_{\lambda \mu}$ we obtain a puzzle $P^{\nu, \overline{\Omega}_T}_{\lambda \mu}$ induced by the atomic refinement $\overline{\Omega}_T$. The system of polynomial equations can be similarly constructed, and the only difference from the Knutson-Tao-Woodward puzzles is the construction of the distinguishing polynomial, for we have two more allowed pieces which break the closeness with rotation. 

To restrict an atomic puzzle piece to lie in $\overline{\Omega}_T$, we need to consider the two cases when the unit triangle of the piece is upward or downward. For each upward unit triangle $\bigtriangleup_k$ inside $\bigtriangleup^{\nu}_{\lambda \mu}$, let $x_{k1}$, $x_{k2}$, and $x_{k3}$ be the variables corresponding to its left, right, and bottom sides. In the same way by solving an underdetermined linear system, we can construct a distinguishing polynomial 
$$f^k(x_{k1}, x_{k2}, x_{k3}) = x_{k1}^2x_{k2} + 2x_{k1}x_{k2}^2 + 2x_{k1}^2 + x_{k1}$$
such that $f^k(\overline{x}_{k1}, \overline{x}_{k2}, \overline{x}_{k3}) = 0$ if the assigned values $\overline{x}_{k1}$, $\overline{x}_{k2}$, and $\overline{x}_{k3}$ correspond to an upward atomic puzzle piece in $\overline{\Omega}_T$ and $f^k(\overline{x}_{k1}, \overline{x}_{k2}, \overline{x}_{k3}) \neq 0$ (for example $=1$) otherwise. Similarly for each downward unit triangle $\bigtriangledown_k$ with $x_{k1}$, $x_{k2}$, and $x_{k3}$ corresponding to its left, right, and top sides, a distinguishing polynomial is 
$$f_k(x_{k1}, x_{k2}, x_{k3}) = x_{k1}x_{k2}^2 + 2 x_{k1}^2 + x_{k1}x_{k2} + 2 x_{k2}^2 + x_{k1} + 2 x_{k2}.$$

From the defining polynomials of the system above we have the atomic puzzle ideal $\leftindex_a I^{\nu, \overline{\Omega}_T}_{\lambda \mu}$, and clearly the inequality $\#\var(\leftindex_a I^{\nu, \overline{\Omega}_T}_{\lambda \mu})  = \#\zero(P^{\nu, \overline{\Omega}_T}_{\lambda \mu}) \geq \#\zero(P^{\nu, \Omega_T}_{\lambda \mu})$ also holds. The fundamental difference lies in the process of stitching atomic puzzle pieces sharing $2$-sides: stitching all the puzzle pieces in $\overline{\Omega}_T \setminus \Omega_T$ to $\fnum_2$-valued rhombuses will result in two new $\fnum_2$-valued puzzle pieces not contained in $\Omega_T$, as shown in Figure~\ref{fig:t-forbidden} below. We call these two new pieces the \emph{forbidden puzzle pieces}, and the existence of forbidden pieces means that any tiling in $\zero(P^{\nu, \overline{\Omega}_T}_{\lambda \mu})$ for which its stitching contains either of the forbidden pieces does not correspond to a tiling in $\zero(P^{\nu, \Omega_T}_{\lambda \mu})$. As a result, we only have a strict inequality $\#\zero(P^{\nu, \overline{\Omega}_T}_{\lambda \mu}) \gneq \#\zero(P^{\nu, \Omega_T}_{\lambda \mu})$ now, and finding $\var(\leftindex_a I^{\nu, \overline{\Omega}_T}_{\lambda \mu})$ is not enough to recover $\zero(P^{\nu, \Omega_T}_{\lambda \mu})$.

\begin{figure}[H]
\begin{center}
    \begin{tikzpicture}
        \drawuptriangle{dlgreen}{0}{0}{1}{0}{};
        \drawdowntriangle{lightgreen}{0.5}{-0.866}{}{1}{0};

        \drawuptriangle{lightgreen}{3}{0}{0}{1}{};
        \drawdowntriangle{dlgreen}{3.5}{-0.866}{}{0}{1};
    \end{tikzpicture}
\end{center}
\caption{The forbidden puzzle pieces for $P^{\nu, \Omega_T}_{\lambda \mu}$}\label{fig:t-forbidden}
\end{figure}
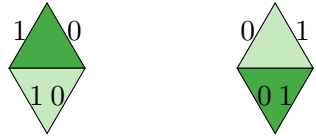

This observation leads us to add a new group of polynomial equations to remove the existence of forbidden pieces in stitching of the tilings in $\zero(P^{\nu, \overline{\Omega}_T}_{\lambda \mu})$. This process is similar to the construction of the distinguishing polynomial, except that now we need to deal with rhombuses. Note that both of the forbidden pieces are of the shape of a vertical rhombus and they are obtained by stitching a pair of  upward and downward atomic puzzle pieces sharing a $2$-side. All the possible $\fnum_3$-valued vertical rhombus pieces consisting of atomic puzzle pieces in $\overline{\Omega}_T$ are demonstrated in the picture below, with the two forbidden puzzle pieces in the rightmost column. 

\begin{figure}[H]
\begin{center}
    \begin{tikzpicture}
        \drawuptriangle{lightred}{0}{0}{0}{0}{0};
        \drawdowntriangle{lightred}{0.5}{-0.866}{0}{0}{0};

        \drawuptriangle{lightred}{1.5}{0}{0}{0}{0};
        \drawdowntriangle{lightgreen}{2}{-0.866}{0}{1}{2};

        \drawuptriangle{lightblue}{3}{0}{1}{1}{1};
        \drawdowntriangle{lightblue}{3.5}{-0.866}{1}{1}{1};
            
        \drawuptriangle{lightblue}{4.5}{0}{1}{1}{1};
        \drawdowntriangle{lightgreen}{5}{-0.866}{1}{2}{0};

        \drawuptriangle{dlgreen}{6}{0}{1}{0}{2};
        \drawdowntriangle{dlgreen}{6.5}{-0.866}{2}{0}{1};

        \drawuptriangle{dlgreen}{7.5}{0}{1}{0}{2};
        \drawdowntriangle{lightgreen}{8}{-0.866}{}{1}{0};

        \drawuptriangle{lightgreen}{0}{-2}{2}{1}{0};
        \drawdowntriangle{lightred}{0.5}{-2.866}{0}{0}{0};

        \drawuptriangle{lightgreen}{1.5}{-2}{2}{1}{0};
        \drawdowntriangle{lightgreen}{2}{-2.866}{0}{1}{2};

        \drawuptriangle{lightgreen}{3}{-2}{0}{2}{1};
        \drawdowntriangle{lightblue}{3.5}{-2.866}{1}{1}{1};
            
        \drawuptriangle{lightgreen}{4.5}{-2}{0}{2}{1};
        \drawdowntriangle{lightgreen}{5}{-2.866}{1}{2}{0};

        \drawuptriangle{lightgreen}{6}{-2}{0}{1}{2};
        \drawdowntriangle{lightgreen}{6.5}{-2.866}{2}{1}{0};

        \drawuptriangle{lightgreen}{7.5}{-2}{0}{1}{2};
        \drawdowntriangle{dlgreen}{8}{-2.866}{}{0}{1};
        
    \end{tikzpicture}
\end{center}
\caption{All the vertical rhombus puzzle pieces consisting of atomic puzzle pieces in $\overline{\Omega}_T$}\label{fig:t-rhombus}
\end{figure}
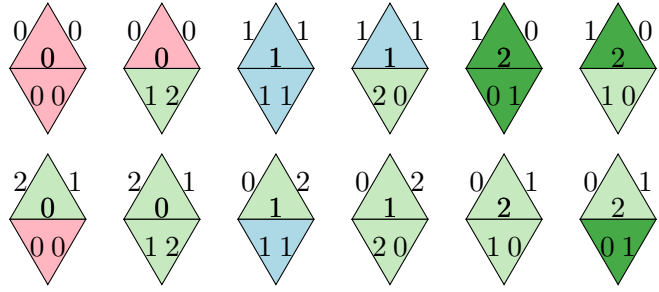

As a result, to ensure that neither of these two forbidden pieces appears in the stitching of a tiling in $\zero(P^{\nu, \overline{\Omega}_T}_{\lambda \mu})$, we need to forbid the assigned values to any vertical rhombus inside $\bigtriangleup^{\nu}_{\lambda \mu}$ to make it either of the forbidden puzzle piece. For each vertical rhombus $\lozenge_{\ell}~(\ell = 1, \ldots, \frac{n(n-1)}{2})$ inside $\bigtriangleup^{\nu}_{\lambda \mu}$, let $x_{\ell 1}$, $x_{\ell 2}$, and $x_{\ell 3}$ be the the variables assigned to the left, right, and bottom sides of its upward atomic puzzle piece, and $x_{\ell 4}$ and $x_{\ell 5}$ be those assigned to the left and right sides of the downward atomic puzzle piece. Then we want to construct a polynomial function $f_{\ell}(x_{\ell 1}, x_{\ell 2}, x_{\ell 3}, x_{\ell 4}, x_{\ell 5})$ over $\fnum_3$ such that 
\begin{equation*}
    \begin{split}
   &f_{\ell}(0,0,0,0,0) = f_{\ell}(0,0,0,1,2) = f_{\ell}(1,1,1,2,0) = f_{\ell}(1,0,2,0,1) = \\
   &f_{\ell}(1,1,1,1,1) = f_{\ell}(2,1,0,0,0) = f_{\ell}(2,1,0,1,2) = f_{\ell}(0,2,1,1,1) = \\
   &f_{\ell}(0,2,1,2,0) = f_{\ell}(0,1,2,1,0) = 0, \quad 
   f_{\ell}(1,0,2,1,0) = f_{\ell}(0,1,2,0,1) =1.
    \end{split}
\end{equation*}
Again this forms an underdetermined linear system with the undetermined coefficients of $f_{\ell}$ as the variables, and solving this system furnishes a choice of $f_{\ell}$ as 
$$ f_{\ell} = x_{\ell 1}^2x_{\ell 4} + x_{\ell 2}^2x_{\ell 4} + 2x_{\ell 1}^2 + 2x_{\ell 1}x_{\ell 2} + 2x_{\ell 2}^2 + x_{\ell 2}x_{\ell 4} + x_{\ell 1} + 2x_{\ell 2}.$$
We call $f_{\ell}$ above a \emph{forbidding polynomial} w.r.t. $\overline{\Omega}_T$. 

Let $\pset{F}_a$ be the defining polynomials of the atomic puzzle ideal $\leftindex_a I^{\nu, \overline{\Omega}_T}_{\lambda \mu}$ and $\pset{F}_f$ be the set of all the forbidding polynomials for all vertical rhombuses inside $\bigtriangleup^{\nu}_{\lambda \mu}$. Now consider the ideal $\bases{\pset{F}_a \cup \pset{F}_f} \subseteq \fnum_3[\p{x}]$. We call this ideal the \emph{forbidding puzzle ideal} of $P^{\nu, \Omega_T}_{\lambda \mu}$ and denote it by $\leftindex_f I^{\nu, \overline{\Omega}_T}_{\lambda \mu}$. Clearly we have $\var(\leftindex_f I^{\nu, \overline{\Omega}_T}_{\lambda \mu}) \subseteq \var(\leftindex_a I^{\nu, \overline{\Omega}_T}_{\lambda \mu})$. Note that $\var(\leftindex_a I^{\nu, \overline{\Omega}_T}_{\lambda \mu})$ one-one corresponds to $\zero(P^{\nu, \overline{\Omega}_T}_{\lambda \mu})$, and by the construction of $\leftindex_f I^{\nu, \overline{\Omega}_T}_{\lambda \mu}$ we know that $\var(\leftindex_f I^{\nu, \overline{\Omega}_T}_{\lambda \mu})$ corresponds to the tilings in $\zero(P^{\nu, \overline{\Omega}_T}_{\lambda \mu})$ such that neither of the two forbidden pieces appears in the stitching of the tilings. Denote by $\zero_f(P^{\nu, \overline{\Omega}_T}_{\lambda \mu})$ this subset of $\zero(P^{\nu, \overline{\Omega}_T}_{\lambda \mu})$. Then we have the following theorem. 

\begin{theorem}\label{thm:t-11}
    Let $P^{\nu, \Omega_T}_{\lambda \mu}$ be an arbitrary $T$-equivariant puzzle with $\lambda$, $\mu$, $\nu$ being partitions in $\abinom{n}{k}$ and $P^{\nu, \overline{\Omega}_T}_{\lambda \mu}$ be the puzzle induced by the atomic refinement $\overline{\Omega}_T$ of $\Omega_T$. Then there exists a one-one correspondence between $\zero(P^{\nu, \Omega_T}_{\lambda \mu})$ and $\zero_f(P^{\nu, \overline{\Omega}_T}_{\lambda \mu})$. 
\end{theorem}

The one-one correspondence in Theorem~\ref{thm:t-11} is exactly the same as the one in Theorem~\ref{thm:KT-11}, stated below that theorem. Let us go back to our treatment of the Knutson-Tao-Woodward puzzle $P^{\nu, \Omega_0}_{\lambda \mu}$ in Section~\ref{sec:kt}. Form the viewpoint of forbidding puzzle ideals, we know that for $P^{\nu, \Omega_0}_{\lambda \mu}$, its forbidding puzzle ideal $\leftindex_f I^{\nu, \Omega_0}_{\lambda \mu}$ is equal to its atomic one $\leftindex_a I^{\nu, \overline{\Omega}_0}_{\lambda \mu}$, for the set of forbidding polynomials is empty. This explains why it suffices to use only the atomic puzzle ideal to solve the Knutson-Tao-Woodward puzzle.

\section{Implying puzzle ideals for $K$-theoretic puzzles}
\label{sec:k}
Either of $\Omega_0$ and $\Omega_T$ in the two previous sections only contains puzzle pieces in the shape of unit triangle and rhombus. In this section we study the $K$-theoretic puzzles, in particular those with the set $\Omega_C$ of puzzle pieces for example. 

Note that for the $C$-piece, there are three tilings of it with the atomic puzzle pieces as shown in Figure~\ref{fig:C-tiling} below. To obtain an atomic refinement of the $C$-piece, we can take any of the tiling, and in the sequel we fix the first one in Figure~\ref{fig:C-tiling}, and the resulting atomic refinement $\overline{\Omega}_C$ induced by this tiling is shown in Figure~\ref{fig:Cpuzzle-refined} below. 

\begin{figure}[H]
    \begin{center}
        \begin{tikzpicture}
            \drawuptriangle{lightblue}{0+4}{0}{1}{1}{1};
            \drawdowntriangle{lightgreen}{0+4}{0}{2}{0}{1};
            \drawdowntriangle{dlgreen}{1+4}{0}{2}{1}{0};
            \drawuptriangle{lightgreen}{-0.5+4}{0.866}{1}{0}{2};
            \drawdowntriangle{lightred}{0.5+4}{0.866}{0}{0}{0};
            \drawuptriangle{dlgreen}{0.5+4}{0.866}{0}{1}{2};

            \drawuptriangle{dlgreen}{0}{0}{2}{0}{1};
            \drawdowntriangle{dlgreen}{0}{0}{1}{0}{2};
            \drawdowntriangle{lightred}{1}{0}{0}{0}{0};
            \drawuptriangle{lightblue}{-0.5}{0.866}{1}{1}{1};
            \drawdowntriangle{lightgreen}{0.5}{0.866}{0}{1}{2};
            \drawuptriangle{lightgreen}{0.5}{0.866}{2}{1}{0};

            \drawuptriangle{lightgreen}{0+8}{0}{0}{2}{1};
            \drawdowntriangle{lightred}{0+8}{0}{0}{0}{0};
            \drawdowntriangle{lightgreen}{1+8}{0}{1}{2}{0};
            \drawuptriangle{dlgreen}{-0.5+8}{0.866}{1}{2}{0};
            \drawdowntriangle{dlgreen}{0.5+8}{0.866}{0}{2}{1};
            \drawuptriangle{lightblue}{0.5+8}{0.866}{1}{1}{1};
        \end{tikzpicture}
    \end{center}
    \caption{Three tiling of the $C$-piece}
    \label{fig:C-tiling}
\end{figure}
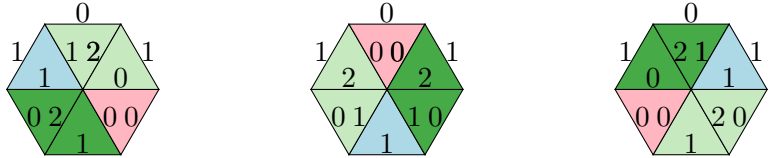

\begin{figure}[H]
\begin{center}
    \begin{tikzpicture}[scale=1]
        \drawuptriangle{lightred}{0}{0}{0}{0}{0};
        \drawuptriangle{lightblue}{1.5}{0}{1}{1}{1};
        \drawuptriangle{lightgreen}{3}{0}{1}{0}{2};
        \drawuptriangle{lightgreen}{4.5}{0}{2}{1}{0};
        \drawuptriangle{lightgreen}{6}{0}{0}{2}{1};
        \drawuptriangle{dlgreen}{7.5}{0}{2}{0}{1};

        \drawdowntriangle{lightred}{0.5}{-2}{0}{0}{0};
        \drawdowntriangle{lightblue}{2}{-2}{1}{1}{1};
        \drawdowntriangle{lightgreen}{3.5}{-2}{2}{0}{1};
        \drawdowntriangle{lightgreen}{5}{-2}{0}{1}{2};
        \drawdowntriangle{lightgreen}{6.5}{-2}{1}{2}{0};
        \drawdowntriangle{dlgreen}{8}{-2}{1}{0}{2};
    \end{tikzpicture}
\end{center}
\caption{Atomic puzzle pieces in $\overline{\Omega}_C$}\label{fig:Cpuzzle-refined}
\end{figure}
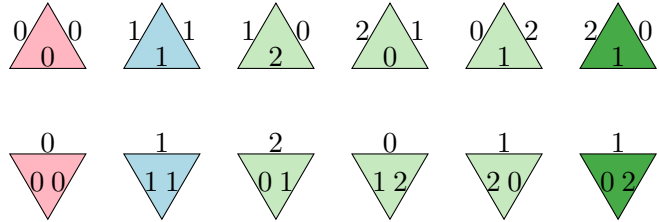

The construction of the four groups of defining polynomials of the atomic puzzle ideal $\leftindex_a I^{\nu, \overline{\Omega}_C}_{\lambda \mu}$ is the same as that for the $T$-equivariant puzzles, and we only present the upward and downward distinguishing polynomials for each triangle inside $\bigtriangleup^{\nu}_{\lambda \mu}$:
\begin{equation*}
    \begin{split}
        f^{k}(x_{k1}, x_{k2}, x_{k3}) &= x_{k2}^2x_{k3} + 2x_{k2}x_{k3}^2 + 2x_{k2}^2 + x_{k2}, \\
        f_k(x_{k1}, x_{k2}, x_{k3}) &= x_{k2}x_{k3}^2 + 2x_{k2}^2 + x_{k2}x_{k3} + 2x_{k3}^2 + x_{k2} + 2x_{k3}.
    \end{split}
\end{equation*}

Next let us stitch all the puzzle pieces in $\overline{\Omega}_C \setminus \Omega_C$ to $\fnum_2$-valued rhombuses and this will result in three new $\fnum_2$-valued puzzle pieces not contained in $\Omega_C$, as shown in Figure~\ref{fig:3newrho} below. But if one takes a further look at the tiling of the $C$-piece we choose in Figure~\ref{fig:C-tiling}, he will see that the last new rhombus puzzle piece is indeed hidden in the tiling: stitching the two atomic puzzle pieces at the bottom by removing their shared 2-side will result in this rhombus puzzle piece. On one hand, this rhombus puzzle piece, though not contained in $\Omega_C$ directly, should not be considered as a forbidden one; on the other, if we want to turn a tiling $\overline{t}$ in $\zero(P^{\nu, \overline{\Omega}_T}_{\lambda \mu})$ to one in $\zero(P^{\nu, \Omega_T}_{\lambda \mu})$, any occurrence of this puzzle piece in the stitching of $\overline{t}$ implies that the $C$-piece containing it must also appear in the tiling of $\bigtriangleup^{\nu}_{\lambda \mu}$ with $\Omega_C$. We call this rhombus puzzle piece an \emph{implicit} one. 

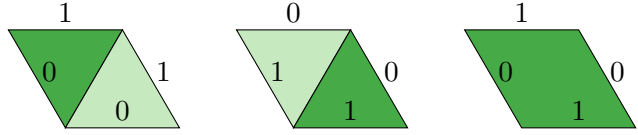
\begin{figure}[H]
    \begin{center}
        \begin{tikzpicture}[scale=1.5]
        \drawdowntriangle{dlgreen}{0}{0}{1}{0}{};
        \drawuptriangle{lightgreen}{0}{0}{}{1}{0};

        \drawdowntriangle{lightgreen}{2}{0}{0}{1}{};
        \drawuptriangle{dlgreen}{2}{0}{}{0}{1};
            \drawNErumbusReverse{dlgreen}{4}{0};
        \end{tikzpicture}
    \end{center}
    \caption{New rhombus pieces after stitching}
    \label{fig:3newrho}
\end{figure}

The first two puzzle pieces in Figure~\ref{fig:3newrho} are the forbidden ones and both of them have their middle interval parallel to the left side of $\bigtriangleup^{\nu}_{\lambda \mu}$, and we call such rhombus left ones. For each left rhombus ~\lrhombus~ inside $\bigtriangleup^{\nu}_{\lambda \mu}$, let $x_{\ell 1}, \ldots, x_{\ell 5}$ be the variables of its unit intervals as shown in Figure~\ref{fig:var4rho} below. Then one forbidding polynomials for ~~\lrhombus~ can be constructed as 
\begin{equation*}
    \begin{split}
f_{\ell}(x_{\ell 1}, &x_{\ell 2},x_{\ell 3},x_{\ell 4},x_{\ell 5}) =\\
&x_{\ell 1}^2x_{\ell 4} + x_{\ell 2}^2x_{\ell 4} + 2x_{\ell 1}^2 + 2x_{\ell 1}x_{\ell 2} + 2x_{\ell 2}^2 + x_{\ell 2}x_{\ell 4} + x_{\ell 1} + 2x_{\ell 2}.        
    \end{split}
\end{equation*}

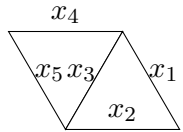
\begin{figure}[H]

\begin{center}
    \begin{tikzpicture}[scale=1.5]

        \drawdowntriangle{white}{4}{0}{$x_4$}{$x_5$}{$x_3$};
        \drawuptriangle{white}{4}{0}{}{$x_1$}{$x_2$};
    \end{tikzpicture}
\end{center}
\caption{The variables for a left rhombus}
\label{fig:var4rho} 
\end{figure}

Now let us work on this particular implicit rhombus puzzle piece. Note that it only appears inside the $C$-piece in the stitching of a tiling of $\bigtriangleup^{\nu}_{\lambda \mu}$ with $\overline{\Omega}_C$. This means that its occurrence implies that of the $C$-piece, and we need to construct a polynomial for this implication. We call it an \emph{implying} polynomial. The $C$-piece is an $\fnum_2$-value puzzle piece of the shape of a hexagon, which contains $12$ unit intervals, and thus the implying polynomial we want to construct has $12$ variables at most. 
Any hexagon contains two left rhombuses inside, and with the atomic puzzle pieces in $\overline{\Omega}_C$, the assigned values of the five unit intervals of each of these two left rhombuses can be determined by two of them, namely by the one assigned to the middle unit interval and any of the two assigned to the unit intervals inside the hexagon. For each hexagon $\hexagon_m$, let $x_{m1}$, $x_{m2}$, $x_{m3}$, and $x_{m4}$ be the variables indicated in the figure below. Then one can observe that if we fix the values of these 4 variables, the values for the remaining variables of an $\fnum_2$-valued hexagon $\hexagon_m$ tiled with $\overline{\Omega}_C$ are also fixed. To see this, one needs to bear in mind that the unit intervals on the boundary of the hexagon are assigned $\fnum_2$-values. 

\begin{figure}[H]
\begin{center}
    \begin{tikzpicture}[scale=1.5]
            \drawuptriangle{white}{0}{0}{$x_4$}{}{};
            \drawdowntriangle{white}{0}{0}{$x_3$}{}{$x_4$};
            \drawdowntriangle{white}{1}{0}{}{}{};
            \drawuptriangle{white}{-0.5}{0.866}{}{$x_2$}{$x_3$};
            \drawdowntriangle{white}{0.5}{0.866}{}{$x_2$}{$x_1$};
            \drawuptriangle{white}{0.5}{0.866}{}{}{};
    \end{tikzpicture}
\end{center}
\caption{The hexagon with $4$ variables}
\label{fig:20} 
\end{figure}
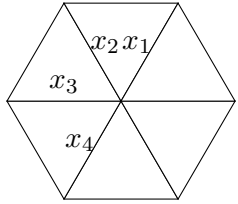

Remember that we want to construct an implying polynomial to represent the condition that the occurrence of the implicit rhombus puzzle piece implies that of the $C$-piece. With the values assigned to an $\fnum_2$-valued $\hexagon_m$ determined by those of the four variables $x_{m1}$, $x_{m2}$, $x_{m3}$, and $x_{m4}$, the above implication is equivalent to that $x_{m3} = 1$ and $x_{m4} = 2$ imply $x_{m1} = 2$ and $x_{m2} = 1$. Then the implying polynomial $g_m(x_{m1}, x_{m2}, x_{m3}, x_{m4})$ for $\hexagon_m$ should satisfies $g_m(2, 1, 1, 2) = 0$ and $g_m(\overline{x}_{m1}, \overline{x}_{m2}, 1, 2) \neq 0$ for any values $(\overline{x}_{m1}, \overline{x}_{m1}) \in \fnum_3^2 \setminus \{(1, 2)\} $. Again by solving a system of linear equations determined by the $9$ constraints above, we can have an implying polynomial for $\hexagon_m$ as 
\begin{equation*}
    \begin{split}
&g_m(x_{m1},x_{m2},x_{m3},x_{m4}) =\\
&2x_{m1}^2x_{m2}^2x_{m3}^2x_{m4}^2 + x_{m1}^2x_{m2}^2x_{m3}^2x_{m4} + 2x_{m1}^2x_{m2}^2x_{m3}x_{m4}^2 + 2x_{m1}^2x_{m2}x_{m3}^2x_{m4}^2 + \\
&x_{m1}x_{m2}^2x_{m3}^2x_{m4}^2 + x_{m1}^2x_{m2}^2x_{m3}x_{m4} + x_{m1}^2x_{m2}x_{m3}^2x_{m4} + 2x_{m1}^2x_{m2}x_{m3}x_{m4}^2 +\\
&2x_{m1}x_{m2}^2x_{m3}^2x_{m4} + x_{m1}x_{m2}^2x_{m3}x_{m4}^2 + x_{m1}x_{m2}x_{m3}^2x_{m4}^2 + x_{m1}^2x_{m2}x_{m3}x_{m4} +\\
&2x_{m1}^2x_{m2}^2x_{m3}x_{m4} + 2x_{m1}x_{m2}x_{m3}^2x_{m4} + x_{m1}x_{m2}x_{m3}x_{m4}^2 + 2x_{m1}x_{m2}x_{m3}x_{m4} +\\
&x_{m3}^2x_{m4}^2 + 2x_{m3}^2x_{m4} + x_{m3}x_{m4}^2 + 2x_{m3}x_{m4}.
    \end{split}
\end{equation*}

Let $\pset{F}_a$ be the defining polynomials of the atomic puzzle ideal $\leftindex_a I^{\nu, \overline{\Omega}_C}_{\lambda \mu}$, $\pset{F}_f$ be the set of all the forbidding polynomials for all left rhombuses inside $\bigtriangleup^{\nu}_{\lambda \mu}$, and $\pset{F}_i$ be the set of all the implying polynomials for all the hexagons inside $\bigtriangleup^{\nu}_{\lambda \mu}$. Now consider the ideal $\bases{\pset{F}_a \cup \pset{F}_f \cup \pset{F}_i} \subseteq \fnum_3[\p{x}]$. We call this ideal the \emph{implying puzzle ideal} of $P^{\nu, \Omega_C}_{\lambda \mu}$ and denote it by $\leftindex_i I^{\nu, \overline{\Omega}_C}_{\lambda \mu}$. Then $\var(\leftindex_i I^{\nu, \overline{\Omega}_C}_{\lambda \mu})$ corresponds to the tilings in $\zero(P^{\nu, \overline{\Omega}_C}_{\lambda \mu})$ such that in their stitchings neither of the two forbidden pieces occurs and the occurrence of the implicit piece implies that of the $C$-piece in the tiling of $\bigtriangleup^{\nu}_{\lambda \mu}$ with $\Omega_C$. Denote by $\zero_i(P^{\nu, \overline{\Omega}_C}_{\lambda \mu})$ this subset of $\zero(P^{\nu, \overline{\Omega}_C}_{\lambda \mu})$. Then we have the following theorem.

\begin{theorem}\label{thm:t-11}
    Let $P^{\nu, \Omega_C}_{\lambda \mu}$ be an arbitrary $K$-theoretic C-puzzle with $\lambda$, $\mu$, $\nu$ being partitions in $\abinom{n}{k}$ and $P^{\nu, \overline{\Omega}_C}_{\lambda \mu}$ be the puzzle induced by the atomic refinement $\overline{\Omega}_C$ of $\Omega_C$. Then there exists a one-one correspondence between $\zero(P^{\nu, \Omega_C}_{\lambda \mu})$ and $\zero_i(P^{\nu, \overline{\Omega}_C}_{\lambda \mu})$. 
\end{theorem}

To send a tiling $\overline{t}$ in $\zero_i(P^{\nu, \overline{\Omega}_C}_{\lambda \mu})$ to one in $\zero(P^{\nu, \Omega_C}_{\lambda \mu})$, one needs to construct the stitching $\tilde{t}$ of $\overline{t}$ and then for each occurrence of the implicit puzzle piece in $\tilde{t}$, find the four additional atomic puzzle pieces which form a $C$-piece with it and stitch them together to the $C$-piece. In this way, the tiling $\tilde{t}$ will be transformed into a tiling in $\zero(P^{\nu, \Omega_C}_{\lambda \mu})$. With the $D$-piece just a rotation of the $C$-piece by $60$ degrees, the method above naturally works with the $K$-theoretic puzzles with $\Omega_D$. 

It is time to summarize what we have done to the three kinds of puzzles before we formulate the method based on puzzle ideals for solving puzzles. The atomic puzzle ideal deals with the case when neither forbidden nor implicit rhombus puzzle piece appears after the process of stitching a tiling with atomic puzzle pieces. In this case, the defining polynomials of the ideals distinguish the allowed atomic puzzle pieces in the atomic refinement by including the distinguishing polynomials. If forbidden rhombus puzzle pieces do appear, the defining polynomials of the ideals need to be expanded with the forbidding polynomials which forbids these rhombus puzzle pieces, and then we have the forbidding puzzle ideals. For all the puzzles with the pieces of the shapes of unit triangle and rhombus, the varieties of their forbidding puzzle ideals one-one correspond to their solutions. The situation becomes trickier when the puzzle pieces are of the shapes of unit triangle, unit rhombus, and another bigger polygon, for we have to deal with the implicit rhombus puzzle pieces which are neither original pieces for the puzzle nor the forbidden pieces after stitching a tiling with atomic puzzle pieces. For the implicit pieces we further introduce the implying polynomials to define the implying puzzle ideals, whose variety again one-one correspond to the solutions of the original puzzles under some conditions we will soon reveal in the next section. 

\section{Puzzle ideal}
\label{sec:puzzle}

In this section, we formulate a general framework to construct the puzzle ideals for the puzzles under certain conditions such that the varieties of the puzzle ideals one-one corresponds to the solutions of the puzzles. 

\subsection{Atomic refinement}

Let $P^{\nu, \Omega}_{\lambda \mu}$ be a puzzle consisting of the $\fnum_2$-valued $\bigtriangleup^{\nu}_{\lambda \mu}$ and the set $\Omega$ of $\fnum_2$-valued puzzle pieces. Recall that an atomic puzzle piece is an $\fnum_3$-valued unit triangle whose assigned values sum up to $0$ in $\fnum_3$. For any puzzle piece $p \in \Omega$, let $\overline{t}_p$ be a tiling of $p$ with a set of atomic puzzle pieces, if it exists (in this process the $\fnum_2$-value of $p$ is naturally embedded in $\fnum_3$). In this case, $p$ is said to be \emph{refinable}, and the set of atomic puzzle pieces effectively appearing in $\overline{t}_p$ is called the \emph{atomic refinement} of $p$ induced by $\overline{t}_p$ and denoted by $\Omega(\overline{t}_p)$. The set $\Omega$ is said to be \emph{refinable} if each $p \in \Omega$ is so. In this case, let $\overline{t}_p$ be such a tiling of $p$. Then the set $\overline{\Omega} = \bigcup_{p \in \Omega} \Omega(\overline{t}_p)$ is called the \emph{atomic refinement} of $\Omega$ induced by $\{\overline{t}_p: p \in \Omega\}$. When we are not particularly interested in the tilings, we say that $\overline{\Omega}$ is an atomic refinement of $\Omega$ for short. Since tilings of the puzzle pieces in $\Omega$ may not be unique, in general the atomic refinement of $\Omega$ is also not unique.

\begin{remark}\em
    For atomic refinement, we have the following straightforward observations. 
    \begin{enumerate}
        \item[(1)] The only $\fnum_2$-valued atomic puzzle pieces are the four illustrated in Figure~\ref{fig:0-piece}, and thus any $\Omega$ containing an $\fnum_2$-value piece of unit triangle other than these four is not refinable. 
        \item[(2)] The six sets $\Omega_0$, $\Omega_T$, $\Omega_A$, $\Omega_B$, $\Omega_C$, and $\Omega_D$ of puzzle pieces are all refinable. In particular, atomic refinement of $\Omega_0$, $\Omega_T$, $\Omega_A$, and $\Omega_B$ is unique and there are three different atomic refinements of $\Omega_C$ and $\Omega_D$. 
        \end{enumerate}
\end{remark}

Clearly our treatment in the previous sections only works with puzzles with refinable sets of puzzle pieces, so hereafter when referring to a puzzle, we assume that it has a refinable set of puzzle pieces. 

Let $P^{\nu, \Omega}_{\lambda \mu}$ be a puzzle and $\overline{\Omega}$ be an atomic refinement of $\Omega$ induced by the tilings $\{\overline{t}_p: p \in \Omega\}$. Then the puzzle $P^{\nu, \overline{\Omega}}_{\lambda \mu}$ consisting of $\bigtriangleup^{\nu}_{\lambda \mu}$ and the set $\overline{\Omega}$ of atomic puzzle pieces is called the \emph{refined puzzle} of $P^{\nu, \Omega}_{\lambda \mu}$ w.r.t. $\overline{\Omega}$. Recall that we use $\zero(P^{\nu, \Omega}_{\lambda \mu})$ to denote the set of all the solutions of $P^{\nu, \Omega}_{\lambda \mu}$, namely all the possible tilings of $\bigtriangleup^{\nu}_{\lambda \mu}$ with $\Omega$. Then for each tiling in $\zero(P^{\nu, \Omega}_{\lambda \mu})$, if we replace any puzzle piece $p \in \Omega$ in this tiling with the tiling $\overline{t}_p$, we will have a tiling in $\zero(P^{\nu, \overline{\Omega}}_{\lambda \mu})$, and thus the following proposition follows. 

\begin{proposition}\label{prop:refined-more}
    Let $P^{\nu, \overline{\Omega}}_{\lambda \mu}$ be the refined puzzle of $P^{\nu, \Omega}_{\lambda \mu}$ w.r.t. an atomic refinement $\overline{\Omega}$ of $\Omega$. Then $\#\zero(P^{\nu, \Omega}_{\lambda \mu}) \leq \#\zero(P^{\nu, \overline{\Omega}}_{\lambda \mu})$.
\end{proposition}

\subsection{Atomic puzzle ideal}
\label{sec:atomic}

Let $P^{\nu, \Omega}_{\lambda \mu}$ be a puzzle of size $n$ and $P^{\nu, \overline{\Omega}}_{\lambda \mu}$ be the refined puzzle of $P^{\nu, \Omega}_{\lambda \mu}$ w.r.t. $\overline{\Omega}$, an atomic refinement of $\Omega$. Next we construct the polynomials to define the atomic puzzle ideal as in Section~\ref{sec:kt}. First we formally define the distinguishing polynomials and prove their existence. 

Let $\Phi$ be the set of all the possible atomic puzzle pieces, with $\Phi_{\bigtriangleup}$ and $\Phi_{\bigtriangledown}$ being its subsets of all the upward and downward ones respectively. Then $\#\Phi_{\bigtriangleup} = \#\Phi_{\bigtriangledown} = 9$. Clearly an upward atomic puzzle piece $p \in \Phi_{\bigtriangleup}$ is determined by the assigned values $\overline{x}$, $\overline{y}$, and $\overline{z}$ of its left, right, and bottom sides. Abusing the notations, we use $(\overline{x}, \overline{y}, \overline{z}) \in \fnum_3^3$ to denote both the vector of assigned values of $p$ and $p$ itself.

\begin{definition}\rm
Let $\overline{\Omega}$ be a set of atomic puzzle pieces. Then a polynomial $f(x, y, z) \in \fnum_3[x, y, z] / \bases{x^3-x, y^3-y, z^3-z}$ is called an \emph{upward distinguishing polynomial} w.r.t.\, $\overline{\Omega}$ if for each $(\overline{x}, \overline{y}, \overline{z}) \in \Phi_{\bigtriangleup}$, we have $f(\overline{x}, \overline{y}, \overline{z}) = 0$ if and only if $(\overline{x}, \overline{y}, \overline{z}) \in \overline{\Omega}$.
\end{definition}

Taking the quotient ring of $\fnum_3[x, y, z]$ modulo $\bases{x^3-x, y^3-y, z^3-z}$ restricts the polynomial $f(x, y, z)$ to take only $\fnum_3$-values. 

\begin{proposition}
    Let $\overline{\Omega}$ be an arbitrary set of atomic puzzle pieces. Then there exists an upward distinguishing polynomial w.r.t.\, $\overline{\Omega}$.
\end{proposition}

\begin{proof}
    A general polynomial in $\fnum_3[x, y, z] / \bases{x^3-x, y^3-y, z^3-z}$ is of the form $f(x, y, z) = \sum_{0\leq d_1, d_2, d_3 \leq 2} c_{d_1d_2d_3}x^{d_1}y^{d_2}z^{d_3}$ and thus has $27$ terms, where $d_1, d_2$, and $d_3$ are all integers. Set $f(\overline{x}, \overline{y}, \overline{z}) = 0$ if $(\overline{x}, \overline{y}, \overline{z}) \in \Phi_{\bigtriangleup}\cap \overline{\Omega}$ and set $f(\overline{x}, \overline{y}, \overline{z}) = 1$ if $(\overline{x}, \overline{y}, \overline{z}) \in \Phi_{\bigtriangleup} \setminus \overline{\Omega}$. These conditions form a system of 9 linear equations over $\fnum_3$ w.r.t. 27 variables $c_{d_1d_2d_3}$ for $0\leq d_1, d_2, d_3 \leq 2$. Then the solution set of this underdetermined linear system is non-empty, and any solution will furnish an upward distinguishing polynomial w.r.t.\, $\overline{\Omega}$. 
\end{proof}

For each downward atomic puzzle piece in $\Phi_{\bigtriangledown}$, it is determined by the assigned values $\overline{x}$, $\overline{y}$, and $\overline{z}$ of its left, top, and right sides. In the exact same way, we can define the downward distinguishing polynomial w.r.t. $\overline{\Omega}$ and prove its existence. 

The triangle $\bigtriangleup^{\nu}_{\lambda \mu}$ is an upward equilateral one of side-length of $n$ units and inside $\bigtriangleup^{\nu}_{\lambda \mu}$ there are $n^2$ unit triangles (including $\frac{n(n+1)}{2}$ upward ones and $\frac{n(n-1)}{2}$ downward ones) with $N = \frac{3n(n+1)}{2}$ unit intervals as their sides. We introduce the variables $x_1, \ldots, x_N$ for the values of the $N$ unit intervals. Now we are ready to construct the polynomials in $\fnum_3[\p{x}]$, where $\p{x} = \{x_1, \ldots, x_N\}$, in the following four groups to define the atomic puzzle ideal.

\begin{enumerate}
    \item Let $\pset{F}_1 = \{x_i^3 - x_i:\, i=1, \ldots, N\}$. They represent the field equations to restrict the variable values in $\fnum_3$. 
    
    \item Let $x_{l1}, \ldots, x_{ln}$ be the variables in $\p{x}$ corresponding to the $n$ intervals of the left side of $\bigtriangleup^{\nu}_{\lambda \mu}$ from left to right  and $x_{r1}, \ldots, x_{rn}, x_{b1}, \ldots, x_{bn}$ be those for the $n$ intervals of the right and bottom sides, also from left to right. Then let $\pset{F}_2 =\{x_{lj} - \lambda_j, x_{rj} - \mu_j, x_{bj} - \nu_j :\,j=1, \ldots, n\}$. 
    
    \item For each unit triangle $\bigtriangleup_k~(k=1, \ldots, n^2)$, let $x_{k1}$, $x_{k2}$, and $x_{k3}$ be the three variables in $\p{x}$ corresponding to its three sides. Then let $\pset{F}_3 = \{x_{k1} + x_{k2} + x_{k3}:\, k=1, \ldots, n^2\}$. 
    
    \item Let $f_{\bigtriangleup}(x,y,z)$ and $f_{\bigtriangledown}(x,y,z)$ be an upward and downward distinguishing polynomial for $\overline{\Omega}$ respectively. For each upward unit triangle $\bigtriangleup_p~(p=1, \ldots, \frac{n(n+1)}{2})$ inside $\bigtriangleup^{\nu}_{\lambda \mu}$, let $x_{p1}$, $x_{p2}$, and $x_{p3}$ be the three variables in $\p{x}$ corresponding to its left, right, and bottom sides, and for each downward $\bigtriangledown_q~(q=1, \ldots, \frac{n(n-1)}{2})$, let $x_{q1}$, $x_{q2}$, and $x_{q3}$ be the three variables corresponding to its left, top, and right sides. Then let
    \begin{equation*}
        \begin{split}
    \pset{F}_4 =& \{f_{\bigtriangleup}(x_{p1}, x_{p2}, x_{p3}):\, p=1, \ldots, \frac{n(n+1)}{2}\} \\
    & \cup \{f_{\bigtriangledown}(x_{q1}, x_{q2}, x_{q3}):\, q=1, \ldots, \frac{n(n-1)}{2}\}.            
        \end{split}
    \end{equation*}
\end{enumerate}

\begin{remark}\rm
  Clearly with the linear polynomials in $\pset{F}_2$, we can remove the polynomials in $\pset{F}_1$ for the variables corresponding to the unit intervals of the three sides of $\bigtriangleup^{\nu}_{\lambda \mu}$. But for the consistency of presentation, we decide to keep $\pset{F}_1$ as it is. Furthermore, from the definition of $\pset{F}_4$ we can find that the upward and downward distinguishing polynomials are indeed templates we apply to each triangle inside $\bigtriangleup^{\nu}_{\lambda \mu}$. 
\end{remark}

\begin{definition}\rm
    Let $P^{\nu, \overline{\Omega}}_{\lambda \mu}$ be a refined puzzle of $P^{\nu, \Omega}_{\lambda \mu}$ and $\pset{F}_1$, $\pset{F}_2$, $\pset{F}_3$, and $\pset{F}_4$ be as defined above with a pair of fixed upward and downward distinguishing polynomials. Then the ideal in $\fnum_3[\p{x}]$ generated by $\bigcup_{i=1, \ldots, 4}\pset{F}_i$ is called the \emph{atomic puzzle ideal} of $P^{\nu, \overline{\Omega}}_{\lambda \mu}$ and denoted by $I^{\nu, \overline{\Omega}}_{\lambda \mu}$.
\end{definition}

Since the defining polynomials in $\pset{F}_4$ are dependent on the choice of upward and downward distinguishing polynomials, the atomic puzzle ideal $I^{\nu, \overline{\Omega}}_{\lambda \mu}$ is not unique. But in fact any atomic puzzle ideal has the desired property in the following theorem, and we decide to be a bit flexible to remove the dependency of atomic puzzle ideals on the choices of upward and downward distinguishing polynomials. Recall that $\var(I^{\nu, \overline{\Omega}}_{\lambda \mu})$ denotes the $\fnum_3$-variety of the ideal $I^{\nu, \overline{\Omega}}_{\lambda \mu}$. 

\begin{theorem}\label{thm:11-atomic}
    Let $P^{\nu, \overline{\Omega}}_{\lambda \mu}$ be a refined puzzle of $P^{\nu, \Omega}_{\lambda \mu}$ and $I^{\nu, \overline{\Omega}}_{\lambda \mu}$ be the atomic puzzle ideal of $P^{\nu, \overline{\Omega}}_{\lambda \mu}$. Then there is a one-one correspondence between $\var(I^{\nu, \overline{\Omega}}_{\lambda \mu})$ and $\zero(P^{\nu, \overline{\Omega}}_{\lambda \mu})$. 
\end{theorem}

\begin{proof}
 Any tiling in $\zero(P^{\nu, \overline{\Omega}}_{\lambda \mu})$ is tiling of $\bigtriangleup^{\nu}_{\lambda \mu}$ with the set $\overline{\Omega}$ of atomic puzzle pieces, and it can be regarded equivalently as assigning values to the $N$ unit intervals in $\bigtriangleup^{\nu}_{\lambda \mu}$ such that (1) the assigned values are in $\fnum_3$; (2) the assigned values to the unit intervals of the left, right, and bottom sides of $\bigtriangleup^{\nu}_{\lambda \mu}$ match $\lambda, \mu, \nu$; (3) the assigned values make each unit triangle in $\bigtriangleup^{\nu}_{\lambda \mu}$ an atomic puzzle piece; and (4) the atomic puzzle pieces constructed in (3) are from $\overline{\Omega}$. Comparing the 4 conditions above with $\pset{F}_1, \ldots, \pset{F}_4$, we know that the assigned values $\overline{x}_1, \ldots, \overline{x}_N$ to $x_1, \ldots, x_N$ satisfying the 4 conditions form a vanishing zero $(\overline{x}_1, \ldots, \overline{x}_N)$ of the ideal $\bases{\cup_{i=1, \ldots, 4}\pset{F}_i}$, and vice versa. This proves the existence of the one-one correspondence. 
\end{proof}

\subsection{Forbidden and implicit puzzle pieces}
From now on, we restrict ourselves to study a refinable set $\Omega$ of puzzle pieces which are of the shapes of unit triangle, unit rhombus, and a bigger polygon (rotation and mirroring are not allowed). Such an $\Omega$ covers the puzzle pieces of all the six puzzles mentioned above. Denote by $\Omega_t$, $\Omega_r$, and $\Omega_b$ the subsets of the puzzle pieces in $\Omega$ of the three shapes above respectively. In the sequel, when we write $\Omega = \Omega_t \cup \Omega_r \cup \Omega_b$, we mean that the puzzle pieces in $\Omega$ are of the three shapes above with $\Omega_t$, $\Omega_r$, and $\Omega_b$ being the corresponding subsets. Clearly for each puzzle piece $p \in \Omega_t \cup \Omega_r$, its atomic refinement is unique and atomic refinement of each $p \in \Omega_b$ is not necessarily unique. 

Let $\overline{\Omega}$ be an atomic refinement of $\Omega = \Omega_t \cup \Omega_r \cup \Omega_b$ induced by $\{\overline{t}_p: p \in \Omega\}$. As shown in Section~\ref{sec:t}, the number of points in $\var(I^{\nu, \overline{\Omega}}_{\lambda \mu})$ may be strictly greater than that of the tilings in $\zero(P^{\nu, \Omega}_{\lambda \mu})$ of our original puzzle. That means, in order to have an ideal whose variety one-one corresponds to the tilings in $\zero(P^{\nu, \Omega}_{\lambda \mu})$, we need to further add more algebraic constraints in the defining polynomials. These constraints indeed come from the stitching process, in which we want to recover a tiling of $\bigtriangleup^{\nu}_{\lambda \mu}$ with $\overline{\Omega}$ back to a tiling with $\Omega$. 

Out of the 18 possible atomic puzzle pieces, the 12 ones which are not $(0, 0, 0)$-, $(1, 1, 1)$-, and $(2, 2, 2)$-valued are said to be \emph{regular}. Clearly any pair of regular atomic puzzle pieces sharing a $2$-side form a rhombus with four sides of $\fnum_2$-values. 

\begin{definition}\rm
Let $\overline{\Omega}$ be an atomic refinement of $\Omega$, $p$ be $\bigtriangleup^{\nu}_{\lambda \mu}$ or a puzzle piece in $\Omega$, and $\overline{t}$ be an arbitrary tiling of $p$ with $\overline{\Omega}$. Then the \emph{stitching} of $\overline{t}$ is the new tiling of $p$ obtained by replacing any pair of regular atomic puzzle pieces sharing a $2$-side with their corresponding $\fnum_2$-valued rhombus.
\end{definition}
From the definition, it is easy to see that the stitching of any tiling of $p$ with $\overline{\Omega}$ is unique and that any remaining $2$-side in a stitching is from a $(2, 2, 2)$-valued atomic puzzle piece. 

In Section~\ref{sec:k} we have seen an example of an implicit rhombus puzzle pieces which is hidden in the puzzle piece in $\Omega_b$. Now we give a formal definition of such pieces. Recall that $\overline{\Omega}$ is an atomic refinement of $\Omega$ induced by $\{\overline{t}_p: p \in \Omega\}$. For each $p \in \Omega_b$, let $\tilde{t}_p$ be the stitching of $\overline{t}_p$. Then the puzzle pieces in $\Omega(\tilde{t}_p) \setminus \Omega$ is said to be \emph{implicit} w.r.t. $\overline{t}_p$. Denote by $\Omega_i(\overline{t}_p)$ the set of implicit puzzle pieces of $p$, and then the puzzle pieces in $\cup_{p \in \Omega}\Omega_i(\overline{t}_p)$ is said to be \emph{implicit} w.r.t. $\overline{\Omega}$. For the particular $(2, 2, 2)$-valued atomic puzzle piece, it is always implicit w.r.t. $\overline{t}_p$ as long as it appears in the tiling $\overline{t}_p$. By definition, an implicit puzzle piece is either the $(2, 2, 2)$-valued atomic one or an $\fnum_2$-valued rhombus one.

\begin{definition}\rm
    Let $\overline{\Omega}$ be an atomic refinement of $\Omega = \Omega_t \cup \Omega_r \cup \Omega_b$ induced by $\{\overline{t}_p:\, p \in \Omega\}$ and $\tilde{t}_p$ be the stitching of $\overline{t}_p$. Then $\overline{\Omega}$ is said to be \emph{separable} if the following two conditions hold:
    \begin{enumerate}
        \item[(1)] for any implicit puzzle piece $\tilde{p}$ w.r.t. $\overline{t}_p$, $\tilde{p}$ occurs only once in $\tilde{t}_p$;
        \item[(2)] for any distinct $p, q \in \Omega_b$, $\Omega_i(\overline{t}_p) \cap \Omega_i(\overline{t}_q) = \emptyset$. 
    \end{enumerate}
\end{definition}

All the six kinds of puzzles above have separable atomic refinements: condition~(2) is trivial because for them $\#\Omega_b = 0$ or $1$, and one can verify that for $\Omega_A$ and $\Omega_B$ whose atomic refinement is unique and for $\Omega_C$ and $\Omega_D$ who have $3$ atomic refinements, condition~(1) also holds. In general we allow an $\Omega$ with $\Omega_b$ containing several puzzle pieces in the same shape of a bigger polygon as long as its atomic refinement is separable. 

As shown in Section~\ref{sec:k}, implicit puzzle pieces can appear in the stitching of some tiling $\overline{t} \in \zero(P^{\nu, \overline{\Omega}}_{\lambda \mu})$ and if we want to recover $\overline{t}$ to some tiling $t \in \zero(P^{\nu, \Omega}_{\lambda \mu})$, then their appearance needs to imply that of the original puzzle pieces in $\Omega_b$ corresponding to them. The conditions of $\overline{\Omega}$ being separable ensure that once an implicit puzzle piece $\tilde{p}$ appears, there is only one puzzle piece $p$ in $\Omega_b$ corresponding to it (condition~(2)) and the relative position of $\tilde{p}$ in $p$ is also determined (condition~(1)), so that we can generalize the method based on implying polynomials to deal with them as in Section~~\ref{sec:k}.

For the atomic refinement $\overline{\Omega}$ of $\Omega = \Omega_t \cup \Omega_r \cup \Omega_b$ induced by $\{\overline{t}_p: p \in \Omega\}$, let $\overline{\Omega}_i$ be the set of all the implicit puzzle pieces w.r.t. $\overline{\Omega}$. Consider the set 
\begin{equation}\label{eq:Omega_r}
\begin{split}
\overline{\Omega}_r = \{\fnum_2\mbox{-valued rhombus puzzle piece } \tilde{p}: \tilde{p} & \mbox{ has a tiling with } \overline{\Omega} \\
&\mbox{ but not with } \Omega\}.    
\end{split}
\end{equation}
Then the puzzle pieces in $\overline{\Omega}_f := \overline{\Omega}_r \setminus \overline{\Omega}_i$ are said to be \emph{forbidden} w.r.t. $\overline{\Omega}$. By definition, forbidden puzzle pieces are all $\fnum_2$-valued rhombus ones. 

Let $\overline{\Omega}$ be a separable atomic refinement of $\Omega = \Omega_t \cup \Omega_r \cup \Omega_b$ induced by $\{\overline{t}_p: p \in \Omega\}$. For any tiling $\overline{t} \in \zero(P^{\nu, \overline{\Omega}}_{\lambda \mu})$, denote the stitching of $\overline{t}$ by $\tilde{t}$, which represents our first step to transform $\overline{t}$ to a tiling in $\zero(P^{\nu, \Omega}_{\lambda \mu})$. Note that by definition, the set $\overline{\Omega}_r$ above covers all the possible $\fnum_2$-valued rhombus puzzle pieces in $\tilde{t}$ which can be tiled with $\overline{\Omega}$ but cannot with $\Omega$. The $\fnum_2$-valued rhombus puzzle pieces in $\overline{\Omega}_r$ are divided into two groups. The subset $\overline{\Omega}_f$ covers those whose appearance will directly disable the transformation, for these pieces do not have tilings with $\Omega$ and cannot be expanded to some bigger puzzle pieces in $\Omega_b$. Or in other words, if we forbid the appearance of all the puzzle pieces in $\overline{\Omega}_f$ in $\tilde{t}$, then $\tilde{t}$ will only have original puzzle pieces in $\Omega_t \cup \Omega_r$ and implicit ones w.r.t. $\overline{\Omega}$, which can all be expanded to bigger bigger puzzle pieces in $\Omega_b$. With $\overline{\Omega}$ separable, each implicit puzzle piece $\tilde{p} \in \overline{\Omega}_i$ corresponds to only one piece $p \in \Omega_b$ and its position in $\tilde{t}_p$ is determined, where $\tilde{t}_p$ is the stitching of $\overline{t}_p$. And if we further restrict the tiling $\tilde{t}$ such that the appearance of each $\tilde{p}$ will force the puzzle pieces around $\tilde{p}$ in $\tilde{t}$ to form $\tilde{t}_p$, then replacing $\tilde{t}_p$ in $\tilde{t}$ with $p$ will result in a tiling in $\zero(P^{\nu, \Omega}_{\lambda \mu})$. 

For all the forbidden puzzle pieces in $\overline{\Omega}_f$, we forbid their appearance in $\tilde{t}$ by ensuring that the assigned values to the edges of all the rhombuses in $\bigtriangleup^{\nu}_{\lambda \mu}$ will not turn them into forbidden puzzle pieces, and this is achieved by associating a forbidding polynomial to each rhombus in $\bigtriangleup^{\nu}_{\lambda \mu}$. For each implicit puzzle piece $\tilde{p}$ appearing in $\tilde{t}$, we ensure the puzzle pieces around $\tilde{p}$ to form the tiling $\tilde{t}_p$ of the corresponding $p \in \Omega_b$ by restricting the assigned values to the edges of the surrounding puzzle pieces to match those of $p$, and this is achieved by associating an implying polynomial to each polygon  in $\bigtriangleup^{\nu}_{\lambda \mu}$ of the shape of $p$.

We first define forbidding polynomials and prove their existence similarly to what we have done for distinguishing polynomials. Let 
$$\Psi = \{\fnum_3\mbox{-valued rhombus puzzle piece } \tilde{p}: \mbox{ there exists a tiling of } \tilde{p} \mbox{ with } \overline{\Omega}\}$$
and $\Psi_l$, $\Psi_r$, and $\Psi_b$ being its subsets of all the rhombus puzzle pieces whose middle unit intervals are parallel to the left, right, and bottom side of $\bigtriangleup^{\nu}_{\lambda \mu}$ respectively. Then one can check that in the extreme case when $\overline{\Omega}$ contains all the 18 possible atomic puzzle pieces, we have $\#\Psi_l = \#\Psi_r = \#\Psi_b = 27$, and thus for general $\overline{\Omega}$ the three numbers are $\leq 27$. For any rhombus puzzle piece $\overline{p} \in \Psi_l$, let $\overline{x}_{1}$ be the $\fnum_2$-value assigned to 
its top side and $\overline{x}_2$, $\overline{x}_3$, and $\overline{x}_4$ be those assigned to the other 3 sides in a clockwise direction. Then $\overline{p}$ is determined by the four assigned values $(\overline{x}_1, \overline{x}_2, \overline{x}_3, \overline{x}_4) \in \fnum_3^4$. Abusing the notations, we use $(\overline{x}_1, \overline{x}_2, \overline{x}_3, \overline{x}_4) \in \fnum_3^4$ to denote both the vector of assigned values of $\overline{p}$ and $\overline{p}$ itself. 

\begin{definition}\rm
    Let $\overline{\Omega}$ be an atomic refinement of $\Omega = \Omega_t \cup \Omega_r \cup \Omega_b$ and $\overline{\Omega}_f$ be the set of all the forbidden puzzle pieces w.r.t. $\overline{\Omega}$. Then a polynomial $f_l(z_1, z_2, z_3, z_4) \in \fnum_3[z_1, z_2, z_3, z_4] / \bases{z_1^3-z_1, z_2^3-z_2, z_3^3-z_3, z_4^3-z_4}$ is called a \emph{left forbidding polynomial} w.r.t.\, $\overline{\Omega}$ if for each $(\overline{x}_1, \overline{x}_2, \overline{x}_3, \overline{x}_4) \in \Psi_l$, we have $f_l(\overline{x}_1, \overline{x}_2, \overline{x}_3, \overline{x}_4) \neq 0$ if and only if $(\overline{x}_1, \overline{x}_2, \overline{x}_3, \overline{x}_4) \in \overline{\Omega}_f$. 
\end{definition}

\begin{proposition}
    Let $\overline{\Omega}$ be an atomic refinement of $\Omega = \Omega_t \cup \Omega_r \cup \Omega_b$ and $\overline{\Omega}_f$ be the set of all the forbidden puzzle pieces w.r.t. $\overline{\Omega}$. Then there exists a left forbidding polynomial w.r.t.\, $\overline{\Omega}$.
\end{proposition}

\begin{proof}
    Comparing the number of coefficients 81 in the general form of $f(z_1, z_2, z_3, z_4)$ and the maximum number of linear equations 27 from the constrains $f(\overline{x}_1, \overline{x}_2, \overline{x}_3, \overline{x}_4) = 0$ or $1$ for $(\overline{x}_1, \overline{x}_2, \overline{x}_3, \overline{x}_4) \in \Psi_l$, we know that the solution set of this underdetermined linear system is non-empty. 
\end{proof}

In the same way we can define the right and bottom forbidding polynomials $f_r(z_1, z_2, z_3, z_4)$ and $f_b(z_1, z_2, z_3, z_4)$ w.r.t. $\overline{\Omega}$ and prove their existence. These forbidding polynomials allow $\fnum_3$-valued rhombus puzzle pieces in $\Psi \setminus \overline{\Omega}_f$, according to their directions, to appear and disallow the appearance of forbidden rhombus puzzle pieces in $\overline{\Omega}_f$. 

To deal with implicit puzzle pieces we need the atomic refinement $\overline{\Omega}$ of $\Omega = \Omega_t \cup \Omega_r \cup \Omega_b$ induced by $\{\overline{t}_p:\, p \in \Omega\}$ to be separable. Let $\overline{\Omega}_i$ be the set of all the implicit puzzle pieces w.r.t. $\overline{\Omega}$. Then for each $\overline{p} \in \overline{\Omega}_i$, with $\overline{\Omega}$ separable, we know that there exists a unique $p \in \Omega_b$ such that $\overline{p}$ is implicit w.r.t. $\overline{t}_p$. Consider the convex polygon of the shape of the puzzle pieces in $\Omega_b$, and assume that there are $m$ unit intervals inside the polygon. Then the tiling $\overline{t}_p$ of $p$ with atomic puzzle pieces essentially assigns $\fnum_3$-values to all the $m$ intervals of the polygon. Let the implicit puzzle piece $\overline{p}$ be determined by the assigned values $\overline{x}_1, \ldots, \overline{x}_i$, where $i=3$ if $\overline{p}$ is the $(2, 2, 2)$-valued atomic puzzle piece and $i=4$ if $\overline{p}$ is an $\fnum_2$-valued rhombus puzzle piece. Then the condition that appearance of $\overline{p}$ in the stitching of a tiling $\overline{t} \in \zero(P^{\nu, \overline{\Omega}}_{\lambda \mu})$ implies the puzzle pieces around $\overline{p}$ to form a tiling of the corresponding $p \in \Omega_b$ means that the assigned values to the surrounding puzzle pieces match those assigned by $\overline{t}_p$ as long as the values $\overline{x}_1, \ldots, \overline{x}_i$ match those assigned to the rhombus by $\overline{t}_p$. This observation motivates us to study the following implying polynomial for each implicit puzzle piece.

Let the puzzle pieces in $\Omega_b$ be of the shape of a convex polygon with $m$ unit intervals inside. Label all the $m$ unit intervals of the polygon by $I_1, \ldots, I_m$. Then the $\fnum_3$-values $\overline{x}_1, \ldots, \overline{x}_m$ assigned to $I_1, \ldots, I_m$ respectively determine an $\fnum_3$-valued puzzle piece. For each implicit puzzle piece $\overline{p}$ which corresponds uniquely to $p \in \Omega_b$, let $\tilde{t}_p$ be the stitching of $\overline{t}_p$. We take a puzzle piece $\overline{p}$ in $\Psi_l$ for example. With $\overline{\Omega}$ induced by $\{\overline{t}_p:\, p \in \Omega\}$ separable, we know that there is only one occurrence of $\overline{p}$ in $\tilde{t}_p$. Let $I_{k_1}$ be the interval in $I_1, \ldots, I_m$ the top side of $\overline{p}$ takes, and $I_{k_2}$, $I_{k_3}$, and $I_{k_4}$ be the intervals the other sides take clockwisely. Then for a vector $\overline{\p{x}} = (\overline{x}_1, \ldots, \overline{x}_m) \in \fnum_3^m$, we use $\overline{\p{x}}^{\ominus}$ to denote the $(m-4)$-long vector obtained by removing $\overline{x}_{k_1}$, $\overline{x}_{k_2}$, $\overline{x}_{k_3}$, and $\overline{x}_{k_4}$. Also for an $(m-4)$-long vector $\overline{\p{y}}_1 = (\overline{y}_{11}, \ldots, \overline{y}_{1,m-4}) \in \fnum_3^{m-4}$ and a $4$-long vector $\overline{\p{y}}_2 = (\overline{y}_{21}, \ldots, \overline{y}_{24}) \in \fnum_3^4$, we use $\overline{\p{y}}_1 \oplus \overline{\p{y}}_2$ to denote the $m$-long vector obtained by inserting $\overline{y}_{2,1}, \ldots, \overline{y}_{2,4}$ in $\overline{\p{y}}_1$ so that $\overline{y}_{2,j}$ occurs at the $k_j$-th position of the resulting vector for $j=1, \ldots, 4$.

\begin{definition}\rm
    Let $\overline{\Omega}$ be a separable atomic refinement of $\Omega = \Omega_t \cup \Omega_r \cup \Omega_b$ induced by $\{\overline{t}_p:\, p \in \Omega\}$ and $\overline{p} = (\overline{y}_1, \ldots, \overline{y}_4) \in \fnum_2^4$ be an implicit puzzle piece w.r.t. $\overline{\Omega}$ which is in $\Psi_l$ and corresponds to $p \in \Omega_b$. Let $I_1, \ldots, I_m$ be the $m$ unit intervals inside $p$, $I_{k_1}, \ldots, I_{k_4}$ being the intervals $\overline{p}$ takes, and $\overline{\p{x}} = (\overline{x}_1, \ldots, \overline{x}_m) \in \fnum_3^m$ be the assigned values of $\overline{t}_p$ to $I_1, \ldots, I_m$. Then a polynomial $f_{\overline{p}}(z_1, \ldots, z_m) \in \fnum_3[z_1, \ldots, z_m] / \bases{z_1^3-z_1, \ldots, z_m^3-z_m}$ is called an \emph{implying polynomial} of $\overline{p}$ w.r.t. $\overline{\Omega}$ if $f(\overline{\p{y}} \oplus \overline{p}) = 0$ when $\overline{\p{y}} = \overline{\p{x}}^{\ominus}$ and $f(\overline{\p{y}} \oplus \overline{p}) = 1$ for any $\overline{\p{y}} \in \fnum_3^{m-4} \setminus \{\overline{\p{x}}^{\ominus}\}$. 
\end{definition}

\begin{proposition}\label{prop:implying}
    Let $\overline{\Omega}$ be a separable atomic refinement of $\Omega = \Omega_t \cup \Omega_r \cup \Omega_b$ induced by $\{\overline{t}_p:\, p \in \Omega\}$ and $\overline{p}$ be any implicit puzzle piece w.r.t. $\overline{\Omega}$ in $\Psi_l$. Then there exists an implying polynomial of $\overline{p}$ w.r.t.\, $\overline{\Omega}$.
\end{proposition}

\begin{proof}
This can be similarly proved by considering the corresponding underdetermined linear system. 
\end{proof}

Implying polynomials for implicit puzzle pieces in $\Psi_r$ and $\Psi_b$ and their existence can be studied in the exactly same way. For that for the $(2, 2, 2)$-valued atomic puzzle piece, the study is similar. 

For each implicit puzzle piece $\overline{p}_1$ w.r.t. $\overline{\Omega}$, there is only one $p \in \Omega_b$ corresponding to it if $\overline{\Omega}$ is a separable, but there may be another implicit puzzle piece $\overline{p}_2$ which corresponds to $p$ too. Then we will have the implying polynomials $f_{1}$ for $\overline{p}_1$ and $f_{2}$ for $\overline{p}_2$ respectively. With some further investigation on the definition of implying polynomials, one will see that $f_1(\overline{\p{x}}) = f_2(\overline{\p{x}}) = 0$ for the assigned values $\overline{\p{x}}$ of $\overline{t}_p$ and that $f_1(\overline{\p{z}}) = f_1(\overline{\p{z}})$ for any vector $\overline{\p{z}} \in \fnum_3^m$ which can be taken for both $\overline{p}_1$ and $\overline{p}_2$. This means that our definition of implying polynomials does not introduce inconsistent polynomial equations for multiple implicit puzzle pieces inside the same piece in $\Omega_b$. 

\subsection{Puzzle ideal}
Now we are ready to adjoin the forbidding and implying polynomials to the defining polynomials of the ideals. Let $\overline{\Omega}$ be a separable atomic refinement of $\Omega = \Omega_t \cup \Omega_r \cup \Omega_b$ induced by $\{\overline{t}_p:\, p \in \Omega\}$, and $\overline{\Omega}_f$ and $\overline{\Omega}_i$ be the sets of forbidden and implicit puzzle pieces w.r.t. $\overline{\Omega}$ respectively. Let $f_l(z_1, z_2, z_3, z_4)$, $f_r(z_1, z_2, z_3, z_4)$, and $f_b(z_1, z_2, z_3, z_4)$ be a left, right, and bottom forbidding polynomial w.r.t. $\overline{\Omega}$ respectively, and for each $\overline{p} \in \overline{\Omega}_i$, let $f_{\overline{p}}(z_1, \ldots, z_{m})$ be the implying polynomial of $\overline{p}$ w.r.t. $\Omega$, where $m$ is the number of unit intervals inside the polygon which is the shape of puzzle pieces in $\Omega_b$. 

For the triangle $\bigtriangleup^{\nu}_{\lambda \mu}$ of size $n$, there are $\frac{n(n-1)}{2}$ left rhombuses ~\lrhombus~, right ones ~~\rrhombus, and bottom ones $\lozenge$ inside $\bigtriangleup^{\nu}_{\lambda \mu}$ respectively. For each category of rhombuses inside $\bigtriangleup^{\nu}_{\lambda \mu}$, order them row by row from top to bottom and from left to right. Then for the $k$-th left rhombus ~\lrhombus~ and $k$-th right one ~~\rrhombus, let $x_{l, k1}, \ldots, x_{l, k4}$ and $x_{r, k1}, \ldots, x_{r, k4}$ be the variables in $\p{x}$ corresponding to the top, right, bottom, and left sides of the rhombus respectively. And for the $k$-th bottom rhombus $\lozenge$, let $x_{b, k1}, \ldots, x_{b, k4}$ be variables for the NW, NE, SE, and SW sides. Then let
\begin{equation*}
\begin{split}
    \pset{F}_5 = \{f_l(x_{l, k1}, \ldots, x_{l, k4}), f_r(x_{r, k1}, \ldots, x_{r, k4}), f_b(x_{b, k1}, \ldots, x_{b, k4}) :&\\
    k = 1, \ldots, \frac{n(n-1)}{2}\}&.
\end{split}
\end{equation*}
They are the forbidding polynomials applied to each rhombus inside $\bigtriangleup^{\nu}_{\lambda \mu}$ and $\#\pset{F}_5 = \frac{3n(n-1)}{2}$.

Assume that the puzzle pieces in $\Omega_b$ are of the shape of a convex polygon and there are $M$ such polygons inside $\bigtriangleup^{\nu}_{\lambda \mu}$. We also order these polygons row by row from top to bottom and from left to right. Then for the $k$-th polygon, let $x_{k1}, \ldots, x_{km}$ be the variables in $\p{x}$ corresponding to the $m$ intervals. Then let 
$$\pset{F} _6 = \{f_{\overline{p}}(x_{k1}, \ldots, x_{km}): \overline{p} \in \overline{\Omega}_i, k= 1, \ldots, M\}.$$
They are the implying polynomials for all the implicit puzzle pieces applied to each polygon inside $\bigtriangleup^{\nu}_{\lambda \mu}$ and $\#\pset{F}_6 = M\#\overline{\Omega}_i$.

\begin{definition}\rm
Let $P^{\nu, \Omega}_{\lambda \mu}$ be a puzzle with $\Omega = \Omega_t \cup \Omega_r \cup \Omega_b$, $\overline{\Omega}$ be a separable atomic refinement of $\Omega$ induced by $\{\overline{t}_p:\, p \in \Omega\}$, and $\pset{F}_1, \ldots, \pset{F}_6$ be as defined above for fixed choices of distinguishing, forbidding, and implying polynomials. Then the ideal in $\fnum_3[\p{x}]$ generated by $\bigcup_{i=1, \ldots, 6}\pset{F}_i$ is called the \emph{puzzle ideal} of $P^{\nu, \Omega}_{\lambda \mu}$ and denoted by $I^{\nu, \Omega}_{\lambda \mu}$.
\end{definition}

Next we prove the one-one correspondence between the variety $\var(I^{\nu, \Omega}_{\lambda \mu})$ and $\zero(P^{\nu, \Omega}_{\lambda \mu})$, which justifies the name of puzzle ideals. 

\begin{theorem}\label{thm:11-S}
    Let $P^{\nu, \Omega}_{\lambda \mu}$ be a puzzle with $\Omega = \Omega_t \cup \Omega_r \cup \Omega_b$, $\overline{\Omega}$ be a separable atomic refinement of $\Omega$ induced by $\{\overline{t}_p:\, p \in \Omega\}$, and $I^{\nu, \Omega}_{\lambda \mu}$ be the puzzle ideal of $P^{\nu, \Omega}_{\lambda \mu}$. For any $\overline{t} \in \zero(P^{\nu, \overline{\Omega}}_{\lambda \mu})$, let $\tilde{t}$ be its stitching and for any $\tilde{p} \in \Omega(\tilde{t}) \cap \overline{\Omega}_i$, let $p$ be the puzzle piece in $\Omega_b$ corresponding to $\tilde{p}$. Denote
    \begin{equation*}
        \begin{split}
    \pset{S} = \{\overline{t} \in \zero(P^{\nu, \overline{\Omega}}_{\lambda \mu}):\, & \Omega(\tilde{t}) \cap \overline{\Omega}_f = \emptyset \mbox{ and for any }\tilde{p} \in \Omega(\tilde{t}) \cap \overline{\Omega}_i, \tilde{p} \mbox{ and its} \\
    & \mbox{ surrounding puzzle pieces form the stitching of } \overline{t}_p\}.
        \end{split}
    \end{equation*}
Then there is a one-one correspondence between $\var(I^{\nu, \Omega}_{\lambda \mu})$ and $\pset{S}$. 
\end{theorem}

\begin{proof}
    Let $\psi: \var(\leftindex_a I^{\nu, \overline{\Omega}}_{\lambda \mu}) \rightarrow \zero(P^{\nu, \overline{\Omega}}_{\lambda \mu})$ be the one-one correspondence in Theorem~\ref{thm:11-atomic}. With $\var(I^{\nu, \Omega}_{\lambda \mu}) \subseteq \var(\leftindex_a I^{\nu, \overline{\Omega}}_{\lambda \mu})$ and $\pset{S} \subseteq \zero(P^{\nu, \overline{\Omega}}_{\lambda \mu})$, next we prove the restriction of $\psi$ on $\var(I^{\nu, \Omega}_{\lambda \mu})$ maps one-one to $\pset{S}$. For any zero $\overline{\p{z}} \in \var(I^{\nu, \Omega}_{\lambda \mu})$, $\psi(\overline{\p{z}}) \in \zero(P^{\nu, \overline{\Omega}}_{\lambda \mu})$. To prove $\psi(\overline{\p{z}}) \in \pset{S}$, it suffices to show that the tiling $\psi(\overline{\p{z}})$ satisfies the two conditions in $\pset{S}$. 

    Recall that $\pset{F}_5 \subseteq I^{\nu, \Omega}_{\lambda \mu}$, and thus $f(\overline{\p{z}}) = 0$ for any $f \in \pset{F}_5$. The polynomials in $\pset{F}_5$ cover all the rhombuses inside $\bigtriangleup^{\nu}_{\lambda \mu}$ and they are the forbidding polynomials applied to the rhombuses. Therefore $f(\overline{\p{z}}) = 0$ is equivalent to that the assigned values of $\overline{\p{z}}$ to the rhombus corresponding to $f$ do not make the rhombus any of the forbidden pieces. The latter condition implies that there is no forbidden piece in the stitching of $\psi(\overline{\p{z}})$: this proves condition (1). 

    Recall that $\pset{F}_6 \subseteq I^{\nu, \Omega}_{\lambda \mu}$, and thus $f(\overline{\p{z}}) = 0$ for any $f \in \pset{F}_6$. For each $\tilde{p} \in \overline{\Omega}_i$, the polynomials in $\pset{F}_6$ for $\tilde{p}$ cover all the polygons inside $\bigtriangleup^{\nu}_{\lambda \mu}$, and they are the implying polynomial of $\tilde{p}$ applied to the polygons. Let $f$ be a polynomial in $\pset{F}_6$ for $\tilde{p}$. Then $f(\overline{\p{z}}) = 0$ is equivalent to that the assigned values of $\overline{\p{z}}$ to the polygon corresponding to $f$ make the tiling $\tilde{t}_p$ of $p \in \Omega_b$ as long as the assigned values of $\overline{\p{z}}$ to the rhombus inside the polygon make it to $\tilde{p}$. The latter condition implies that, if in the stitching of $\psi(\overline{\p{z}})$ there exists some implicit puzzle $\tilde{p}$, then for any occurrence of $\tilde{p}$ in the stitching of $\psi(\overline{\p{z}})$, the polygon in which $\tilde{p}$'s position coincide with the position in $\tilde{t}_p$ becomes the tiling $\overline{t}_p$ of $p$ with the assigned values of $\overline{\p{z}}$, and thus $\tilde{p}$ and its surrounding puzzle pieces form $\tilde{t}_p$. This proves condition (2). 

    Now we prove that the restriction of $\psi$ on $\var(I^{\nu, \Omega}_{\lambda \mu})$ is onto $\pset{S}$. Assume that there is a tiling $\overline{t} \in \pset{S}$ such that no points in $\var(I^{\nu, \Omega}_{\lambda \mu})$ maps to it. With $\overline{t}$ in $\pset{S} \subseteq \zero(P^{\nu, \overline{\Omega}}_{\lambda \mu})$, from the one-one correspondence $\psi$ we know that there exists some points $\overline{\p{z}} \in \var(\leftindex_a I^{\nu, \overline{\Omega}}_{\lambda \mu})$ such that $\overline{\p{z}} = \psi^{-1}(\overline{t})$. But $\overline{\p{z}} \not \in \var(I^{\nu, \Omega}_{\lambda \mu})$, and thus there exists a polynomial $f \in \pset{F}_5 \cup \pset{F}_6$ such that $f(\overline{\p{z}}) \neq 0$. 

    If $f \in \pset{F}_5$, then $f$ is a forbidding polynomial for some rhombus inside $\bigtriangleup^{\nu}_{\lambda \mu}$, and $f(\overline{\p{z}}) \neq 0$ is equivalent to that the assigned values of $\overline{\p{z}}$ to this rhombus makes it a forbidden puzzle piece. This implies that $\Omega(\tilde{t}) \cap \overline{\Omega}_f \neq \emptyset$: a contradiction with condition (1). 
    
    If $f \in \pset{F}_6$, then $f$ is an implying polynomial for some implicit puzzle piece $\tilde{p}$ in $\Omega(\tilde{t})$ applied to a polygon inside $\bigtriangleup^{\nu}_{\lambda \mu}$. In this case $f(\overline{\p{z}}) \neq 0$ is equivalent to that the assigned values of $\overline{\p{z}}$ to the surrounding unit intervals inside the corresponding polygon do not math those of $p$, and thus together with $\tilde{p}$ they do not form $\tilde{t}_p$: a contradiction with condition (2).     
\end{proof}

Theorem~\ref{thm:11-S} above identifies the tilings in $\zero(P^{\nu, \overline{\Omega}}_{\lambda \mu})$ which the variety $\var(I^{\nu, \Omega}_{\lambda \mu})$ corresponds to. Next we show that these tilings in $\pset{S}$ are those which can be recovered to tilings in $\zero(P^{\nu, \Omega}_{\lambda \mu})$.

\begin{lemma}\label{lem:stitching}
Let $P^{\nu, \Omega}_{\lambda \mu}$ be a puzzle with $\Omega = \Omega_t \cup \Omega_r \cup \Omega_b$ and $\overline{\Omega}$ be a separable atomic refinement of $\Omega$. For the stitching $\tilde{t}$ of any $\overline{t} \in \zero(P^{\nu, \overline{\Omega}}_{\lambda \mu})$, let $\Omega(\tilde{t})$ be the set of puzzle pieces appearing in $\tilde{t}$. If $\Omega(\tilde{t})$ does not contain any forbidden puzzle piece, then $\Omega(\tilde{t}) \subseteq \Omega_t \cup \Omega_r \cup \overline{\Omega}_i$, where $\overline{\Omega}_i$ is the set of implicit puzzle pieces w.r.t. $\overline{\Omega}$. 
\end{lemma}

\begin{proof}
   For any $\tilde{p} \in \Omega(\tilde{t})$, we know that $\tilde{p}$ is of the shape of a unit triangle or rhombus. 

   If $\tilde{p}$ is of the shape of a unit triangle, then by the fact that there is no regular atomic puzzle pieces in $\tilde{t}$, we know that $\tilde{p} \in \Omega_t$ or $\tilde{p}$ is the $(2, 2, 2)$-valued atomic one, which is always in $\overline{\Omega}_i$. 

   If $\tilde{p}$ is of the shape of a unit rhombus, then it comes from a pair of upward and downward regular atomic puzzle pieces sharing a $2$-side, which implies that $\tilde{p}$ is an $\fnum_2$-valued rhombus puzzle piece. If $\tilde{p} \in \Omega_r$, then we finish the proof. Otherwise $\tilde{p}$ is in the stitching $\tilde{t}_q$ for some $q \in \Omega_b$, and thus $\tilde{p} \in \overline{\Omega}_r$ in \eqref{eq:Omega_r}. From the assumption we know that $\tilde{p}$ is not a forbidden puzzle piece, which implies that $\tilde{p}$ is an implicit one. 
\end{proof}

\begin{theorem}\label{thm:main}
    Let $P^{\nu, \Omega}_{\lambda \mu}$ be a puzzle with $\Omega = \Omega_t \cup \Omega_r \cup \Omega_b$, $\overline{\Omega}$ be a separable atomic refinement of $\Omega$ induced by $\{\overline{t}_p:\, p \in \Omega\}$, and $I^{\nu, \Omega}_{\lambda \mu}$ be the puzzle ideal of $P^{\nu, \Omega}_{\lambda \mu}$. Then there exists a one-one correspondence between $\var(I^{\nu, \Omega}_{\lambda \mu})$ and $\zero(P^{\nu, \Omega}_{\lambda \mu})$. 
\end{theorem}

\begin{proof}
Let $\pset{S}$ be as defined in Theorem~\ref{thm:11-S}. With the one-one correspondence between $\pset{S}$ and $\var(I^{\nu, \Omega}_{\lambda \mu})$ in that theorem, it suffices to find a one-one correspondence between $\pset{S}$ and $\zero(P^{\nu, \Omega}_{\lambda \mu})$.

For any tiling $t \in \zero(P^{\nu, \Omega}_{\lambda \mu})$, replacing any puzzle piece $p \in \Omega$ with the tiling $\overline{t}_p$ will give a tiling $\overline{t} \in \zero(P^{\nu, \overline{\Omega}}_{\lambda \mu})$. Let $\tilde{t}$ be the stitching of $\overline{t}$. We first prove that there is no forbidden puzzle piece in $\tilde{t}$. Otherwise assume that $\tilde{p}$ is a forbidden puzzle piece in $\tilde{t}$. Note that $\tilde{t}$ can be regarded as the union of the stitching of $\overline{t}_p$ for $p \in \Omega$, and thus there exists $p \in \Omega$ such that $\tilde{p}$ is in the stitching of $\overline{t}_p$. But by definition, this means that $\tilde{p}$ is implicit: a contradiction. This implies that $\overline{t}$ corresponding to $t \in \zero(P^{\nu, \Omega}_{\lambda \mu})$ satisfies condition (1) of $\pset{S}$. For each implicit puzzle piece $\tilde{p} \in \Omega(\tilde{t})$, let $p \in \Omega_b$ be the puzzle piece corresponding to $\tilde{p}$. Then we know that $\tilde{p} \in \Omega(\tilde{t}_p)$, where $\tilde{t}_p$ is the stitching of $\overline{t}_p$ with $\tilde{p}$. Note that $\overline{t}$ is obtained by replacing $p$ in $t$ by $\overline{t}_p$, and thus the surrounding puzzles form $\tilde{t}_p$. This proves condition (2) of $\pset{S}$, and thus the tiling $\overline{t}$ constructed from $t$ belongs to $\pset{S}$.

Let $\overline{t}$ be any tiling in $\pset{S}$ and $\tilde{t}$ be its stitching. Then there is no forbidden puzzle pieces in $\tilde{t}$, and thus by Lemma~\ref{lem:stitching} we know that $\Omega(\tilde{t}) \subseteq \Omega_t \cup \Omega_r \cup \overline{\Omega}_i$. This means that apart from the original puzzle pieces in $\Omega_t \cup \Omega_r$, there are only implicit puzzle pieces w.r.t. $\overline{\Omega}$. For each $\tilde{p} \in \Omega(\tilde{t}) \cap \overline{\Omega}_i$, let $p \in \Omega_b$ be the puzzle piece corresponding to $\tilde{p}$. Then by condition (2) of $\pset{S}$, the surrounding puzzle pieces form the stitching of $\overline{t}_p$ with $\tilde{p}$. Replacing the stitching of $\overline{t}_p$ with $p$ for each $\tilde{p}\in \Omega(\tilde{t})$ will turn the tiling $\tilde{t}$ to a tiling of $\bigtriangleup^{\nu}_{\lambda \mu}$ with puzzle pieces from only $\Omega_t \cup \Omega_r \cup \Omega_b$, which means a tiling in $\zero(P^{\nu, \Omega}_{\lambda \mu})$. This finishes the proof. 
\end{proof}

For any puzzle $P^{\nu, \Omega}_{\lambda \mu}$ whose puzzle pieces are of the shapes of unit triangle, rhombus, and a bigger polygon such that $\Omega$ is refinable with a separable atomic refinement, one can construct the puzzle ideal $I^{\nu, \Omega}_{\lambda \mu}$ for it and its variety $\var(I^{\nu, \Omega}_{\lambda \mu})$ one-one corresponds to $\zero(P^{\nu, \Omega}_{\lambda \mu})$. This one-one correspondence enables us to compute all the tilings in $\zero(P^{\nu, \Omega}_{\lambda \mu})$ by solving the defining polynomial system of $I^{\nu, \Omega}_{\lambda \mu}$, say by using the \grobner basis. The puzzle pieces of all the six puzzles for Grassmannians in Section~\ref{sec:pre} satisfy the conditions above and thus our treatment with the puzzle ideal $I^{\nu, \Omega}_{\lambda \mu}$ works for all of them.

\section{Side-free puzzle ideal}
\label{sec:side-free}

Recall that for $\lambda, \mu, \nu \in \abinom{n}{k}$, $\bigtriangleup^{\nu}_{\lambda \mu}$ is the upright equilateral triangle whose left, right, and bottom sides are assigned the values of $\lambda$, $\mu$, and $\nu$ and $P^{\nu, \Omega}_{\lambda \mu}$ is the puzzle consisting of $\bigtriangleup^{\nu}_{\lambda \mu}$ with a set $\Omega$ of puzzle pieces. Then $\zero(P^{\nu, \Omega}_{\lambda \mu})$ consists of all the tilings of $\bigtriangleup^{\nu}_{\lambda \mu}$ with $\Omega$. For an upward equilateral triangle of side-length of $n$ units, if two of its sides are assigned $\fnum_2$-values and the remaining one is not,  say the left and right sides are assigned values of $\lambda, \mu \in \abinom{n}{k}$ but the bottom one is not, then we denote it by $\bigtriangleup^{\emptyset}_{\lambda \mu}$. The puzzle consisting of $\bigtriangleup^{\emptyset}_{\lambda \mu}$ and the set $\Omega$ of $\fnum_2$-valued puzzle pieces is called the $\nu$-free puzzle w.r.t. $\Omega$ and denoted by $P^{\emptyset, \Omega}_{\lambda \mu}$. Similarly one can define $P^{\nu, \Omega}_{\emptyset \mu}$ and $P^{\nu, \Omega}_{\lambda \emptyset}$, and they are called \emph{side-free puzzles}. In the sequel we work with $P^{\emptyset, \Omega}_{\lambda \mu}$ for example, and the discussions for the other side-free puzzles $P^{\nu, \Omega}_{\emptyset \mu}$ and $P^{\nu, \Omega}_{\lambda \emptyset}$ are the same. 

Similarly, the set of tilings of $\bigtriangleup^{\emptyset}_{\lambda \mu}$ with $\Omega$ is denoted by $\zero(P^{\emptyset, \Omega}_{\lambda \mu})$. Fix a vector $\overline{\nu} \in \fnum_2^n$, we have the puzzle $P^{\overline{\nu}, \Omega}_{\lambda \mu}$ and clearly $\zero(P^{\overline{\nu}, \Omega}_{\lambda \mu}) \subseteq \zero(P^{\emptyset, \Omega}_{\lambda \mu})$. In particular, for $\overline{\nu}_1, \overline{\nu}_2 \in \abinom{n}{k}$ with $\overline{\nu}_1 \neq \overline{\nu}_2$, we have 
$\zero(P^{\overline{\nu}_1, \Omega}_{\lambda \mu}) \cap \zero(P^{\overline{\nu}_2, \Omega}_{\lambda \mu}) \neq \emptyset$. This means that $\zero(P^{\emptyset, \Omega}_{\lambda \mu}) = \bigcup_{\overline{\nu} \in  \abinom{n}{k}} \zero(P^{\overline{\nu}, \Omega}_{\lambda \mu})$, and indeed $\{\zero(P^{\overline{\nu}, \Omega}_{\lambda \mu}): \overline{\nu} \in  \abinom{n}{k}\}$ forms a partition of $\zero(P^{\emptyset, \Omega}_{\lambda \mu})$. Take $\Omega = \Omega_0$ for the Knutson-Tao-Woodward puzzles for example, then $\#\zero(P^{\overline{\nu}, \Omega_0}_{\lambda \mu})$ is equal to the Littlewood-Richardson coefficient $c^{\overline{\nu}}_{\lambda \mu}$, and thus $\#\zero(P^{\emptyset, \Omega_0}_{\lambda \mu}) = \sum_{\overline{\nu} \in  \abinom{n}{k}} \#\zero(P^{\overline{\nu}, \Omega_0}_{\lambda \mu}) = \sum_{\overline{\nu} \in \abinom{n}{k}} c^{\overline{\nu}}_{\lambda \mu}$. Note that when $|\overline{\nu}| \neq |\lambda| + |\mu|$, we have $c^{\overline{\nu}}_{\lambda \mu} = 0$ and thus $\#\zero(P^{\overline{\nu}, \Omega_0}_{\lambda \mu}) = 0$. This means that we can safely replace $\abinom{n}{k}$ in the summand by $\fnum_2^n$. 

Like the puzzle ideal $I^{\nu, \Omega_0}_{\lambda \mu}$ whose varieties one-one correspond to $\zero(P^{\nu, \Omega_0}_{\lambda \mu})$, for the side-free puzzles we can define the side-free puzzle ideals as follows. Take the $\nu$-free puzzle $\zero(P^{\emptyset, \Omega}_{\lambda \mu})$ for example, we only need to change the defining polynomials in $\pset{F}_1$ and $\pset{F}_2$ for the puzzle ideals $I^{\nu, \Omega}_{\lambda \mu}$ to loose the restrictions on the values of the bottom side of $\bigtriangleup^{\emptyset}_{\lambda \mu}$. Let $x_{b1}, \ldots, x_{bn}$ be the variables assigned to the $n$ intervals of the bottom side of $\bigtriangleup^{\emptyset}_{\lambda \mu}$. Then in $\pset{F}_1$ we replace the field equations $\{x_{bj}^3 - x_{bj}:\,j=1, \ldots, n\}$ for $\fnum_3$ by $\{x_{bj}^2 - x_{bj}:\,j=1, \ldots, n\}$ for $\fnum_2$ to have a new set $\tilde{\pset{F}}_1$ and in $\pset{F}_2$ we remove the polynomials $\{x_{bj} - \nu_j :\,j=1, \ldots, n\}$ to have a new set $\tilde{\pset{F}}_2$.

\begin{definition}\rm
Let $P^{\emptyset, \Omega}_{\lambda \mu}$ be a $\nu$-free puzzle with $\Omega = \Omega_t \cup \Omega_r \cup \Omega_b$, $\overline{\Omega}$ be a separable atomic refinement of $\Omega$ induced by $\{\overline{t}_p:\, p \in \Omega\}$, and $\tilde{\pset{F}}_1, \tilde{\pset{F}}_2, \pset{F}_3, \ldots, \pset{F}_6$ be as defined above for fixed choices of distinguishing, forbidding, and implying polynomials. Then the ideal in $\fnum_3[\p{x}]$ generated by $\tilde{\pset{F}}_1 \cup  \tilde{\pset{F}}_2 \cup \pset{F}_3 \cup \cdots \cup \pset{F}_6$ is called the \emph{$\nu$-free puzzle ideal} of $P^{\emptyset, \Omega}_{\lambda \mu}$ and denoted by $I^{\emptyset, \Omega}_{\lambda \mu}$.
\end{definition}

The following theorem can be proved in almost the same way as the proof of Theorem~\ref{thm:main}, and the one-one correspondence is also the same as that in Theorem~\ref{thm:main}. 

\begin{theorem}\label{thm:side-free}
    Let $P^{\emptyset, \Omega}_{\lambda \mu}$ be a $\nu$-free puzzle with $\Omega = \Omega_t \cup \Omega_r \cup \Omega_b$, $\overline{\Omega}$ be a separable atomic refinement of $\Omega$, and $I^{\emptyset, \Omega}_{\lambda \mu}$ be the puzzle ideal of $P^{\nu, \Omega}_{\lambda \mu}$. Then there exists a one-one correspondence between $\var(I^{\emptyset, \Omega}_{\lambda \mu}) = \zero(P^{\emptyset, \Omega}_{\lambda \mu})$. 
\end{theorem}

Take the Littlewood-Richardson coefficients in Knutson-Tao-Woodward puzzles for example. Let $\p{x}_{b}$ be the variables in $\p{x}$ corresponding to the $n$ intervals of the bottom side of $\bigtriangleup^{\emptyset}_{\lambda \mu}$. Recall that $\#\zero(P^{\emptyset, \Omega_0}_{\lambda \mu}) = \sum_{\overline{\nu} \in \fnum_2^n} c^{\overline{\nu}}_{\lambda \mu}$, and thus to recover the coefficients $c^{\overline{\nu}}_{\lambda \mu}$ for all possible $\overline{\nu} \in \fnum_2^n$, we can compute the variety $\var(I^{\emptyset, \Omega_0}_{\lambda \mu})$ by solving the corresponding polynomial system and then partition all the solutions according to the values restricted to $\p{x}_b$. 

Let $I^{\emptyset, \Omega}_{\lambda \mu} \subseteq \fnum_{3}[\p{x}]$ be a $\nu$-free puzzle ideal and $\p{x}_{b}$ be the set of variables corresponding to  the $n$ intervals of the bottom side of $\bigtriangleup^{\emptyset}_{\lambda \mu}$. Then the elimination ideal $I^{\emptyset, \Omega_0}_{\lambda \mu} \cap \fnum_3[\p{x}_{b}]$ is denoted by $I^{\emptyset, \Omega_0}$, and it contains all the information we need for the values $\overline{\nu} \in \fnum_2^n$ for the coefficients $c^{\overline{\nu}}_{\lambda \mu}$, as shown in the proposition below. 

\begin{proposition}
    Let $P^{\emptyset, \Omega}_{\lambda \mu}$ be a $\nu$-free puzzle with $\Omega = \Omega_t \cup \Omega_r \cup \Omega_b$, $\overline{\Omega}$ be a separable atomic refinement of $\Omega$, $I^{\emptyset, \Omega}_{\lambda \mu}$ be the $\nu$-free puzzle ideal of $P^{\emptyset, \Omega}_{\lambda \mu}$, and $I^{\emptyset, \Omega}$ be the elimination ideal of $I^{\emptyset, \Omega}_{\lambda \mu}$ w.r.t. $\p{x}_b$. Then the following statements hold. 
    \begin{enumerate}
        \item[(1)] $I^{\emptyset, \Omega}$ is radical.
        \item[(2)] Let $I^{\emptyset, \Omega} = \cap_{k}p_k$ be the prime decomposition of $I^{\emptyset, \Omega}$. Then each $p_k$ is of the form $\bases{x_{b1} - \overline{\nu}_{k1}, \ldots, x_{bn} - \overline{\nu}_{kn}}$, where $\overline{\nu}_{ki} \in \fnum_2$ for $i=1, \ldots, n$.
        \item[(3)] Let $\overline{\nu}_k = (\overline{\nu}_{k1}, \ldots, \overline{\nu}_{kn}) \in \var(I^{\emptyset, \Omega})$ and $p_k$ be the prime ideal corresponding to it. Then $I^{\overline{\nu}_k, \Omega}_{\lambda \mu} = I^{\emptyset, \Omega}_{\lambda \mu} + p_k$.
    \end{enumerate} 
\end{proposition}

\begin{proof}
(1) With $\{x_{bi}^2 - x_{bi}:\,i=1, \ldots, n\} \subseteq I^{\emptyset, \Omega_0}$, clearly $I^{\emptyset, \Omega_0}$ is radical.

(2) Let $J$ be the ideal generated by the field polynomials $\{x_{bi}^2 - x_{bi}:\,i=1, \ldots, n\} $. Then consider the quotient ideal $ I^{\emptyset, \Omega_0} / J \subseteq \fnum_3[\p{x}_{b}] / J$. The quotient ring $\fnum_3[\p{x}_{b}] / J$ is a finite commutative ring so that every prime ideal is maximal. Therefore, the prime ideal in prime decomposition of $I^{\emptyset, \Omega_0} / J$ is of the form $\bases{x_{b1} - \overline{\nu}_{k1}, \ldots, x_{bn} - \overline{\nu}_{kn}}$, where $\overline{\nu}_{ki} \in \fnum_2$ for $i=1, \ldots, n$. The same also holds for the ideal $ I^{\emptyset, \Omega_0} $.

(3) It follows directly from the definitions of puzzle ideals and side-free puzzle ideals.
\end{proof}

In the remaining part of this section we present a non-trivial example of side-free puzzle ideals and demonstrate how to obtain all the information of the Littlewood-Richardson coefficients for given $\lambda$ and $\mu$ with computation of \grobner bases for side-free puzzle ideals. 

Consider $\lambda = (1, 0, 1, 0, 1, 0, 1, 0, 1, 0, 1, 0, 1, 0, 1, 0)$ and $\mu = (0, 0, 1, 0, 1, 0, 1, 0, 1, 0, 1, 0, 1, 1, 0, 1) \in \abinom{16}{8}$ which correspond to the partitions $(8, 7, 6, 5, 4, 3, 2, 1)$ and $(6, 5, 4, 3, 2, 1, 1)$ respectively. Now we construct the side-free puzzle ideal $I^{\emptyset, \Omega_0}_{\lambda \mu} $ for the Knutson-Tao-Woodward puzzle $P^{\emptyset, \Omega_0}_{\lambda \mu}$ and compute all the Littlewood-Richardson coefficients for $\lambda$ and $\mu$. 

Inside $\bigtriangleup^{\emptyset}_{\lambda \mu}$ there are $256$ unit triangles with $408$ unit intervals. The defining polynomials of $I^{\emptyset, \Omega_0}_{\lambda \mu}$ can be constructed in the four groups discussed in Section~\ref{sec:atomic}. Fix a variable order in which the variables $\p{x}_{b}$ corresponding to $\nu$ are ordered smaller than the others, then the lexicographic Gr\"{o}bner basis $\pset{G}$ of $I^{\emptyset, \Omega_0}_{\lambda \mu}$ contains $544$ polynomials, among which $232$ are linear, $222$ are quadratic, and $90$ are cubic. For the elimination ideal $I^{\emptyset, \Omega_0}$, its lexicographic Gr\"{o}bner basis 
$\pset{G} \cap \fnum_3[\p{x}_b]$ contains $61$ polynomials, and the prime decomposition of $I^{\emptyset, \Omega_0}$ tells us that the corresponding variety $\var(I^{\emptyset, \Omega_0})$ has $11$ points, which means that there are $11$ non-zero Littlewood-Richardson coefficients in the expansion $ s_{\lambda} s_{\mu} = \sum_{\nu} c_{\lambda,\mu}^{\nu} s_{\nu}$ for $\nu \in \abinom{16}{8}$. These $11$ points are listed below. 

\begin{equation*}
    \begin{split}
    &\var(I^{\emptyset, \Omega_0}) = \\
    \{&(1, 1, 1, 1, 1, 1, 0, 1, 0, 0, 0, 0, 1, 0, 0, 0), (1, 1, 0, 1, 1, 1, 1, 1, 1, 0, 0, 0, 0, 0, 0, 0) \\
    &(1, 1, 1, 0, 1, 1, 1, 1, 0, 1, 0, 0, 0, 0, 0, 0), (1, 1, 1, 1, 1, 0, 1, 1, 0, 0, 0, 1, 0, 0, 0, 0) \\
    &(1, 1, 1, 1, 1, 1, 1, 0, 0, 0, 0, 0, 0, 1, 0, 0), (1, 1, 1, 1, 1, 1, 0, 0, 1, 0, 0, 1, 0, 0, 0, 0) \\
    &(1, 1, 1, 1, 0, 1, 1, 1, 0, 0, 1, 0, 0, 0, 0, 0), (1, 1, 1, 1, 1, 0, 1, 0, 1, 0, 1, 0, 0, 0, 0, 0) \\
    &(1, 1, 1, 1, 1, 1, 0, 0, 0, 1, 1, 0, 0, 0, 0, 0), (1, 1, 1, 1, 0, 1, 1, 0, 1, 1, 0, 0, 0, 0, 0, 0) \\
    &(1, 1, 1, 1, 1, 0, 0, 1, 1, 1, 0, 0, 0, 0, 0, 0)\}.
    \end{split}
\end{equation*}

Take one point $ \overline{\nu} = (1, 1, 1, 1, 1, 1, 0, 1, 0, 0, 0, 0, 1, 0, 0, 0) \in \var(I^{\emptyset, \Omega_0})$ for example and let $p_{\overline{\nu}} $ be the prime ideal in the prime decomposition of $I^{\emptyset, \Omega_0}$ which corresponds to $\overline{\nu}$. By computing the Gr\"{o}bner basis of $I^{\overline{\nu}, \Omega_0}_{\lambda \mu} = I^{\emptyset, \Omega_0}_{\lambda \mu} + p_{\overline{\nu}}$ we know that the variety $\var(I^{\overline{\nu}, \Omega_0}_{\lambda \mu})$ contains $5$ points, which implies that $c^{\overline{\nu}}_{\lambda \mu} = 5$. One out of the $5$ tilings in $P^{\overline{\nu}, \Omega_0}_{\lambda \mu}$ is shown in Figure~\ref{fig:bigpuzzle} below. 

\begin{figure}[H]
    \begin{center}
        \begin{tikzpicture}[scale=0.8]
            \foreach \x in {0,...,15} {
            \draw (0+\x/2,0.866*\x) --(16-\x/2,0.866*\x);
            }
            \foreach \x in {0,...,15} {
            \draw (\x,0) -- (8 +0.5*\x,0.866*16 - 0.866*\x);
            }
            \foreach \x in {0,...,15} {
            \draw (\x*0.5+0.5,0.866*\x+0.866) -- (1+\x,0);
            }

            \drawuptriangle{lightblue}{0}{0}{1}{1}{1};
            \drawdowntriangle{lightblue}{1}{0}{1}{1}{1};
            \drawuptriangle{lightblue}{1}{0}{1}{1}{1};
            \drawdowntriangle{lightblue}{2}{0}{1}{1}{1};
            \drawuptriangle{lightblue}{2}{0}{1}{1}{1};
            \drawdowntriangle{lightblue}{3}{0}{1}{1}{1};
            \drawuptriangle{lightblue}{3}{0}{1}{1}{1};
            \drawdowntriangle{lightblue}{4}{0}{1}{1}{1};
            \drawuptriangle{lightblue}{4}{0}{1}{1}{1};
            \drawdowntriangle{lightblue}{5}{0}{1}{1}{1};
            \drawuptriangle{lightblue}{5}{0}{1}{1}{1};
            \drawNErumbus{lightgreen}{6}{0};
            \drawdowntriangle{lightblue}{7}{0}{1}{1}{1};
            \drawuptriangle{lightblue}{7}{0}{1}{1}{1};
            \drawNErumbus{lightgreen}{8}{0};
            
            \drawNErumbus{lightgreen}{9}{0};
            \drawNErumbus{lightgreen}{10}{0};
            \drawNErumbus{lightgreen}{11}{0};
            \drawdowntriangle{lightblue}{12}{0}{1}{1}{1};
            \drawuptriangle{lightblue}{12}{0}{1}{1}{1};
            \drawNErumbus{lightgreen}{13}{0};
            \drawNErumbus{lightgreen}{14}{0};
            \drawNErumbus{lightgreen}{15}{0};

            \drawSWrumbus{lightgreen}{0.5}{0.866};
            \drawSWrumbus{lightgreen}{0.5+1}{0.866};
            \drawSWrumbus{lightgreen}{0.5+2}{0.866};
            \drawSWrumbus{lightgreen}{0.5+3}{0.866};
            \drawSWrumbus{lightgreen}{0.5+4}{0.866};
            \drawuptriangle{lightred}{0.5+5}{0.866}{0}{0}{0};
            \drawdowntriangle{lightred}{0.5+6}{0.866}{0}{0}{0};
            \drawSWrumbus{lightgreen}{0.5+6}{0.866};
            \drawuptriangle{lightred}{0.5+7}{0.866}{0}{0}{0};
            \drawdowntriangle{lightred}{0.5+8}{0.866}{0}{0}{0};
            \drawuptriangle{lightred}{0.5+8}{0.866}{0}{0}{0};
            \drawdowntriangle{lightred}{0.5+9}{0.866}{0}{0}{0};
            \drawuptriangle{lightred}{0.5+9}{0.866}{0}{0}{0};
            \drawdowntriangle{lightred}{0.5+10}{0.866}{0}{0}{0};
            \drawuptriangle{lightred}{0.5+10}{0.866}{0}{0}{0};
            \drawdowntriangle{lightred}{0.5+11}{0.866}{0}{0}{0};
            \drawSWrumbus{lightgreen}{0.5+11}{0.866};
            \drawuptriangle{lightred}{0.5+12}{0.866}{0}{0}{0};
            \drawdowntriangle{lightred}{0.5+13}{0.866}{0}{0}{0};
            \drawuptriangle{lightred}{0.5+13}{0.866}{0}{0}{0};
            \drawdowntriangle{lightred}{0.5+14}{0.866}{0}{0}{0};
            \drawuptriangle{lightred}{0.5+14}{0.866}{0}{0}{0};

            \drawuptriangle{lightblue}{1+0}{0.866*2}{1}{1}{1};
            \drawdowntriangle{lightblue}{1+1}{0.866*2}{1}{1}{1};
            \drawuptriangle{lightblue}{1+1}{0.866*2}{1}{1}{1};
            \drawdowntriangle{lightblue}{1+2}{0.866*2}{1}{1}{1};
            \drawuptriangle{lightblue}{1+2}{0.866*2}{1}{1}{1};
            \drawdowntriangle{lightblue}{1+3}{0.866*2}{1}{1}{1};
            \drawuptriangle{lightblue}{1+3}{0.866*2}{1}{1}{1};
            \drawdowntriangle{lightblue}{1+4}{0.866*2}{1}{1}{1};
            \drawuptriangle{lightblue}{1+4}{0.866*2}{1}{1}{1};
            \drawNErumbus{lightgreen}{1+5}{0.866*2};
            \drawdowntriangle{lightblue}{1+6}{0.866*2}{1}{1}{1};
            \drawuptriangle{lightblue}{1+6}{0.866*2}{1}{1}{1};
            \drawNErumbus{lightgreen}{1+7}{0.866*2};
            \drawNErumbus{lightgreen}{1+8}{0.866*2};
            \drawNErumbus{lightgreen}{1+9}{0.866*2};
            \drawNErumbus{lightgreen}{1+10}{0.866*2};
            \drawdowntriangle{lightblue}{1+11}{0.866*2}{1}{1}{1};
            \drawuptriangle{lightblue}{1+11}{0.866*2}{1}{1}{1};
            \drawNErumbus{lightgreen}{1+12}{0.866*2};
            \drawNErumbus{lightgreen}{1+13}{0.866*2};

            \drawSWrumbus{lightgreen}{1.5+0}{0.866*3};
            \drawSWrumbus{lightgreen}{1.5+1}{0.866*3};
            \drawSWrumbus{lightgreen}{1.5+2}{0.866*3};
            \drawSWrumbus{lightgreen}{1.5+3}{0.866*3};
            \drawuptriangle{lightred}{1.5+4}{0.866*3}{0}{0}{0};
            \drawuprumbus{lightgreen}{1.5+5}{0.866*3};
            \drawuptriangle{lightblue}{1.5+5}{0.866*3}{1}{1}{1};
            \drawNErumbus{lightgreen}{1.5+6}{0.866*3};
            \drawNErumbus{lightgreen}{1.5+7}{0.866*3};
            \drawNErumbus{lightgreen}{1.5+8}{0.866*3};
            \drawNErumbus{lightgreen}{1.5+9}{0.866*3};
            \drawdowntriangle{lightblue}{1.5+10}{0.866*3}{1}{1}{1};
            \drawuptriangle{lightblue}{1.5+10}{0.866*3}{1}{1}{1};
            \drawNErumbus{lightgreen}{1.5+11}{0.866*3};
            \drawNErumbus{lightgreen}{1.5+12}{0.866*3};

            \drawuptriangle{lightblue}{2+0}{0.866*4}{1}{1}{1};
            \drawdowntriangle{lightblue}{2+1}{0.866*4}{1}{1}{1};
            \drawuptriangle{lightblue}{2+1}{0.866*4}{1}{1}{1};
            \drawdowntriangle{lightblue}{2+2}{0.866*4}{1}{1}{1};
            \drawuptriangle{lightblue}{2+2}{0.866*4}{1}{1}{1};
            \drawdowntriangle{lightblue}{2+3}{0.866*4}{1}{1}{1};
            \drawuptriangle{lightblue}{2+3}{0.866*4}{1}{1}{1};
            \drawdowntriangle{lightblue}{2+4}{0.866*4}{1}{1}{1};
            \drawdowntriangle{lightred}{2+5}{0.866*4}{0}{0}{0};
            \drawuptriangle{lightred}{2+5}{0.866*4}{0}{0}{0};
            \drawdowntriangle{lightred}{2+6}{0.866*4}{0}{0}{0};
            \drawuptriangle{lightred}{2+6}{0.866*4}{0}{0}{0};
            \drawdowntriangle{lightred}{2+7}{0.866*4}{0}{0}{0};
            \drawuptriangle{lightred}{2+7}{0.866*4}{0}{0}{0};
            \drawdowntriangle{lightred}{2+8}{0.866*4}{0}{0}{0};
            \drawuptriangle{lightred}{2+8}{0.866*4}{0}{0}{0};
            \drawdowntriangle{lightred}{2+9}{0.866*4}{0}{0}{0};
            \drawSWrumbus{lightgreen}{2+9}{0.866*4};
            \drawuptriangle{lightred}{2+10}{0.866*4}{0}{0}{0};
            \drawdowntriangle{lightred}{2+11}{0.866*4}{0}{0}{0};
            \drawuptriangle{lightred}{2+11}{0.866*4}{0}{0}{0};

            \drawSWrumbus{lightgreen}{2.5+0}{0.866*5};
            \drawSWrumbus{lightgreen}{2.5+1}{0.866*5};
            \drawSWrumbus{lightgreen}{2.5+2}{0.866*5};
            \drawSWrumbus{lightgreen}{2.5+3}{0.866*5};
            \drawuptriangle{lightred}{2.5+4}{0.866*5}{0}{0}{0};
            \drawdowntriangle{lightred}{2.5+5}{0.866*5}{0}{0}{0};
            \drawuptriangle{lightred}{2.5+5}{0.866*5}{0}{0}{0};
            \drawdowntriangle{lightred}{2.5+6}{0.866*5}{0}{0}{0};
            \drawuptriangle{lightred}{2.5+6}{0.866*5}{0}{0}{0};
            \drawdowntriangle{lightred}{2.5+7}{0.866*5}{0}{0}{0};
            \drawuptriangle{lightred}{2.5+7}{0.866*5}{0}{0}{0};
            \drawdowntriangle{lightred}{2.5+8}{0.866*5}{0}{0}{0};
            \drawuptriangle{lightred}{2.5+8}{0.866*5}{0}{0}{0};
            \drawuprumbus{lightgreen}{2.5+9}{0.866*5};
            \drawuptriangle{lightblue}{2.5+9}{0.866*5}{1}{1}{1};
            \drawNErumbus{lightgreen}{2.5+10}{0.866*5};

            \drawuptriangle{lightblue}{3+0}{0.866*6}{1}{1}{1};
            \drawdowntriangle{lightblue}{3+1}{0.866*6}{1}{1}{1};
            \drawuptriangle{lightblue}{3+1}{0.866*6}{1}{1}{1};
            \drawdowntriangle{lightblue}{3+2}{0.866*6}{1}{1}{1};
            \drawuptriangle{lightblue}{3+2}{0.866*6}{1}{1}{1};
            \drawdowntriangle{lightblue}{3+3}{0.866*6}{1}{1}{1};
            \drawuptriangle{lightblue}{3+3}{0.866*6}{1}{1}{1};
            \drawNErumbus{lightgreen}{3+4}{0.866*6};
            \drawNErumbus{lightgreen}{3+5}{0.866*6};
            \drawNErumbus{lightgreen}{3+6}{0.866*6};
            \drawNErumbus{lightgreen}{3+7}{0.866*6};
            \drawdowntriangle{lightblue}{3+8}{0.866*6}{1}{1}{1};
            \drawdowntriangle{lightred}{3+9}{0.866*6}{0}{0}{0};
            \drawuptriangle{lightred}{3+9}{0.866*6}{0}{0}{0};

            \drawSWrumbus{lightgreen}{3.5+0}{0.866*7};
            \drawSWrumbus{lightgreen}{3.5+1}{0.866*7};
            \drawSWrumbus{lightgreen}{3.5+2}{0.866*7};
            \drawuptriangle{lightred}{3.5+3}{0.866*7}{0}{0}{0};
            \drawdowntriangle{lightred}{3.5+4}{0.866*7}{0}{0}{0};
            \drawuptriangle{lightred}{3.5+4}{0.866*7}{0}{0}{0};
            \drawdowntriangle{lightred}{3.5+5}{0.866*7}{0}{0}{0};
            \drawuptriangle{lightred}{3.5+5}{0.866*7}{0}{0}{0};
            \drawdowntriangle{lightred}{3.5+6}{0.866*7}{0}{0}{0};
            \drawuptriangle{lightred}{3.5+6}{0.866*7}{0}{0}{0};
            \drawuprumbus{lightgreen}{3.5+7}{0.866*7};
            \drawuptriangle{lightblue}{3.5+7}{0.866*7}{1}{1}{1};
            \drawNErumbus{lightgreen}{3.5+8}{0.866*7};

            \drawuptriangle{lightblue}{4+0}{0.866*8}{1}{1}{1};
            \drawdowntriangle{lightblue}{4+1}{0.866*8}{1}{1}{1};
            \drawuptriangle{lightblue}{4+1}{0.866*8}{1}{1}{1};
            \drawdowntriangle{lightblue}{4+2}{0.866*8}{1}{1}{1};
            \drawuptriangle{lightblue}{4+2}{0.866*8}{1}{1}{1};
            \drawNErumbus{lightgreen}{4+3}{0.866*8};
            \drawNErumbus{lightgreen}{4+4}{0.866*8};
            \drawNErumbus{lightgreen}{4+5}{0.866*8};
            \drawdowntriangle{lightblue}{4+6}{0.866*8}{1}{1}{1};
            \drawdowntriangle{lightred}{4+7}{0.866*8}{0}{0}{0};
            \drawuptriangle{lightred}{4+7}{0.866*8}{0}{0}{0};

            \drawSWrumbus{lightgreen}{4.5+0}{0.866*9};
            \drawSWrumbus{lightgreen}{4.5+1}{0.866*9};
            \drawuptriangle{lightred}{4.5+2}{0.866*9}{0}{0}{0};
            \drawdowntriangle{lightred}{4.5+3}{0.866*9}{0}{0}{0};
            \drawuptriangle{lightred}{4.5+3}{0.866*9}{0}{0}{0};
            \drawdowntriangle{lightred}{4.5+4}{0.866*9}{0}{0}{0};
            \drawuptriangle{lightred}{4.5+4}{0.866*9}{0}{0}{0};
            \drawuprumbus{lightgreen}{4.5+5}{0.866*9};
            \drawuptriangle{lightblue}{4.5+5}{0.866*9}{1}{1}{1};
            \drawNErumbus{lightgreen}{4.5+6}{0.866*9};

            \drawuptriangle{lightblue}{5+0}{0.866*10}{1}{1}{1};
            \drawdowntriangle{lightblue}{5+1}{0.866*10}{1}{1}{1};
            \drawuptriangle{lightblue}{5+1}{0.866*10}{1}{1}{1};
            \drawNErumbus{lightgreen}{5+2}{8.66};
            \drawNErumbus{lightgreen}{5+3}{8.66};
            \drawdowntriangle{lightblue}{5+4}{0.866*10}{1}{1}{1};
            \drawdowntriangle{lightred}{5+5}{8.66}{0}{0}{0};
            \drawuptriangle{lightred}{5+5}{8.66}{0}{0}{0};

            \drawSWrumbus{lightgreen}{5.5+0}{0.866*11};
            \drawuptriangle{lightred}{5.5+1}{0.866*11}{0}{0}{0};
            \drawdowntriangle{lightred}{5.5+2}{0.866*11}{0}{0}{0};
            \drawuptriangle{lightred}{5.5+2}{0.866*11}{0}{0}{0};
            \drawuprumbus{lightgreen}{5.5+3}{0.866*11};
            \drawuptriangle{lightblue}{5.5+3}{0.866*11}{1}{1}{1};
            \drawNErumbus{lightgreen}{5.5+4}{0.866*11};

            \drawuptriangle{lightblue}{6+0}{0.866*12}{1}{1}{1};
            \drawNErumbus{lightgreen}{6+1}{0.866*12};
            \drawdowntriangle{lightblue}{6+2}{0.866*12}{1}{1}{1};
            \drawdowntriangle{lightred}{6+3}{0.866*12}{0}{0}{0};
            \drawuptriangle{lightred}{6+3}{0.866*12}{0}{0}{0};

            \drawuptriangle{lightred}{6.5+0}{0.866*13}{0}{0}{0};
            \drawuprumbus{lightgreen}{6.5+1}{0.866*13};
            \drawuptriangle{lightblue}{6.5+1}{0.866*13}{1}{1}{1};
            \drawNErumbus{lightgreen}{6.5+2}{0.866*13};

            \drawdowntriangle{lightred}{7+1}{0.866*14}{0}{0}{0};
            \drawuptriangle{lightred}{7+1}{0.866*14}{0}{0}{0};

            \drawuptriangle{lightred}{7.5}{0.866*15}{0}{0}{0};
        \end{tikzpicture}
    \end{center}
    \caption{One tiling of a puzzle of size 16.}
    \label{fig:bigpuzzle}
\end{figure}
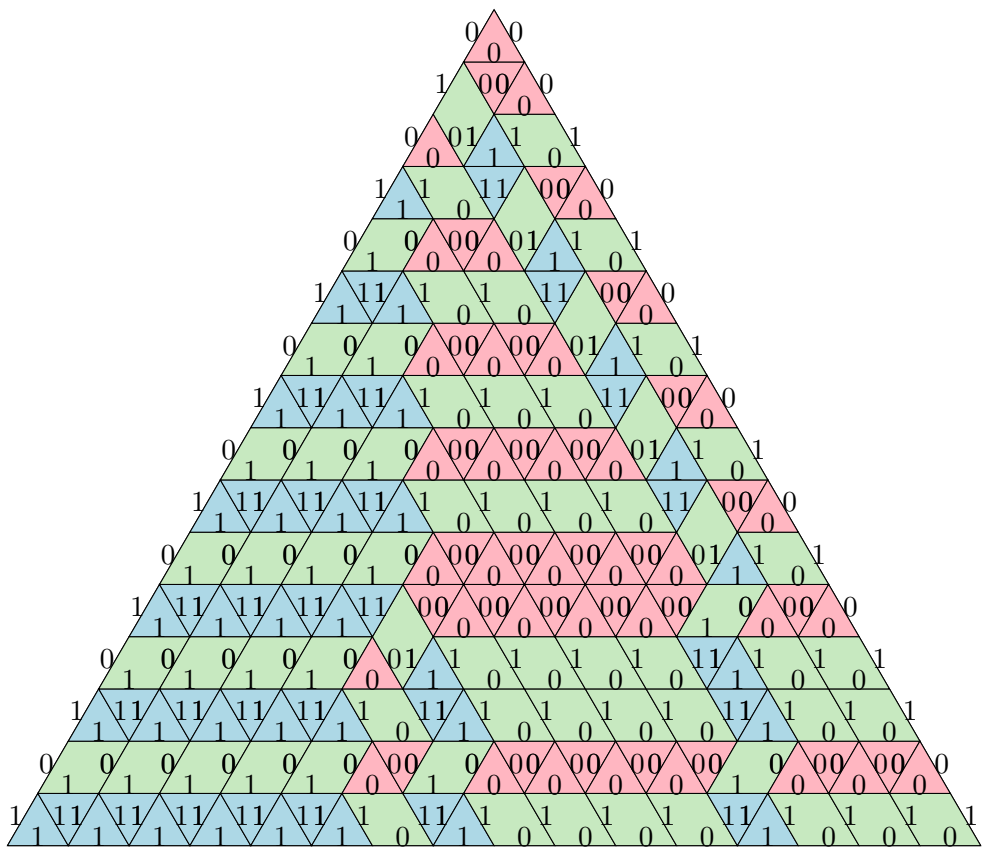

\bibliographystyle{amsplain}
\bibliography{references.bib}

\end{document}